\documentclass[a4paper, 11pt]{article}
\usepackage{amssymb,amsmath,graphicx}
\usepackage{amsfonts}
\usepackage{subfig}
\usepackage{epsfig,latexsym,graphicx}
\usepackage{setspace}
\usepackage[svgnames]{xcolor}
\usepackage{listings} 
\RequirePackage[OT1]{fontenc}
\RequirePackage{amsthm,amsmath}
\RequirePackage[numbers,sort,compress]{natbib}

\newtheorem{theorem}{Theorem}[section]

\newtheorem{lemma}[theorem]{Lemma}
\newtheorem{proposition}[theorem]{Proposition}

\newtheorem{remark}{Remark}[section]
\newtheorem{definition}{Definition}[section]

\theoremstyle{plain}

\begin{document}
	
\title{\bf{General Robust Bayes Pseudo-Posteriors: 
		Exponential Convergence Results with Applications}}
\author{
	Abhik Ghosh$^\dagger$, Tuhin Majumder$^{\ddagger}$ and Ayanendranath Basu$^{\dagger\ast}$\\
$^\dagger$ Indian Statistical Institute, India. 
	$^\ddagger$ North Carolina State University, USA.\\
$^\ast$Corresponding author: ayanbasu@isical.ac.in
}
\date{}
\maketitle


\begin{abstract}
Although Bayesian inference is an immensely popular paradigm among a large segment of scientists including statisticians,
most  applications consider objective priors and need critical investigations \citep{Efron:2013}. 
While it has several optimal properties, a major drawback of Bayesian inference is 
the lack of robustness against data contamination and model misspecification, 
which becomes pernicious in the use of objective priors. 
This paper presents the general formulation of a Bayes pseudo-posterior distribution yielding robust inference. 
Exponential convergence results related to the new pseudo-posterior and the corresponding Bayes estimators 
are established under the general parametric set-up and 
illustrations are provided for the independent stationary as well as 
non-homogeneous models. 
Several additional details and properties of the procedure are described, 
including 
the estimation under fixed-design regression models. 

\end{abstract}


\section{Introduction}
\label{SEC:intro}

Bayesian analysis is arguably one of the most popular statistical paradigms with applications 
across different scientific disciplines.
It is widely preferred by many non-statisticians due to its nice interpretability 
and incorporation of prior knowledge. 
From a statistical  point of view, it is widely accepted
even among many non-Bayesians, because of its nice optimal (asymptotic) properties. 
Bayesian inference is built on the famous `Bayes theorem', the celebrated 1763 paper of Thomas Bayes,
which combines prior knowledge with experimental evidence to produce posterior conclusions.  
However, over these 250 years, Bayesian inference has also been subject to several criticisms 
and some of these debates are still ongoing \citep{Efron:2013}.
Other than the controversies about its internal logic \citep{Halpern:1999a,Dupre/Tipler:2009},
a major practical drawback of Bayesian inference is 
its non-robust nature against misspecification in models (including data contamination and outliers) and priors,
as has been extensively  observed in the literature; 
see \cite{Berk:1966,Wasserman:1996,Millar/Stewart:2007,DeBlasi/Walker:2012,
Owhadi/etc:2015a} and the references therein.
The optimal solution to this problem has been developed mainly for prior misspecifications  
\citep{Berger:1994,Berger/Berliner:1986, gd91,db94,Delampady/Dey:1994,gw95,gds06};
the Bayesians traditionally viewed the model to be perfect for the given data.
Thus the possibility of model misspecification and data contamination has been generally ignored for a long time 
till the appearance of some very recent publications, some of which we describe later in this section.

In applying Bayesian inference to the complicated datasets of the present era, 
we need to use complex and sophisticated models which are highly prone to misspecification or data contamination.
In reality, where ``All models are wrong", 
the Bayesian philosophy of refining the fixed model adaptively \citep{Gelman/etc:1996} often fails 
to handle complex scenarios or leads to ``a model as complex as the data'' \citep{Wang/Blei:2016}.
Data contamination can lead to erroneous posterior conclusions. 
The problem becomes more clear but pernicious in case of inference with objective or reference priors.
For example, the Bayes estimate of the mean of a normal model, 
with any objective prior and symmetric loss function, is the highly non-robust sample mean. 
What is a matter of greater concern, as noted by \citet{Efron:2013}, 
is that most of the recent applications of Bayesian inference 
hinge on objective priors and so they always need to be scrutinized carefully, 
sometimes even from a frequentist perspective.
The posterior non-robustness against model misspecification and data contamination makes the process vulnerable 
and we clearly need a solution to this problem. 

From a true Bayesian perspective, there are only few solutions to the problem of model misspecification
\citep{Ritov:1985,Ritov:1987,Sivaganesan:1993,Shyamalkumar:2000}.
However, most of them, if not all, assume that the perturbation in the model is known beforehand,
such as gross error contaminated models with known contamination proportion $\epsilon$.
For modern complex datasets, this is rarely meaningful.
There has been several recent publications which are motivated by the need to safeguard Bayes inference against model misspecification 
by relying on a generalized (pseudo) posterior 
which is expressed in terms of a loss function and a tuning parameter $\eta$
\cite{Alquier/Lounici:2011,Catoni:2007,Jiang/Tanner:2008, Walker/Hjort:2001,Kleijn/VanderVaart:2006,Gruenwald/vanOmmen:2017,Holmes/Walker:2017,Ramamoorthi/etc:2015,DeBlasi/Walker:2012}.
This approach, referred to as the PAC-Bayesian approach generated from Gibb's posterior, 
has been quite successful in  regression and other supervised classification problems with misspecified model assumptions. 
But the resulting inference is not robust against outliers with respect 
to a specified model which is correct for the majority of the data. 
This is because every sample observation, including outliers, receives equal weight in PAC-Bayesian approach
and hence it closely resembles the model robust non-parametric analysis;
see \cite{Ghosh/Basu:2016}.

To achieve robustness against data contamination (outliers) in Bayesian inference, 
some attempts have been made to develop alternative solutions by linking Bayesian inference suitably to the frequentist concept of robustness.
In the frequentist sense, there are two major approaches to achieve robustness,
namely the use of heavy tailed distributions (e.g., $t$-distribution in place of normal),
or new (robust) inference methodologies \citep{Hampel/etc:1986,Basu/etc:2011}.
The first one has been adopted by some Bayesian scientists;
see \cite{Andrade/OHagan:2006,Andrade/OHagan:2011} and \cite{Desgagne:2013} among others.
However, the difficulty with this approach is the availability of appropriate heavy tailed alternatives in complex scenarios
and it indeed does not solve the non-robustness of Bayesian inference for a specified model (which might be of a lighter tail).
The second approach of frequentist robustness
serves the purpose but differs in the strictest probabilistic sense from  the Bayesian philosophy,
since one needs to alter the posterior density appropriately to achieve robustness 
against data contamination or model misspecification;
the resulting modified posteriors are generally referred to as pseudo-posterior densities.
Different such pseudo-posteriors have been proposed by 
\cite{Greco/etc:2008,Agostinelli/Greco:2013,hv11, Ghosh/Basu:2016, Danesi/etc:2016,Atkinson/etc:2017,Nakagawa/Hashimoto:2017}; 
but all of them have primarily considered independent stationary models and have different pros and cons.
Another recent attempt, in the borderline of these two approaches, has been proposed by \cite{Wang/Blei:2016},
who have transformed the given model to a localized model involving hyperparameters 
to be estimated through the empirical Bayes approach.

\subsection{Background: $R^{(\alpha)}$-posterior for IID set-up}

We consider a particular pseudo-posterior originally proposed by \cite{Ghosh/Basu:2016}
in the independently and identically distributed (IID) set-up.
This choice has been motivated by its several nice properties 
and its potential for extension to more general set-ups.
As a brief description, consider $n$ IID random variables 
$X_1, \ldots, X_n$ taking values in a measurable space $(\chi, \mathcal{B})$. 
Assume that there is an underlying true probability space $(\Omega, \mathcal{B}_\Omega, P)$
such that, for $i= 1, \ldots, n$, $X_i$ is $\mathcal{B} \slash \Omega$ measurable,
independent with respect to $P$ and it's induced distribution $G(x)$ has an absolutely continuous density $g(x)$ 
with respect to a dominating $\sigma$-finite measure $\lambda(dx)$.
We model $G$ by a parametric family 
$\{F_{\boldsymbol{\theta}} : \boldsymbol{\theta}\in\Theta \subseteq \mathbb{R}^p \}$
which is assumed to be absolutely continuous with respect to $\lambda$ having density $f_{\boldsymbol{\theta}}$.
Consider a prior density for $\boldsymbol{\theta}$ over the parameter space $\Theta$ given by $\pi(\boldsymbol{\theta})$. 
\citet{Ghosh/Basu:2016} defined a robust pseudo-posterior density,  namely the $R^{(\alpha)}$-posterior density
of $\boldsymbol{\theta}$, given the sample $\underline{\boldsymbol{x}}_n = (x_1, \ldots, x_n)^T$ 
on the random variable $\underline{\boldsymbol{X}}_n = (X_1, \ldots, X_n)^T$, as
\begin{equation}
\pi_n^{(\alpha)}(\boldsymbol{\theta}|\underline{\boldsymbol{x}}_n) 
= \frac{\exp(q_n^{(\alpha)}(\underline{\boldsymbol{x}}_n|\boldsymbol{\theta}))\pi(\boldsymbol{\theta})}{
	\int \exp(q_n^{(\alpha)}(\underline{\boldsymbol{x}}_n|\boldsymbol{\theta}'))\pi(\boldsymbol{\theta}') d\boldsymbol{\theta}'},
~~~~\alpha\geq 0,
\label{EQ:R_post_density}
\end{equation}  

\noindent
where  $q_n^{(\alpha)}(\underline{\boldsymbol{x}}_n|\boldsymbol{\theta})$ is the $\alpha$-likelihood 
of $\underline{\boldsymbol{x}}_n$ given by 
\begin{eqnarray}
q_n^{(\alpha)}(\underline{\boldsymbol{x}}_n|\boldsymbol{\theta}) 
&=& \frac{1}{\alpha} \sum_{i=1}^n f_{\boldsymbol{\theta}}^\alpha (x_i) 
- \frac{n}{1+\alpha} \int  f_{\boldsymbol{\theta}}^{1+\alpha} - \frac{n}{\alpha} 
= \sum_{i=1}^n q_{\boldsymbol{\theta}}^{(\alpha)}(x_i), 
\label{EQ:alpha-likelihood}
\end{eqnarray}

\noindent
with $G_n$ being the empirical distribution based on the data and
\begin{eqnarray}
q_{\boldsymbol{\theta}}^{(\alpha)}(y) =  \frac{1}{\alpha}\left(f_{\boldsymbol{\theta}}^\alpha (y) -1\right)
- \frac{1}{1+\alpha} \int  f_{\boldsymbol{\theta}}^{1+\alpha}.
\label{EQ:q_IID}
\end{eqnarray}

\noindent
In a limiting sense,
$q_n^{(0)}(\underline{\boldsymbol{x}}_n|\boldsymbol{\theta}) 
=   \sum_{i=1}^n \left(\log(f_{\boldsymbol{\theta}}(x_i))-1\right)$,
which is the usual log-likelihood (plus a constant); 
so the $R^{(0)}$-posterior is just the ordinary Bayes posterior.
The idea came from a frequentist robust estimator, 
the minimum density power divergence (DPD) estimator (MDPDE) of \cite{Basu/etc:1998},
which has proven to be a useful robust generalization of the maximum likelihood estimator (MLE);
see \cite{Ghosh/Basu:2016} for details.
The similarity of this approach (at $\alpha>0$) with the usual Bayes posterior (at $\alpha=0$) is that, 
it does not require nonparametric smoothing like some other pseudo-posteriors
and it is additive in the data so that the posterior update is easy with new observations.
In \cite{Ghosh/Basu:2016}, its robustness is demonstrated 
and a Bernstein-von Mises type limiting result is proved under the IID set-up. 


\subsection{The Contribution of This Paper}

We provide a generalization of the $R^{(\alpha)}$-posterior density for a completely general 
parametric model set-up beyond IID data, through a suitable structural definition of the $\alpha$-likelihood function,
and derive the exponential convergence results associated with the new pseudo-posterior for the general set-up.
These, in fact, generalize the corresponding results for the usual Bayes posterior \citep{Barron:1988} 
for the  $R^{(\alpha)}$-posterior and their advantages are illustrated through several applications.
Our major contribution in the present paper can be summarized as follows.

\begin{itemize}
	\item This paper is the first to define a robust pseudo-posterior 
	for the general class of parametric models with a finite set of parameters. 
	 All the previous literature on pseudo-posterior are confined to the IID set-up or 
	 a particular example of a non-IID case.
	 Our model set-up is extremely general to cover the IID case
	 as well as every type of non-homogeneous and dependent observations  
	 provided the inference is to be performed based on a finite set of parameters.
	 We have defined a robust $R^{(\alpha)}$-posterior and the associated estimators 
	 for such a general class of  statistical inference problems covering enormous applications.

	\item To illustrate the wide applicability of our proposal, we have explicitly presented the forms of 
	the $R^{(\alpha)}$-posterior or the  $\alpha$-likelihood function for several important cases
	like the independent non-homogeneous data including linear and logistic regressions, 
	time series and Markov models, diffusion processes, etc. 
	Our $R^{(\alpha)}$-posteriors also contain the usual Bayes posterior at $\alpha\rightarrow 0$
	and hence provides a direct generalization of the latter at $\alpha>0$.
	
	\item All the previous pseudo-posteriors currently available in the literature sacrifice
	the conditional probability interpretation 	of the usual Bayes theory. 
	In this paper, for the first time, we discuss a pseudo-posterior,
	namely the $R^{(\alpha)}$-posterior, that retains this conditional probability interpretation 
	with respect to a suitably modified model and modified prior;  
	the $R^{(\alpha)}$-posterior indeed becomes the ordinary Bayes posterior for such a modified set-up (Remark \ref{REM:star}).
	We also introduce the $R^{(\alpha)}$-marginal density of data,
	a robust generalization of the usual marginal.
 
	\item Beyond the methodological proposals, we also establish the theoretical properties of the proposed  $R^{(\alpha)}$-posterior 
	under the fully general parametric set-up.
	We study the asymptotic properties of the  $R^{(\alpha)}$-marginal 
	and the corresponding joint density of data and parameters. 
	We also derive the exponential convergence of the $R^{(\alpha)}$-posterior probabilities 
	and hence the exponential consistency of the associated  $R^{(\alpha)}$-Bayes estimators 
	under the fully general set-up.
	As per our knowledge, such an optimal asymptotic property is not available for any other 
	pseudo-posterior.
	
	\item The assumptions needed for our theoretical derivations are indeed extensions of those 
	required for the classical Bayes theory \citep{Barron:1988}; they are based on the usual concepts of 
	information denseness, merging of distributions in probability, (modified) prior negligibility 
	and the existence of uniform exponential consistent tests.
	We have further simplified these conditions for the IID and the non-homogeneous set-ups.
	They are verified for common examples like linear regression with known or unknown error variance
	and logistic regression models. Although the initial set of conditions under the general parametric models
	look more stringent than the current literature, we have illustrated that they indeed hold under very mild conditions in common examples; e.g., for linear or logistic regressions they are seen to hold only under 
	the boundedness conditions on the fixed design matrix and the positive definiteness of the associated variance matrix. 
	
	\item We have also separately studied the interesting cases of discrete priors under IID set-up, 
	and the associated maximum $R^{(\alpha)}$-posterior estimator with their exponential consistency.
	
	\item Finally, to bridge the gap between the theoretical developments with their practical applicability,
	 we also discuss several important practical issues  
	like the computation of the $R^{(\alpha)}$-posterior and associated estimates
	and  the choice of the tuning parameter $\alpha$. 
	The usefulness of our proposal is illustrated numerically for the 
	linear regression with known and unknown error variance and logistic regression
	along with the corresponding algorithms and R codes.
\end{itemize}

%

For brevity, all proofs and the R-codes are given in the Online Supplement. 

\section{A general form of the $R^{(\alpha)}$-posterior distribution}
\label{SEC:Asymp_Post_setup}

In order to extend the $R^{(\alpha)}$-posterior density to a more general set-up,
let us assume that the random variable $\underline{\boldsymbol{X}}_n$ 
is defined on a general measurable space $(\boldsymbol{\chi}_n, \mathcal{B}_n)$ for each $n$ (sample size).
Also assume that there is an underlying true probability space $(\Omega, \mathcal{B}_\Omega, P)$
such that, for each $n\geq 1$, $\underline{\boldsymbol{X}}_n$ is $\mathcal{B}_n \slash \Omega$ 
measurable and its induced distribution $G^n(\underline{\boldsymbol{x}}_n)$ is absolutely continuous with respect to 
some $\sigma$-finite measure $\lambda^n(d\underline{\boldsymbol{x}}_n)$ 
having ``true" probability density  $g^n(\underline{\boldsymbol{x}}_n)$.
We wish to model it by a parametric family of distributions 
$\mathcal{F}_n = \{F^n(\cdot|\boldsymbol{\theta}) : \boldsymbol{\theta}\in\Theta_n\subseteq \mathbb{R}^p\}$
where the elements of $\mathcal{F}_n$ are assumed to be absolutely continuous with respect to $\lambda^n$ 
having density $f^n(\underline{\boldsymbol{x}}_n|\boldsymbol{\theta})$ for each $n$.
Note that, we have not assumed the parameter space $\Theta_n$ to be independent of the sample size $n$. 
Similarly, the  prior measure $\pi_n(\boldsymbol{\theta})$ on $\Theta_n$ may be $n$-dependent 
with $\pi_n(\Theta_n)\leq 1$. Consider a $\sigma$-field $\mathcal{B}_{\Theta_n}$ on the parameter space $\Theta_n$.
Generalizing from (\ref{EQ:alpha-likelihood}), we propose to define the $\alpha$-likelihood function 
$q_n^{(\alpha)}(\underline{\boldsymbol{x}}_n|\boldsymbol{\theta})$ 
in such a way that ensures
\begin{eqnarray}
~~~~~~~~~~~~~q_n^{(0)}(\underline{\boldsymbol{x}}_n|\boldsymbol{\theta}) := 
\lim\limits_{\alpha\downarrow 0} q_n^{(\alpha)}(\underline{\boldsymbol{x}}_n|\boldsymbol{\theta}) 
= \log f^n(\underline{\boldsymbol{x}}_n|\boldsymbol{\theta}) - n,
\mbox{ for all }\underline{\boldsymbol{x}}_n\in \chi_n.
\label{EQ:alpha-likelihood_gen}
\end{eqnarray}

Our definition should guarantee that the $\alpha$-likelihood, as a function of $\boldsymbol{\theta}$, 
is $\mathcal{B}_{\Theta_n}$ measurable for each $\underline{\boldsymbol{x}}_n$
and jointly $\mathcal{B}_n\times\mathcal{B}_{\Theta_n}$ measurable when both 
$\underline{\boldsymbol{X}}_n$ and $\boldsymbol{\theta}$  are random. 
Then, for this general set-up, we define the corresponding $R^{(\alpha)}$-posterior probabilities as
\begin{equation}
	\pi_n^{(\alpha)}\left(A_n|\underline{\boldsymbol{x}}_n\right) 
	= \frac{\int_{A_n}\exp(q_n^{(\alpha)}(\underline{\boldsymbol{x}}_n|\boldsymbol{\theta}))\pi_n(\boldsymbol{\theta})d\boldsymbol{\theta}}{
		\int_{\Theta_n} \exp(q_n^{(\alpha)}(\underline{\boldsymbol{x}}_n|\boldsymbol{\theta})) \pi_n(\boldsymbol{\theta}) d\boldsymbol{\theta}},
	~~~~~~A_n \in \mathcal{B}_{\Theta_n},
	\label{EQ:R_post_densityGen}
\end{equation}  
whenever the denominator is finitely defined and is positive;
otherwise we may define it arbitrarily, e.g., $\pi_n^{(\alpha)}\left(A_n|\underline{\boldsymbol{x}}_n\right)=\pi_n(A_n)$. 
Definition (\ref{EQ:alpha-likelihood_gen}) ensures that $\pi_n^{(0)}$ is the usual Bayes posterior.

For an useful alternative representation, we define
$
Q_n^{(\alpha)}(S_n|\boldsymbol{\theta}) 
:= \int_{S_n} \exp(q_n^{(\alpha)}(\underline{\boldsymbol{x}}_n|\boldsymbol{\theta})) d\underline{\boldsymbol{x}}_n,
$\\
$M_n^{(\alpha)}(S_n, A_n) 
:= \int_{A_n} Q_n^{(\alpha)}(S_n|\boldsymbol{\theta})\pi_n(\boldsymbol{\theta}) d\boldsymbol{\theta}
$
and 
$M_n^{(\alpha)}(S_n) := M_n^{(\alpha)}(S_n, \Theta_n)/M_n^{(\alpha)}(\boldsymbol{\chi}_n, \Theta_n)$,
for $S_n\in \mathcal{B}_n$ and $A_n\in \mathcal{B}_{{\Theta}_n}$.
In the following, we will assume that the model and priors are chosen to satisfy 
$0 < M_n^{(\alpha)}(\boldsymbol{\chi}_n, \Theta_n) < \infty$.
Then, the last two measures have densities with respect to 
$\lambda^n(d\underline{\boldsymbol{x}}_n)$ given by  
$$
m_n^{(\alpha)}(\underline{\boldsymbol{x}}_n, A_n) 
= \int_{A_n} \exp(q_n^{(\alpha)}(\underline{\boldsymbol{x}}_n|\boldsymbol{\theta}))\pi_n(\boldsymbol{\theta}) 
d\boldsymbol{\theta},
$$
and 
$m_n^{(\alpha)}(\underline{\boldsymbol{x}}_n) 
=m_n^{(\alpha)}(\underline{\boldsymbol{x}}_n, \Theta_n)/M_n^{(\alpha)}(\boldsymbol{\chi}_n, \Theta_n)$, 
respectively. 
Clearly, $m_n^{(\alpha)}(\underline{\boldsymbol{x}}_n)$ is a proper probability density,
which we refer to as the $R^{(\alpha)}$-marginal density of $\underline{\boldsymbol{X}}_n$;
the associated $R^{(\alpha)}$-marginal distribution is $M_n^{(\alpha)}(\cdot)$.
At $\alpha>0$, it provides a robust version of the ordinary Bayes marginal $m_n^{(0)}(\underline{\boldsymbol{x}}_n)$. 
Whenever $0< m_n^{(\alpha)}(\underline{\boldsymbol{x}}_n)<\infty$, 
we can re-express the $R^{(\alpha)}$-posterior probabilities (\ref{EQ:R_post_densityGen}) 
in terms of this $R^{(\alpha)}$-marginal density as
$\pi_n^{(\alpha)}\left(A_n|\underline{\boldsymbol{x}}_n\right) 
	= \frac{m_n^{(\alpha)}(\underline{\boldsymbol{x}}_n, A_n)}{m_n^{(\alpha)}(\underline{\boldsymbol{x}}_n, \Theta_n)}
	= \frac{m_n^{(\alpha)}(\underline{\boldsymbol{x}}_n, A_n)/M_n^{(\alpha)}(\boldsymbol{\chi}_n, \Theta_n)}{
		m_n^{(\alpha)}(\underline{\boldsymbol{x}}_n)}$, 
for $A_n \in \mathcal{B}_{\Theta_n}$.
Then the $R^{(\alpha)}$-Bayes joint posterior law of the parameter $\boldsymbol{\theta}$ 
and the data $\underline{\boldsymbol{X}}_n$ is defined as
\begin{eqnarray}
~~~~~~L_n^{(\alpha) Bayes}\left(d\boldsymbol{\theta}, d\underline{\boldsymbol{x}}_n\right) = 
\pi_n^{(\alpha)}\left(d\boldsymbol{\theta}|\underline{\boldsymbol{x}}_n\right)
M_n^{(\alpha)}\left(d\underline{\boldsymbol{x}}_n\right)
= \frac{M_n^{(\alpha)}(d\underline{\boldsymbol{x}}_n, d\boldsymbol{\theta})}{M_n^{(\alpha)}(\boldsymbol{\chi}_n, \Theta_n)}.
\label{EQ:R_post_jointdensityGen}
\end{eqnarray}

\noindent
This provides a nice interpretation of the quantity $M_n^{(\alpha)}(S_n, A_n)$,
when properly normalized, as the product measure associated with 
the $R^{(\alpha)}$-Bayes joint posterior distribution of $\boldsymbol{\theta}$ and $\underline{\boldsymbol{X}}_n$.
At $\alpha=0$, all these again simplify to the ordinary Bayes measures.

\bigskip
\noindent
\textbf{Example \ref{SEC:Asymp_Post_setup}.1 [Independent Stationary Data]:}\\ 
The simplest possible set-up is that of IID observations as described in Section \ref{SEC:intro}.
In terms of the general notation presented above, we have $\underline{\boldsymbol{X}}_n = (X_1, \ldots, X_n)$
with its observed value $\underline{\boldsymbol{x}}_n=(x_1, \ldots, x_n)$ and
the general measurable space $(\boldsymbol{\chi}_n, \mathcal{B}_n)$ is the $n$-fold product of $(\chi, \mathcal{B})$.
Additionally, we have $G^n(\underline{\boldsymbol{x}}_n) = \prod_{i=1}^{n} G(x_i)$,
$g^n(\underline{\boldsymbol{x}}_n) = \prod_{i=1}^{n} g(x_i)$,
$\lambda^n(d\underline{\boldsymbol{x}}_n)=\prod_{i=1}^{n} \lambda(dx_i)$,
$F^n(\underline{\boldsymbol{x}}_n|\boldsymbol{\theta})=\prod_{i=1}^{n} F_{\boldsymbol{\theta}}(x_i)$,
$f^n(\underline{\boldsymbol{x}}_n|\boldsymbol{\theta})=\prod_{i=1}^{n} f_{\boldsymbol{\theta}}(x_i)$ 
and so $\mathcal{F}_n$ is also the $n$-fold product of the family of individual distributions $F_{\boldsymbol{\theta}}$.
Under these notations, the $\alpha$-likelihood $q_n^{(\alpha)}(\underline{\boldsymbol{x}}_n|\boldsymbol{\theta})$, 
which is given by (\ref{EQ:alpha-likelihood}), satisfies the required measurability assumptions 
along with the condition in (\ref{EQ:alpha-likelihood_gen}).

Then, under suitable assumptions on the prior distribution as before,
the corresponding $R^{(\alpha)}$-posterior distribution is defined by (\ref{EQ:R_post_densityGen})
which is now equivalent to (\ref{EQ:R_post_density}) and can be written as a product of 
stationary independent terms corresponding to each $x_i$ (additivity).
Other related measures can be defined from these quantities; we will come back to them again in  Section \ref{SEC:Asymp_IID}.  
\hfill{$\square$}

\bigskip
\noindent
\textbf{Example \ref{SEC:Asymp_Post_setup}.2 [Independent Non-homogeneous Data]:}\\ 
Suppose $X_1, \ldots, X_n$ are independently but not identically distributed random variables,
where each $X_i$ is defined on a measurable space $(\chi^i, \mathcal{B}^i)$ for $i=1, \ldots, n$. 
Considering an underlying common probability space $(\Omega, \mathcal{B}_\Omega, P)$, 
the random variable $X_i$ is assumed to be $\mathcal{B}^i \slash \Omega$ measurable,
independent with respect to $P$ and its induced distribution $G_i(x)$ has an absolutely continuous density $g_i(x)$ 
with respect to some common dominating $\sigma$-finite measure $\lambda(dx)$, for each $i=1,\ldots,n$.
For each $i$, the true distribution $G_i$ is to be modeled by a parametric family 
$\mathcal{F}^i = \{F_{i,\boldsymbol{\theta}} : \boldsymbol{\theta}\in\Theta\subseteq \mathbb{R}^p \}$
which is absolutely continuous with respect to $\lambda$ having density $f_{i, \boldsymbol{\theta}}$.
Note that, although the densities are potentially different for each $i$, they are assumed to share the common 
unknown parameter $\boldsymbol{\theta}$ leaving us with enough degrees of freedom for estimation of $\boldsymbol{\theta}$.

This set-up of independent non-homogeneous (INH) observations
covers many interesting practical problems, the most common one being the regression with fixed design.
Suppose $\boldsymbol{t}_1, \ldots, \boldsymbol{t}_n$ be $n$ fixed, $k$-variate design points. 
For each $i=1, \ldots, n$, given $\boldsymbol{t}_i$ we independently observe $x_i$
which has the parametric model density 
$f_{i, \boldsymbol{\theta}}(x_i)=f(x_i; \boldsymbol{t}_i, \boldsymbol{\theta})$ 
depending on $\boldsymbol{t}_i$ through a regression structure.
This can, for example, have the form 
\begin{equation}\label{EQ:Reg_Gen}
~~~~E(X_i) = \psi(\boldsymbol{t}_i, \boldsymbol{\beta}), ~~~~i=1, \ldots, n,
\end{equation}
where 
$\boldsymbol{\beta}\subseteq \boldsymbol{\theta}$ is the unknown regression coefficients
and $\psi$ is a suitable link function. 
In general, the unknown parameter $\boldsymbol{\theta}=(\boldsymbol{\beta}, \sigma)$
may additionally contain some variance parameter $\sigma$.
For the subclass of generalized linear models, we take 
$\psi(\boldsymbol{t}_i, \boldsymbol{\beta})=\psi(\boldsymbol{t}_i^T\boldsymbol{\beta})$
and $f$ from the exponential family of distributions. For normal linear regression, 
we have $\psi(\boldsymbol{t}_i, \boldsymbol{\beta})=\boldsymbol{t}_i^T\boldsymbol{\beta}$
and $f$ is the normal density with mean $\boldsymbol{t}_i^T\boldsymbol{\beta}$ and variance $\sigma^2$.
Here, the underlying random variables $X_i$s, associated with observations $x_i$s,
have the INH structure with the common parameter 
$\boldsymbol{\theta}=(\boldsymbol{\beta}, \sigma)$ and the different densities $f_{i, \boldsymbol{\theta}}$.
We can further extend this set-up to include the heterogeneous variances
(by taking different $\sigma_i$ for different $f_{i, \boldsymbol{\theta}}$ but involving some common unknown parameters) 
as a part of our INH set-up.
In terms of the general notation, the random variable $\underline{\boldsymbol{X}}_n = (X_1, \ldots, X_n)$
is defined on the measurable space $(\boldsymbol{\chi}_n, \mathcal{B}_n)=\otimes_{i=1}^n (\chi^i, \mathcal{B}^i)$,
and we have
$G^n(\underline{\boldsymbol{x}}_n) = \prod_{i=1}^{n} G_i(x_i)$,
$g^n(\underline{\boldsymbol{x}}_n) = \prod_{i=1}^{n} g_i(x_i)$,
$\lambda^n(d\underline{\boldsymbol{x}}_n)=\prod_{i=1}^{n} \lambda(x_i)$,
$F^n(\underline{\boldsymbol{x}}_n|\boldsymbol{\theta})=\prod_{i=1}^{n} F_{i,\boldsymbol{\theta}}(x_i)$
and $f^n(\underline{\boldsymbol{x}}_n|\boldsymbol{\theta})=\prod_{i=1}^{n} f_{i,\boldsymbol{\theta}}(x_i)$ 
so that $\mathcal{F}_n=\otimes_{i=1}^n\mathcal{F}^i$.

Now, under this INH set-up, we can define the $R^{(\alpha)}$-posterior 
by suitably extending the definition of the  $\alpha$-likelihood function 
$q_n^{(\alpha)}(\underline{\boldsymbol{x}}_n|\boldsymbol{\theta})$ from its IID~version in (\ref{EQ:alpha-likelihood})
keeping in mind the general requirement (\ref{EQ:alpha-likelihood_gen}). 
Borrowing ideas  from \cite{Ghosh/Basu:2013}, who have developed the MDPDE for the INH set-up, 
and following the intuition behind the construction of the $\alpha$-likelihood (\ref{EQ:alpha-likelihood})
of \cite{Ghosh/Basu:2016}, one possible extended definition for the $\alpha$-likelihood 
in the INH case can be given by 
\begin{eqnarray}
~~~~~~ q_n^{(\alpha)}(\underline{\boldsymbol{x}}_n|\boldsymbol{\theta}) 
&=&  \sum_{i=1}^n \left[\frac{1}{\alpha}f_{i,\boldsymbol{\theta}}^\alpha (x_i) 
- \frac{1}{1+\alpha} \int  f_{i,\boldsymbol{\theta}}^{1+\alpha}\right] - \frac{n}{\alpha} 
= \sum_{i=1}^n q_{i,\boldsymbol{\theta}}^{(\alpha)}(x_i), 
\label{EQ:alpha-likelihood_INH}
\end{eqnarray}
with $q_{i,\boldsymbol{\theta}}^{(\alpha)}(y) =  \frac{1}{\alpha}\left(f_{i, \boldsymbol{\theta}}^\alpha (y) -1\right)
- \frac{1}{1+\alpha} \int  f_{i, \boldsymbol{\theta}}^{1+\alpha}$. 
Note that, we have $q_n^{(0)}(\underline{\boldsymbol{x}}_n|\boldsymbol{\theta}) 
=  \sum_{i=1}^n \left(\log(f_{i,\boldsymbol{\theta}}(x_i))-1\right)$,
satisfying the required condition in (\ref{EQ:alpha-likelihood_gen}).
So, assuming a suitable prior for $\boldsymbol{\theta}$,
the $R^{(\alpha)}$-posterior for the INH observations is defined through (\ref{EQ:R_post_densityGen}) 
with $q_n^{(\alpha)}(\underline{\boldsymbol{x}}_n|\boldsymbol{\theta})$
being given by (\ref{EQ:alpha-likelihood_INH}).
Note that, the resulting posterior is again a product of independent but non-homogeneous terms.
We will discuss their properties
in detail in Section \ref{SEC:Asymp_INH}.
\hfill{$\square$}

\begin{remark}
	\label{REM:star}
In the first introduction of the $R^{(\alpha)}$-posterior under IID set-up \citep{Ghosh/Basu:2016}, 
it was noted that its only drawback is the loss of the probabilistic interpretation.  
Here also, so far, we have defined the $R^{(\alpha)}$-posterior differently than 
the conditional probability approach of the usual Bayes theory and called it a pseudo-posterior. 
But, in fact, it can also be interpreted as an ordinary Bayes posterior 
under a suitably modified model and prior, 
which becomes prominent in our general set-up.
To see this, define an $\alpha$-modified model density 
$\widetilde{q}_n^{(\alpha)}(\underline{\boldsymbol{x}}_n|\boldsymbol{\theta})
= \frac{\exp(q_n^{(\alpha)}(\underline{\boldsymbol{x}}_n|\boldsymbol{\theta}))}{
	Q_n^{(\alpha)}(\boldsymbol{\chi}_n|\boldsymbol{\theta})}$
and the $\alpha$-modified prior density  
$\widetilde{\pi}_n^{(\alpha)}(\boldsymbol{\theta})
= \frac{Q_n^{(\alpha)}(\boldsymbol{\chi}_n|\boldsymbol{\theta})\pi_n(\boldsymbol{\theta})}{
	M_n^{(\alpha)}(\boldsymbol{\chi}_n, {\Theta}_n)}$.
Both are proper densities and satisfy the required measurability assumptions
whenever the relevant integrals exist finitely. 
Further, $\widetilde{\pi}_n^{(\alpha)}(\boldsymbol{\theta})$ is a function of $\boldsymbol{\theta}$ only 
(independent of the data) and hence may be used as a prior density in Bayesian inference;
but it depends on $\alpha$ and the model. In particular, at $\alpha=0$, 
$\widetilde{\pi}_n^{(0)}(\boldsymbol{\theta})=\pi_n(\boldsymbol{\theta})$ and 
$\widetilde{q}_n^{(0)}(\underline{\boldsymbol{x}}_n|\boldsymbol{\theta})=f^n(\underline{\boldsymbol{x}}_n|\boldsymbol{\theta})$
so they indeed represent modifications of the model and the prior, respectively, 
in order to achieve robustness against data contamination. 
Now, for any measurable $A_n \in \mathcal{B}_{\Theta_n}$, 
the standard Bayes (conditional)  posterior probability of $A_n$ 
with respect to the ($\alpha$-modified) model family 
$\mathcal{F}_{n,\alpha} =\left\{\widetilde{q}_n^{(\alpha)}(\cdot|\boldsymbol{\theta}) 
: \boldsymbol{\theta}\in \Theta_n \right\}$ and the ($\alpha$-modified) prior 
$\widetilde{\pi}_n^{(\alpha)}(\boldsymbol{\theta})$  is given by 
$\frac{\int_{A_n}\widetilde{q}_n^{(\alpha)}(\underline{\boldsymbol{x}}_n|\boldsymbol{\theta})
	\widetilde{\pi}_n^{(\alpha)}(\boldsymbol{\theta})d\boldsymbol{\theta}}{
	\int_{\Theta_n} \widetilde{q}_n^{(\alpha)}(\underline{\boldsymbol{x}}_n|\boldsymbol{\theta})
	\widetilde{\pi}_n^{(\alpha)}(\boldsymbol{\theta}) d\boldsymbol{\theta}}$,
which simplifies to $\pi_n^{(\alpha)}\left(A_n|\underline{\boldsymbol{x}}_n\right)$
as in (\ref{EQ:R_post_densityGen}).
\end{remark}

\bigskip
In the following we briefly present the forms of the $\alpha$-likelihood 
for some other practically important model set-ups, but their detailed investigations are kept for the future.

\noindent
 \textbf{Example \ref{SEC:Asymp_Post_setup}.3 [Time Series Data]:}\\ 
Consider the true probability space $(\Omega, \mathcal{B}_\Omega, P)$ and an index set $T$. 
A measurable time series $X_t(\omega)$ is a function defined on $T\times\Omega$, 
which is a random variable on $(\Omega, \mathcal{B}_\Omega, P)$ for each $t\in T$. 
Given a time series $\left\{X_t(\omega) : t\in T \right\}$, 
they are assumed to be associated with an increasing sequence of sub $\sigma$-fields $\{\mathcal{G}_t\}$
and have absolute continuous densities $g(X_t|\mathcal{G}_t)$ for $t\in T$. 
For a stationary time series, one might take $\mathcal{G}_t=\mathcal{F}_{t-1}$, 
the $\sigma$-field generated by $\{X_{t-1}, X_{t-2}, \ldots\}$, for each $t\in T$.
In parametric inference, we model $g(X_t|\mathcal{G}_t)$ by a parametric density $f_{\boldsymbol{\theta}}(X_t|\mathcal{F}_{t-1})$ 
and try to infer about the unknown parameter $\boldsymbol{\theta}$ from an observed sample 
$\underline{\boldsymbol{x}}_n=\left\{x_t : t\in \{1, 2, \ldots, n\} \right\}$ of size $n$.
For example, in a Poisson autoregressive model, we assume $f_{\boldsymbol{\theta}}(x_t|\mathcal{F}_{t-1})$
to be a Poisson density with mean $\lambda_t=h_{\boldsymbol{\theta}}(\lambda_{t-1}, X_{t-1})$ for all $t\in T=\mathbb{Z}$
and some known function $h_{\boldsymbol{\theta}}$ involving the unknown parameter 
$\boldsymbol{\theta} \in \Theta\subseteq \mathbb{R}^p$.
In the Bayesian paradigm, we additionally assume a prior density $\pi(\boldsymbol{\theta})$
and update it to get inference based on the posterior density of $\boldsymbol{\theta}$ given the observed sample data.
We can develop the robust Bayesian inference for any such time series model through the proposed $R^{(\alpha)}$-posterior 
density provided a suitable $\alpha$-likelihood function can be defined. 
Following the construction of the MDPDE in such time series models 
\citep[][among others]{Kim/Lee:2011,Kim/Lee:2013,Kang/Lee:2014}, 
we can define the corresponding $\alpha$-likelihood function as 
\begin{eqnarray}
q_n^{(\alpha)}(\underline{\boldsymbol{x}}_n|\boldsymbol{\theta}) 
&=& \sum_{t=1}^n \left[ \frac{1}{\alpha}  f_{\boldsymbol{\theta}}^\alpha(x_t|\mathcal{F}_{t-1}) 
- \frac{1}{1+\alpha} \int  f_{\boldsymbol{\theta}}^{1+\alpha}(x|\mathcal{F}_{t-1})dx\right] - \frac{n}{\alpha}. 
\end{eqnarray}

We have $q_n^{(0)}(\underline{\boldsymbol{x}}_n|\boldsymbol{\theta}) 
 =  \sum_{i=1}^n \left(\log(f_{\boldsymbol{\theta}}(x_t|\mathcal{F}_{t-1}))-1\right)$,
which satisfies the required Condition (\ref{EQ:alpha-likelihood_gen}).
Robust $R^{(\alpha)}$-posterior inference about $\boldsymbol{\theta}$ 
can be developed using this $\alpha$-likelihood function. 
\hfill{$\square$}

\smallskip
\noindent
\textbf{Example \ref{SEC:Asymp_Post_setup}.4 [Markov Process]:}\\ 
Example \ref{SEC:Asymp_Post_setup}.3 can be easily generalized to Markov processes with stationary transitions.
Consider the random variables $X_1, \ldots, X_n$ defined on the underlying true probability space 
$(\Omega, \mathcal{B}_\Omega, P)$ having true transition probabilities $g(X_{k+1}|X_k)$, $k=0,1, 2, \ldots, n-1$,
with $X_0$ being the initial value of the process. 
We model it by a parametric family of stationary probabilities $f_{\boldsymbol{\theta}}(X_{k+1}|X_k)$
depending on some unknown parameter $\boldsymbol{\theta}\in \Theta\subseteq \mathbb{R}^p$.
Then, the $\alpha$-likelihood function given the sample $\underline{\boldsymbol{x}}_n=(x_1, \ldots, x_n)$
can be defined as 
$$q_n^{(\alpha)}(\underline{\boldsymbol{x}}_n|\boldsymbol{\theta}) 
= \sum_{k=1}^n \left[\frac{1}{\alpha}  f_{\boldsymbol{\theta}}^\alpha(x_{k+1}|x_k)
- \frac{1}{1+\alpha} \int  f_{\boldsymbol{\theta}}^{1+\alpha}(x|x_k)dx\right] - \frac{n}{\alpha}.$$ 
Clearly it satisfies Condition (\ref{EQ:alpha-likelihood_gen}) and
it is possible to perform robust $R^{(\alpha)}$-Bayes inference about $\boldsymbol{\theta}$ under this set-up.
 \hfill{$\square$}

 \bigskip
 \noindent
 \textbf{Example \ref{SEC:Asymp_Post_setup}.5 [Diffusion Process]:}\\ 
Consider again a (true) probability space $(\Omega, \mathcal{B}_\Omega, P)$ and an index set  $T$. 
A measurable random variable $X_t$ defined on $T\time\Omega$ follows a diffusion process if  
 $dX_t = a(X_t, \boldsymbol{\mu})dt + b(X_t, \sigma)dW_t$, $t\geq 0$,
 with $X_0=x_0$ and two known functions $a$ and $b$, where $\left\{W_t : t \geq 0 \right\}$ is a standard Wiener process
 and the parameter of interest is $\boldsymbol{\theta} = (\boldsymbol{\mu}, \sigma)^T \in \Theta$, 
 a convex compact subset of $\mathbb{R}^p\times\mathbb{R}^{+}$.
 This model has important applications in finance,
 where some inference about $\boldsymbol{\theta}$ is desired based on discretized observations 
 $X_{t_i^n}$, $i=1, \ldots, n$, from the above diffusion process.
We generally assume $t_{i}^n=i h_n$ with $h_n\rightarrow 0 $ and $n h_n \rightarrow\infty$ as $n\rightarrow\infty$.
Robust (frequentist) MDPDEs of $\boldsymbol{\theta}$ based on such observations
are developed for two of its special cases, $a(X_t, \mu)=a(X_t)$
and $b(X_t, \sigma)=\sigma$, respectively, by \cite{Song/etc:2007} and \cite{Lee/Song:2013}.
However, whenever we have some prior knowledge about $\boldsymbol{\theta}$, 
quantified through a prior $\pi(\boldsymbol{\theta})$, 
one would apply the Bayesian approach. 
A robust Bayes inference can be done by using our $R^{(\alpha)}$-posterior. 
For this purpose, we note that 
$X_{t_i^n} = X_{t_{i-1}^n} + a(X_{t_{i-1}^n}, \boldsymbol{\mu})h_n + b(X_{t_{i-1}^n}, \sigma)\sqrt{h_n}Z_{n,i}+\Delta_{n,i}, ~~~i=1, \ldots, n$,
where
$\Delta_{n,i} = \int_{t_{i-1}^n}^{t_i^n} \left[a(X_{s}, \boldsymbol{\mu}) - a(X_{t_{i-1}^n}, \boldsymbol{\mu})\right]ds$
$+$ $\int_{t_{i-1}^n}^{t_i^n}  \left[b(X_{s}, \sigma) - b(X_{t_{i-1}^n}, \sigma)\right]dW_s$
and
$Z_{n,i} = h_n^{-1/2}\left(W_{t_i^n} - W_{t_{i-1}^n}\right)$. 
Clearly, $Z_{n,i}$ are IID standard normal variables for $i=1, \ldots, n$.
Therefore, whenever $\Delta_{n,i}$ can be ignored in $P$-probability, for large enough $n$,
$X_{t_i^n}|\mathcal{G}_{i-1}^n$, $i= 1, \ldots, n$, 
behave as INH variables with densities 
$
f_{i,\boldsymbol{\theta}}(\cdot|\mathcal{G}_{i-1}^n) \equiv 
N\left(X_{t_{i-1}^n}+a(X_{t_{i-1}^n}, \boldsymbol{\mu})h_n,  b(X_{t_{i-1}^n}, \sigma)^2h_n\right),
$
where $\mathcal{G}_{i-1}^n$  is the $\sigma$-field generated by $\left\{W_s : s\leq t_i^n \right\}$.
Then, the corresponding $\alpha$-likelihood function based on the observed data 
$\underline{\boldsymbol{x}}_n=(x_{t_1^n}, \ldots, x_{t_n^n})$ can be derived as in Example \ref{SEC:Asymp_Post_setup}.3. 
It satisfies the general requirement (\ref{EQ:alpha-likelihood_gen}) and 
has the simplified form, 
$q_n^{(\alpha)}(\underline{\boldsymbol{x}}_n|\boldsymbol{\theta}) 
=\sum_{i=1}^n q_{i,\boldsymbol{\theta}}^{(\alpha)}(x_{t_i^n})$,
with
\begin{eqnarray}
q_{i,\boldsymbol{\theta}}^{(\alpha)}(x_{t_i^n}) 
&=& \left\{ \begin{array}{l l}
\frac{1}{\left(2\pi b(x_{t_{i-1}^n}, \sigma)^2h_n\right)^{\alpha/2}}\left[
	\frac{1}{\alpha}e^{-\frac{\alpha\left(x_{t_i^n} - x_{t_{i-1}^n} - a(x_{t_{i-1}^n}, \boldsymbol{\mu})h_n\right)^2}{
			2 b(x_{t_{i-1}^n}, \sigma)^2h_n}} 
 - \frac{1}{(1+\alpha)^{3/2}}\right] - \frac{1}{\alpha}, 
&\mbox{if } \alpha>0,  \nonumber\label{EQ:alpha-likelihood_diff}\\
 -\frac{\alpha\left(x_{t_i^n} - x_{t_{i-1}^n} - a(x_{t_{i-1}^n}, \boldsymbol{\mu})h_n\right)^2}{
	2b(x_{t_{i-1}^n}, \sigma)^2h_n}
- \frac{1}{2}\log\left(2\pi b(x_{t_{i-1}^n}, \sigma)^2h_n\right) - 1, 
&\mbox{if } \alpha=0. 
\end{array}
\right.\nonumber
 \end{eqnarray}
The robust $R^{(\alpha)}$-posterior 
can be easily obtained using this $\alpha$-likelihood function. 
 \hfill{$\square$}

\section{Exponential Convergence Results under the General Set-up}
\label{SEC:Asymp_Gen}

Exponential consistency is an important property of posterior (Bayes) inference;
it was first demonstrated in \cite{Barron:1988} and  later refined by several authors 
\citep[see][among others]{Ghosal/etc:2000,Walker:2004,Ghosal/vanderVart:2007,Walker/etc:2007}.
We follow the approach of \cite{Barron:1988} to show 
that of our new robust $R^{(\alpha)}$-posterior probabilities and the corresponding parameter estimates 
also enjoy such asymptotic optimality properties.

\subsection{Properties of the Joint and Marginal $R^{(\alpha)}$-Bayes distributions}
\label{SEC:Merging}

Let us recall the general set-up of Section \ref{SEC:Asymp_Post_setup} 
along with the $\alpha$-modified model and prior densities 
$\widetilde{q}_n^{(\alpha)}(\cdot|\boldsymbol{\theta})$ and $\widetilde{\pi}_n^{(\alpha)}(\boldsymbol{\theta})$ 
as defined in Remark \ref{REM:star}.
Consider the Kullback-Leibler divergence  between two absolutely continuous densities
$f_1$ and $f_2$ with respect to the common $\sigma$-finite measure $\lambda$ defined as 
$KLD(f_1, f_2) = \int f_1 \log\left(\frac{f_1}{f_2}\right)d\lambda$,
and put 
$D_n^{(\alpha)}(\boldsymbol{\theta}) = \frac{1}{n} KLD\left(g^n(\cdot), \widetilde{q}_n^{(\alpha)}(\cdot|\boldsymbol{\theta})\right)$.
We define a joint (frequentist) law of $\boldsymbol{\theta}$ and $\underline{\boldsymbol{X}}$ given by
$L_n^{*(\alpha)}\left(d\boldsymbol{\theta}, d\underline{\boldsymbol{x}}_n\right) = 
\pi_n^{*(\alpha)}\left(d\boldsymbol{\theta}\right)G_n\left(d\underline{\boldsymbol{x}}_n\right)$,
where the probability distribution $\pi_n^{*(\alpha)}$ of $\boldsymbol{\theta}$ on $\Theta_n$ is defined as
$
\pi_n^{*(\alpha)}\left(d\boldsymbol{\theta}\right) = \frac{e^{-nD_n^{(\alpha)}(\boldsymbol{\theta})}\widetilde{\pi}_n^{(\alpha)}(d\boldsymbol{\theta})}{c_n},
$
with 
$c_n = \int e^{-nD_n^{(\alpha)}(\boldsymbol{\theta})}\widetilde{\pi}_n^{(\alpha)}(d\boldsymbol{\theta})$. 
We show that this joint law $L_n^{*(\alpha)}$ provides a frequentist large-deviation approximation 
to the joint $R^{(\alpha)}$-Bayes distribution  (\ref{EQ:R_post_jointdensityGen}) 
of $\boldsymbol{\theta}$ and $\underline{\boldsymbol{X}}_n$; 
to quantify their closeness we consider the concept of 
``\textit{merging}" of probability distributions \cite{Barron:1988}.

%

\begin{definition}
	\label{DEF:merge_in_prob}
	Consider two probability distributions $G_1^n$ and $G_2^n$ of $\underline{\boldsymbol{X}}_n$,
	having densities $g_1^n$ and $g_2^n$ respectively with respect to $\lambda^n$. 
\begin{itemize}
	\item They are said to	{\it merge in probability} if for all $\epsilon>0$,
	$
	\lim\limits_{n\rightarrow\infty} P\left(\frac{g_2^n(\underline{\boldsymbol{X}}_n)}{g_1^n(\underline{\boldsymbol{X}}_n)}
	> e^{-n\epsilon}\right) =1.
	$		
	\item 	They {\it merge with probability one} if for every $\epsilon>0$,
	$
	P\left(\frac{g_2^n(\underline{\boldsymbol{X}}_n)}{g_1^n(\underline{\boldsymbol{X}}_n)}
	> e^{-n\epsilon} ~\mbox{ for all large }n\right) =1.
	$	
\end{itemize}	
\end{definition} 

%
An application of Markov's inequality shows that Definition \ref{DEF:merge_in_prob}
is equivalent to the conditions 
$
\lim\limits_{n\rightarrow\infty} \frac{1}{n} \log \frac{g_2^n(\underline{\boldsymbol{X}}_n)}{g_1^n(\underline{\boldsymbol{X}}_n)} =0
$
in probability or with probability one, respectively. 
See \citet[][Section 4]{Barron:1988} for more results on merging. 
Additionally we assume the following condition.

\noindent
\textbf{Assumption (M1):} For any $\epsilon, r>0$, 
there exists a positive integer $N$ such that \\
$\widetilde{\pi}_n^{(\alpha)}\left(\left\{\boldsymbol{\theta} : D_n^{(\alpha)}(\boldsymbol{\theta}) < \epsilon \right\}\right)
\geq e^{-nr}$, for all $n \geq N$.


\begin{theorem}
	Under Assumption (M1), we have the following results.
\begin{itemize}
	\item[a)] 	$\lim\limits_{n\rightarrow\infty} \frac{1}{n} KLD\left(L_n^{*(\alpha)}, L_n^{(\alpha)Bayes}\right) = 0,
	$
  and $	\lim\limits_{n\rightarrow\infty} \frac{1}{n} 
	E_{G^n}\left[KLD\left(\pi_n^{*(\alpha)}(\cdot), \pi_n^{(\alpha)}(\cdot|\underline{\boldsymbol{X}}_n)\right)\right] = 0.
	$
	\item[c)] $\lim\limits_{n\rightarrow\infty} \frac{1}{n} KLD(g^n, m_n^{(\alpha)}) = 0,$
	so that $G^n$ and $M_n^{(\alpha)}$ merge in probability.
\end{itemize}
	\label{THM:Merging1}
\end{theorem}

Although Assumption (M1) might look a bit complicated, 
it can be further simplified in terms of the common notion of \textit{information denseness} 
of priors $\pi_n$ with respect to a suitable  family of model densities.
This notion of information denseness is frequently used in large sample analyses of usual Bayesian methods 
and is precisely defined below for our context.

\begin{definition}
	\label{DEF:info_dense}
	Suppose $\Theta_n=\Theta$ is independent of $n$ and we define 
	$\bar{D}^{(\alpha)}(\boldsymbol{\theta}) =\limsup\limits_{n\rightarrow\infty}D_n^{(\alpha)}(\boldsymbol{\theta})$.
	%
	Then, the prior sequence $\pi_n$ is said to be information dense at $G^n$ with respect to 
	$\mathcal{F}_{n,\alpha} =\left\{\widetilde{q}_n^{(\alpha)}(\cdot|\boldsymbol{\theta}) 
	: \boldsymbol{\theta}\in \Theta_n \right\}$
	if there exists a finite measure $\widetilde{\pi}$ 
such that 
	$\widetilde{\pi}\left(\left\{\boldsymbol{\theta} : \bar{D}^{(\alpha)}(\boldsymbol{\theta}) < \epsilon \right\}\right)> 0$,
	for all $\epsilon>0$, and 
\begin{eqnarray}
\liminf_{n\rightarrow\infty} e^{nr} \frac{d\widetilde{\pi}_n^{(\alpha)}}{d\widetilde{\pi}}(\boldsymbol{\theta}) \geq 1,
~~\mbox{ for all } r>0, \boldsymbol{\theta}\in \Theta.
\label{EQ:infodense_prior}
\end{eqnarray}	
\end{definition} 

\begin{theorem}
	If the prior is information dense with respect to $\mathcal{F}_{n,\alpha}$ as in Definition \ref{DEF:info_dense}, 
	then Assumption (M1) holds and hence the three results of Theorem \ref{THM:Merging1} also hold.
	\label{THM:Merging2}
\end{theorem}
%
%
%

\subsection{Consistency of the $R^{(\alpha)}$-Posterior Probabilities}
\label{SEC:Asymp_PD}

We now prove the exponential convergence results for our robust $R^{(\alpha)}$-posterior probabilities. 
For measurable sets $A_n, B_n, C_n \subseteq \Theta_n$ 
and constants $b_n, c_n$, we assume the following.

\begin{itemize}
\item[(A1)] $A_n$, $B_n$ and $C_n$ together complete $\Theta_n$, i.e., $A_n\cup B_n \cup C_n =\Theta_n$, for each $n\geq1$.

\item[(A2)] $B_n$ satisfies $\widetilde{\pi}_n^{(\alpha)}(B_n) 
= \frac{M_n^{(\alpha)}(\boldsymbol{\chi}_n, B_n)}{M_n^{(\alpha)}(\boldsymbol{\chi}_n, \Theta_n)}\leq b_n$, for each $n\geq1$.

\item[(A3)] $\{C_n\}$ is such that there exists $S_n \in \mathcal{B}_n$ satisfying  
$
\lim\limits_{n\rightarrow\infty} G^n\left(S_n\right) = 0,
\sup\limits_{\boldsymbol{\theta}\in C_n} \frac{Q_n^{(\alpha)}(S_n^c|\boldsymbol{\theta})}{Q_n^{(\alpha)}(\boldsymbol{\chi}_n|\boldsymbol{\theta})}
\leq c_n.
$

\item[(A3)$^\ast$] $\{C_n\}$ is such that there exists $S_n \in \mathcal{B}_n$ satisfying  
$P\left(\underline{\boldsymbol{X}}_n \in S_n~\mbox{ i.o.}\right)$ $=0$
and\\ 
$\sup\limits_{\boldsymbol{\theta}\in C_n} \frac{Q_n^{(\alpha)}(S_n^c|\boldsymbol{\theta})}{Q_n^{(\alpha)}(\boldsymbol{\chi}_n|\boldsymbol{\theta})}
\leq c_n,$
where i.o.~denotes ``{\it infinitely often}".
\end{itemize}

Here we need either Condition (A3) or Condition (A3)$^\ast$ which, respectively, help us to prove 
the convergence results in probability or with probability one. 
Condition (A3)$^\ast$ is stronger and imply (A3),
but (A3) is sufficient in most practices yielding a convergence in probability type result.
Also, if Condition (A3) holds with $c_n = e^{-nr}$ for some $r>0$, 
then it ensures the existence of a uniformly exponentially consistent (UEC) test 
for $G^n$ against the family of $\alpha$-modified probability distributions
$\left\{\frac{Q_n^{(\alpha)}(\cdot|\boldsymbol{\theta})}{Q_n^{(\alpha)}(\chi_n|\boldsymbol{\theta})} 
: \boldsymbol{\theta} \in C_n \right\}$
corresponding to the $\alpha$-modified model density $\widetilde{q}_n^{(\alpha)}(\cdot|\boldsymbol{\theta})$
defined in Remark \ref{REM:star}.
%
Although complex looking, these conditions
are straightforward extensions of the conditions used by \cite{Barron:1988} 
for proving the exponential convergence of ordinary Bayes posterior probabilities;
they indeed coincide  at $\alpha=0$. 
In particular, at $\alpha=0$, Condition (A2) simplifies to 
$\pi_n(B_n) \leq b_n$, i.e., $B_n$ have negligible prior probabilities if $b_n \rightarrow 0$, 
and (A3) assumes the existence of a UEC test against the models with $\boldsymbol{\theta} \in C_n$. 
Under these conditions, along with the concept of merging (Subsection \ref{SEC:Merging}), 
we have the following main theorem.  

\begin{theorem}
{[Exponential Consistency of $R^{(\alpha)}$-posterior probabilities]}
\begin{enumerate}
\item[(1)] Suppose that $G^n$ and $M_n^{(\alpha)}(\cdot)$ merge in probability and 
let $A_n \in \mathcal{B}_{\Theta_n}$ be any sequence of sets.
Then,
$\limsup\limits_{n\rightarrow\infty} P\left(\pi_n^{(\alpha)}\left(A_n^c|\underline{\boldsymbol{X}}_n\right)< e^{-nr} \right)
=1$,
for some $r>0$, if and only if there exist  $r_1, r_2 >0$ and sets $B_n, C_n \in \mathcal{B}_{\Theta_n}$
such that (A1)--(A3) are satisfied with $b_n = e^{-nr_1}$ and $c_n = e^{-nr_2}$.

\item[(2)] Suppose that $G^n$ and  $M_n^{(\alpha)}(\cdot)$ merge with probability one and 
let $A_n \in \mathcal{B}_{\Theta_n}$ be any sequence of sets.
Then,
$P\left(\pi_n^{(\alpha)}\left(A_n^c|\underline{\boldsymbol{X}}_n\right) \geq e^{-nr} \mbox{ i.o.}\right)=0$,
for some $r>0$, if and only if there exists constants $r_1, r_2 >0$ and sets $B_n, C_n \in \mathcal{B}_{\Theta_n}$
such that  Assumptions (A1), (A2) and (A3)$^\ast$ are satisfied with $b_n = e^{-nr_1}$ and $c_n = e^{-nr_2}$.
\end{enumerate}
\label{THM:PostConv}
\end{theorem}

Note that, for $\alpha=0$, Theorem \ref{THM:PostConv} coincides with the classical exponential convergence results 
of ordinary Bayes posterior probabilities as proved in \cite{Barron:1988}. 
Our theorem generalizes it for the robust $R^{(\alpha)}$-posterior probabilities under suitable conditions. 
Hence, the $R^{(\alpha)}$-posterior distribution, besides yielding robust results under data contamination, 
is asymptotically optimal in exactly the same exponential rate as the ordinary posterior for all $\alpha\geq 0$.
%

\subsection{Consistency of the $R^{(\alpha)}$-Bayes Estimators}
\label{SEC:Asymp_BE}

Let us now examine the asymptotic properties of the  $R^{(\alpha)}$-Bayes estimators 
associated with the $R^{(\alpha)}$-posterior distribution (\ref{EQ:R_post_densityGen}) 
under the general set-up of Section \ref{SEC:Asymp_Post_setup}. 
In the decision-theoretic framework, we consider the problem of estimation of 
a functional $\phi_P :=\phi(P)$ of the true probability $P$;
for example $\phi_P$ could be the probability density of $P$, or any summary measure (like mean) of $P$. 
For the given parametric family $F^n(\cdot|\boldsymbol{\theta})$, let us denote 
$\phi_{\boldsymbol{\theta}} := \phi_{F^n(\cdot|\boldsymbol{\theta})}$.
Then, our action space is $\Phi = \left\{ \phi_Q : Q \mbox{ is a probability measure on } (\Omega, \mathcal{B}_\Omega)\right\}$.
Consider a non-negative loss function $L_n(\phi, \widehat{\phi})$ on $\Phi \times \Phi$ 
denoting the loss in estimating $\phi$ by $\widehat{\phi}$;
let $L_n(\phi_{\boldsymbol{\theta}}, \phi) $ is $\mathcal{B}_{\Theta_n}$ measurable for each $\phi\in \Phi$.
Then the general $R^{(\alpha)}$-Bayes estimator 
$\widehat{\phi} = \widehat{\phi}(\cdot; \underline{\boldsymbol{x}}_n)$ of $\phi$ 
is defined as 
\begin{eqnarray}
\widehat{\phi} = \arg\min_{\phi \in \Phi} \int L_n(\phi_{\boldsymbol{\theta}}, \phi) 
\pi_n^{(\alpha)}\left(d\boldsymbol{\theta}|\underline{\boldsymbol{x}}_n\right),
\label{EQ:R-Est_Gen}
\end{eqnarray} 

\noindent
provided the minimum is attained. 
In particular, the $R^{(\alpha)}$-Bayes estimator of $\phi_{\boldsymbol{\theta}} = \boldsymbol{\theta}$ is 
the mean of the $R^{(\alpha)}$-posterior distribution for squared error loss provided it exists finitely, 
or a median of the $R^{(\alpha)}$-posterior distribution for absolute error loss. 

However, if the minimum in (\ref{EQ:R-Est_Gen})  is not attained,
we may define the approximate $R^{(\alpha)}$-Bayes estimator $\widehat{\phi}$ of $\phi$ through the relation
$\int L_n(\phi_{\boldsymbol{\theta}}, \widehat{\phi})\pi_n^{(\alpha)}\left(d\boldsymbol{\theta}|\underline{\boldsymbol{x}}_n\right) 
\leq \inf\limits_{\phi \in \Phi} \int L_n(\phi_{\boldsymbol{\theta}}, \phi)
\pi_n^{(\alpha)}\left(d\boldsymbol{\theta}|\underline{\boldsymbol{x}}_n\right) + \delta_n$,
with $ \lim\limits_{n\rightarrow\infty}\delta_n = 0$.
An useful example is the approximate mode of the $R^{(\alpha)}$-posterior for discrete parameter space,
which is an approximate $R^{(\alpha)}$-Bayes estimator under 0-1 loss.
Also, note that, if the $R^{(\alpha)}$-Bayes estimator exists, 
it is also an approximate $R^{(\alpha)}$-Bayes estimator.

%

\begin{definition}
A loss function $L_n$ on $\Phi \times \Phi$ is said to be bounded if there exists $\bar{L}<\infty$ such that
$L_n(\phi_{\boldsymbol{\theta}}, \phi_P) \leq \bar{L}$ for all $n$ and all $\boldsymbol{\theta}\in \Theta_n$.
\label{DEF:loss_bounded}
\end{definition}

\begin{definition}
A loss $L_n$ on $\Phi \times \Phi$ is said to be equivalent to a pseudo-metric $d_n$ on $\Phi\times \Phi$ 
if there exist two strictly increasing functions $h_1$ and $h_2$ on $[0, \infty)$ 
that are continuous at 0 with $h_1(0)= h_2(0)=0$ and satisfy
$L_n \leq h_1(d_n)$ and $d_n\leq h_2(L_n)$ on $\Phi\times \Phi$ and for all $n$.	
\label{DEF:loss_equiv}
\end{definition}

Note that, Definition \ref{DEF:loss_equiv} indicates
 $\displaystyle\lim\limits_{n\rightarrow\infty} L_n(\phi_n, \widehat{\phi}_n) =0$
if and only if $\displaystyle\lim\limits_{n\rightarrow\infty} d_n(\phi_n, \widehat{\phi}_n) =0$.
As an example, the squared Hellinger loss is bounded and equivalent to 
the $L_1$-distance. Also, the absolute error ($L_1$) loss is equivalent to itself and 
bounded by twice the Hellinger loss.

We now establish the asymptotic consistency of $R^{(\alpha)}$-Bayes and 
approximate $R^{(\alpha)}$-Bayes estimators of $\phi_{\boldsymbol{\theta}}$ to the true value $\phi_P$
for such loss. 
The proof mimics that of Lemma 12 in \cite{Barron:1988}.

\begin{theorem}[Consistency of $R^{(\alpha)}$-Bayes Estimators]
Given any sample data $\underline{\boldsymbol{x}}_n$, let $\widehat{\phi}_n = \widehat{\phi}(\cdot; \underline{\boldsymbol{x}}_n)$
be an approximate $R^{(\alpha)}$-Bayes estimator (or the $R^{(\alpha)}$-Bayes estimator) of $\phi_P$
with respect to a loss function $L_n$ that is bounded 
and equivalent to a pseudo-metric $d_n$.
Also, for any $\epsilon>0$, define 
$A_{\epsilon,n} = \left\{\boldsymbol{\theta} : d_n(\phi_P, \phi_{\boldsymbol{\theta}})\leq \epsilon\right\}$.
Then, we have 
$d_n(\phi_P, \widehat{\phi}_n) \leq \epsilon + h_2\left(\frac{\epsilon +\bar{L} 
	\pi_n^{(\alpha)}\left(A_{h_1^{-1}(\epsilon),n}^c|\underline{\boldsymbol{x}}_n\right) }{
	1 - \pi_n^{(\alpha)}\left(A_{\epsilon,n}^c|\underline{\boldsymbol{x}}_n\right)}\right).$
Consequently, 
if $ \displaystyle\lim\limits_{n\rightarrow\infty}\pi_n^{(\alpha)}\left(A_{\epsilon,n}^c|\underline{\boldsymbol{X}}_n\right)=0$
in probability or with probability one for all $\epsilon>0$, 
then $\displaystyle\lim\limits_{n\rightarrow\infty} d_n(\phi_P, \widehat{\phi}_n) =0$ 
in probability or with probability one, respectively.
\label{THM:Asymp_R-Est}
\end{theorem}

In simple language, Theorem \ref{THM:Asymp_R-Est} states that whenever the target $\phi_P$ is close 
enough to the model value $\phi_{\boldsymbol{\theta}}$ in the pseudo-metric $d_n$ asymptotically 
under the $R^{(\alpha)}$-posterior probability,
the corresponding $R^{(\alpha)}$-Bayes estimator with respect to $L_n$ is asymptotically consistent for $\phi_P$ in $d_n$.
But, Theorem \ref{THM:PostConv} yields
$\displaystyle\lim\limits_{n\rightarrow\infty}\pi_n^{(\alpha)}\left(A_{\epsilon,n}^c|\underline{\boldsymbol{X}}_n\right)=0$
under appropriate conditions and hence the corresponding $R^{(\alpha)}$-Bayes estimators 
are consistent in suitable $d_n$. 
In particular, Theorem \ref{THM:Asymp_R-Est} applies to the $R^{(\alpha)}$-Bayes estimators
with respect to the squared Hellinger loss and the $L_1$-loss.
 to deduce their $L_1$ consistency.

\section{Application (I): Independent Stationary Models}
\label{SEC:Asymp_IID}

\subsection{$R^{(\alpha)}$-Posterior Convergence}
\label{SEC:AsympR_IID}

Consider the set-up of the independent stationary model as in Example \ref{SEC:Asymp_Post_setup}.1.
Let us study the conditions required for the exponential convergence of the $R^{(\alpha)}$-posterior for this particular set-up.
First, to verify the merging of $G^n$ and $M_n^{(\alpha)}$, 
we define the individual $\alpha$-modified density as $\widetilde{q}^{(\alpha)}(\cdot|\boldsymbol{\theta}) = 
\exp\left({q}_{\boldsymbol{\theta}}^{(\alpha)}(\cdot)\right)/Q^{(\alpha)}({\chi}|\boldsymbol{\theta})$
and the $\alpha$-modified prior $\widetilde{\pi}_n^{(\alpha)}$ as in Remark \ref{REM:star} with $\pi_n=\pi$. 
Then we consider the information denseness of the prior $\pi$ under independent stationary models 
with respect to 
$\mathcal{F}_\alpha =\left\{\widetilde{q}^{(\alpha)}(\cdot|\boldsymbol{\theta}) : \boldsymbol{\theta}\in \Theta \right\}$
defined as follows.

\begin{definition}
The prior $\pi$ under the IID model is information dense at $G$ with respect to $\mathcal{F}_\alpha$ if  there exists a finite measure $\widetilde{\pi}$ satisfying (\ref{EQ:infodense_prior})
and 
$\widetilde{\pi}\left(\left\{\boldsymbol{\theta} : KLD(g, \widetilde{q}^{(\alpha)}(\cdot|\boldsymbol{\theta})) 
< \epsilon \right\}\right)> 0$  for all $\epsilon>0$.
\label{DEF:info_dense_iid}
\end{definition}

Note that, the above definition is equivalent to the general notion of information denseness given in 
Definition \ref{DEF:info_dense}. Thus, in view of Theorem \ref{THM:Merging2}, 
it implies the merging of $G^n$ and $M_n^{(\alpha)}$ in probability for  independent stationary models. 
Then, Theorem \ref{THM:PostConv} may be restated as follows.

\begin{proposition}
Consider the set-up of independent stationary models and assume that the prior $\pi$ is independent of $n$ 
and is information dense at $g$ with respect to $\mathcal{F}_\alpha$ as per Definition \ref{DEF:info_dense_iid}.
Take any sequence of measurable parameter sets $A_n \subset \Theta$. Then, 
$\pi_n^{(\alpha)}\left(A_n^c|\underline{\boldsymbol{X}}_n\right)$ is exponentially small with $P$-probability tending to one,
if and only if there exists constants $r_1, r_2 >0$ and sets $B_n, C_n \in \mathcal{B}_{\Theta}$
such that 
such that (A1)--(A3) are satisfies with $b_n = e^{-nr_1}$ and $c_n = e^{-nr_2}$. 
\label{PROP:PostConv_iid1}
\end{proposition}

Next note that, for the present case, (A3) holds under the assumption of the existence of a UEC test 
for $G$ against the family $\left\{\frac{Q^{(\alpha)}(\cdot|\boldsymbol{\theta})}{Q^{(\alpha)}(\chi|\boldsymbol{\theta})} : \boldsymbol{\theta} \in C_n \right\}$.
We can further simplify it by using a necessary and sufficient condition 
for the existence of UEC from \cite{Barron:1989} which states that,
``\textit{for every $\epsilon>0$ there exists a sequence of UEC tests for the hypothesized distribution $P$ versus 
	the family of distributions $\left\{Q : d_{T_n}(P, Q) > \epsilon/2\right\}$ if and only if the sequence of partitions $T_n$ 
	has effective cardinality (eff. card.) of order $n$ with respect to $P$}"; here, for any measurable partition $T$,  
$d_{T}$ denotes the $T$-variation norm $d_T(P, Q) = \sum_{A\in T}\left|P(A) - Q(A)\right|$.
Using this, we show that the $R^{(\alpha)}$-posterior asymptotically concentrates 
on the $L_1$ model neighborhood of the true density $g$.
Define, for any density $p$ and any partition $T$, the ``theoretical histogram" density $p^T$ as
$p^T(x) = \frac{1}{\lambda(A)} \int_A p(y)\lambda(dy)$, 
for $x\in A\in T,$
whenever $\lambda(A)\neq 0$, and $p^T=0$ otherwise. We call a sequence of partitions $T_n$ to be ``\textit{rich}" 
if the corresponding sequence of densities $g^{T_n}$ converges to $g$ in $L_1$-distance. 
Also, define $B_{\epsilon}^{T_n} = \left\{\boldsymbol{\theta} : 
d_1\left(f_{\boldsymbol{\theta}}, \widetilde{q}^{(\alpha)T_n}(\cdot|\boldsymbol{\theta})\right) > \epsilon \right\}$ 
for any $\epsilon>0$ and sequence of partition $T_n$, where $d_1$ denotes the $L_1$ distance,
and consider the following assumption.

\noindent
\textbf{Assumption (B):}\textit{ For $\epsilon>0$,  
$\widetilde{\pi}_n^{(\alpha)} (B_\epsilon^{T_n}) = 
\frac{M_n^{(\alpha)}(\boldsymbol{\chi}_n, B_\epsilon^{T_n})}{M_n^{(\alpha)}(\boldsymbol{\chi}_n, \Theta)}$
is exponentially small for a rich sequence of partitions $T_n$ with eff. card. of order $n$}.

Note that, Assumption (B) implies Assumption (A2) for $B_\epsilon^{T_n}$, 
or any smaller subset of it. So, applying it with 
$B_n=\left\{\boldsymbol{\theta} : d_1(g, f_{\boldsymbol{\theta}}) \geq \epsilon, 
d_{T_n}\left(G, \frac{Q^{(\alpha)}(\cdot|\boldsymbol{\theta})}{Q^{(\alpha)}(\chi|\boldsymbol{\theta})}\right)<\epsilon/2 \right\}
\subset B_{\epsilon/4}^{T_n}$ and the existence result of UEC tests with $C_n=\left\{\boldsymbol{\theta} : 
d_{T_n}\left(G, \frac{Q^{(\alpha)}(\cdot|\boldsymbol{\theta})}{Q^{(\alpha)}(\chi|\boldsymbol{\theta})}\right)>\epsilon/2 \right\}$, 
Proposition \ref{PROP:PostConv_iid1} yields the asymptotic exponential concentration 
of the $R^{(\alpha)}$-posterior probability in the $L_1$-neighborhood 
$A_n =\left\{\boldsymbol{\theta} : d_1(g, f_{\boldsymbol{\theta}}) < \epsilon \right\}$.
Note that, clearly $A_n\cup B_n \cup C_n=\Theta_n$ for these choices.

\begin{theorem}\label{PROP:PostConv_iid2}
Consider the set-up of IID models and assume that the prior $\pi$ is independent of $n$ 
and information dense at $g$ with respect to $\mathcal{F}_\alpha$ as per Definition \ref{DEF:info_dense_iid}.
If Assumption (B) holds then, for every $\epsilon>0$, 
$\pi_n^{(\alpha)}\left( \left\{\boldsymbol{\theta} : d_1(g, f_{\boldsymbol{\theta}}) \geq \epsilon \right\}
|\underline{\boldsymbol{X}}_n\right)$ 
is exponentially small with $P$-probability one.
\end{theorem}

Note that the final Assumption (B) is easy to verify for model and priors belonging 
to the standard exponential family of distributions with exponentially decaying tails.
However, if Assumption (B) does not hold, we can deduce a weaker conclusion 
in terms of $T_n$-variance distance in place of the $L_1$ distance.
The idea goes back to \cite{Barron:1988} for a similar result in case of the ordinary posterior;
an extended version for the $R^{(\alpha)}$-posterior is given in the following.

\begin{theorem}
Consider the set-up of IID models and assume that the prior $\pi$ is independent of $n$ 
and information dense at $g$ with respect to $\mathcal{F}_\alpha$ as per Definition \ref{DEF:info_dense_iid}.
Then, for any sequence of partitions $T_n$	with effective card. of order $n$, 
$
\pi_n^{(\alpha)}\left( \left\{\boldsymbol{\theta} : 
d_{T_n}\left(G, \frac{Q_n^{(\alpha)}(\cdot|\boldsymbol{\theta})}{Q_n^{(\alpha)}(\chi_n|\boldsymbol{\theta})}\right) \geq \epsilon \right\}
\bigg|\underline{\boldsymbol{X}}_n\right)$ is exponentially small
 with $P$-probability one.
\label{PROP:PostConv_iid3}
\end{theorem}


\subsection{The Cases of Discrete Priors: Maximum $R^{(\alpha)}$-Posterior Estimator}
\label{SEC:Asymp_discrete}

We can derive the exponential consistency of the $R^{(\alpha)}$-Bayes estimators 
with respect to the bounded loss functions from Theorem \ref{THM:Asymp_R-Est} 
along with Proposition \ref{PROP:PostConv_iid1}--\ref{PROP:PostConv_iid3}. 
Let us now consider, in more detail, the particular case of discrete priors 
and the maximum $R^{(\alpha)}$-posterior estimator.

Consider the set-up of IID models, but now with a countable $\Theta$.
On this countable parameter space, we consider a sequence of discrete priors $\pi_n(\boldsymbol{\theta})$ 
which are sub-probability mass functions, i.e., $\sum_{\boldsymbol{\theta}}\pi_n(\boldsymbol{\theta}) \leq 1$.
The most common loss-function to consider under this set-up is the 0-1 loss function,
for which the resulting $R^{(\alpha)}$-Bayes estimator is the (global) mode of the $R^{(\alpha)}$-posterior density;
we call this estimator of $\boldsymbol{\theta}$ as the ``maximum $R^{(\alpha)}$-posterior estimator (MRPE)".
When this mode is not attained, we consider an approximate version
$\widehat{\boldsymbol{\theta}}_\alpha$, to be referred to as an 
``approximate maximum $R^{(\alpha)}$-posterior estimator (AMRPE)", 
defined by the relation
\begin{eqnarray}\label{EQ:AMRPE}
\widetilde{\pi}_n^{(\alpha)}(\widehat{\boldsymbol{\theta}}_\alpha) \widetilde{q}_n^{(\alpha)}(\underline{\boldsymbol{x}}_n|\widehat{\boldsymbol{\theta}}_\alpha)
> \sup_{\boldsymbol{\theta}}\widetilde{\pi}_n^{(\alpha)}(\boldsymbol{\theta}) \widetilde{q}_n^{(\alpha)}(\underline{\boldsymbol{x}}_n|\boldsymbol{\theta})
e^{-n\delta_n},
\end{eqnarray}

with $\lim_{n\rightarrow\infty}\delta_n = 0$ where $\widetilde{q}_n^{(\alpha)}(\cdot|\boldsymbol{\theta})$ 
and $\widetilde{\pi}_n^{(\alpha)}(\boldsymbol{\theta})$ 
are the $\alpha$-modified model and prior densities (see Remark \ref{REM:star}).
This definition follows from the fact that the $R^{(\alpha)}$-posterior density is proportional to 
$\widetilde{\pi}_n^{(\alpha)}(\boldsymbol{\theta}) \widetilde{q}_n^{(\alpha)}(\underline{\boldsymbol{x}}_n|\boldsymbol{\theta})$.
Note that, if the MRPE exists, then it is also an AMRPE.
Assume that this estimator  
$\widehat{\boldsymbol{\theta}}_\alpha=\widehat{\boldsymbol{\theta}}_\alpha(\underline{\boldsymbol{x}}_n)$,
as a function of data $\underline{\boldsymbol{x}}_n$, is measurable,
and consider such prior sequence that satisfies 
\begin{eqnarray}
~~~~~~~~~~\liminf_{n\rightarrow\infty} e^{nr}\widetilde{\pi}_n^{(\alpha)}(\boldsymbol{\theta}) \geq 1, ~~~~~
\mbox{ for all } r>0, ~\boldsymbol{\theta}\in \Theta.
\label{EQ:AsmpP1}
\end{eqnarray}

Assumption (\ref{EQ:AsmpP1}) signifies that the ($\alpha$-modified) prior probabilities 
are not exponentially small anywhere in  $\Theta$.
Then, we have the following theorems.

\begin{theorem}
	Consider the set-up of IID models with fixed countable $\Theta_n=\Theta$
	and discrete prior sequence $\pi_n$ satisfying Assumption (\ref{EQ:AsmpP1}).
	Suppose $\pi_n$ is information dense at the true probability mass function $g$ with respect to $\mathcal{F}_{\alpha}$ 
	as in Definition \ref{DEF:info_dense_iid} and 
	$\pi_n^{(\alpha)}\left(A_n^c|\underline{\boldsymbol{X}}_n\right)$ is exponentially small with probability one 
	for a sequence of measurable subsets $A_n \subseteq \Theta$.
	Then any AMRPE $ \widehat{\boldsymbol{\theta}}_\alpha \in A_n$ for all sufficiently large $n$ with probability one. 
	\label{THM:Consistency_AMRPE}
\end{theorem}

\begin{theorem}
	Consider the set-up of stationary independent models with fixed countable $\Theta_n=\Theta$
	and a discrete prior sequence $\pi_n$ satisfying Assumption (\ref{EQ:AsmpP1}).
	Then, for any true density $g$ which is an information limit of the (countable) family 
	$\left\{\widetilde{q}^{(\alpha)}(\cdot|\boldsymbol{\theta}) : \boldsymbol{\theta}\in \Theta_n \right\}$
	and for  any $\epsilon>0$, we have
	$
	\pi_n^{(\alpha)}\left( \left\{\boldsymbol{\theta} : d_1(g, f_{\boldsymbol{\theta}}) \geq \epsilon \right\} 
	|\underline{\boldsymbol{X}}_n\right)$ is exponentially small with probability one.
	So
	$
	\lim\limits_{n\rightarrow\infty} d_1(g, f_{\widehat{\boldsymbol{\theta}}_\alpha} )=0, 
	$
	with probability one for any AMRPE $ \widehat{\boldsymbol{\theta}}_\alpha$.
	\label{THM:Consistency_AMRPD2}
\end{theorem}

\begin{remark}
	Theorem \ref{THM:Consistency_AMRPD2}, in a special case $\alpha=0$, yields 	a stronger version of 
Theorem 15 of \cite{Barron:1988}. Our result requires fewer assumptions  than required by Barron's result.
	\label{REM:Asymp_AMRPD}
\end{remark}

\section{Application (II): Independent Non-homogeneous Models}
\label{SEC:Asymp_INH}

\subsection{Convergences of $R^{(\alpha)}$-Posterior and $R^{(\alpha)}$-Bayes estimators}
\label{SEC:Asymp_INH_post}

Let us now consider the set-up of independent but non-homogeneous (INH) models 
as described in Example \ref{SEC:Asymp_Post_setup}.2 of Section \ref{SEC:Asymp_Post_setup},
and simplify the exponential convergence results for the $R^{(\alpha)}$-posterior probabilities under this INH set-up.
Note that, in this case, $q_n^{(\alpha)}(\underline{\boldsymbol{x}}_n|\boldsymbol{\theta}) 
= \sum_{i=1}^n q_{i,\boldsymbol{\theta}}^{(\alpha)}(x_i)$ for any observed data $\underline{\boldsymbol{x}}_n=(x_1, \ldots, x_n)$,
and hence  $Q_n^{(\alpha)}(S_n|\boldsymbol{\theta}) = \prod_{i=1}^{n} Q^{(i,\alpha)}(S^i|\boldsymbol{\theta})$
for any $S_n = S^1\times S^2 \times \cdots \times S^n \in \mathcal{B}_n$ with $S^i\in \mathcal{B}^i$ for all $i$ 
and $Q^{(i,\alpha)}(S^i|\boldsymbol{\theta}) = \int_{S^i} \exp(q_{i,\boldsymbol{\theta}}^{(\alpha)}(y)) dy$.
Assume that $\Theta_n=\Theta$ and $\pi_n=\pi$ are independent of $n$.
Then, we have 
$$\widetilde{q}_n^{(\alpha)}(\underline{\boldsymbol{x}}_n|\boldsymbol{\theta}) = 
\frac{\prod_{i=1}^n\exp\left({q}_{i,\boldsymbol{\theta}}^{(\alpha)}(x_i)\right)}{
	Q_n^{(\alpha)}(\boldsymbol{\chi}_n|\boldsymbol{\theta})} = \prod_{i=1}^n\widetilde{q}^{(i,\alpha)}(x_i|\boldsymbol{\theta})$$
with $ \widetilde{q}^{(i,\alpha)}(x_i|\boldsymbol{\theta})
= \frac{\exp\left({q}_{i,\boldsymbol{\theta}}^{(\alpha)}(x_i)\right)}{Q^{(i,\alpha)}(\chi^i|\boldsymbol{\theta})}$.
Thus, in the notation of Section \ref{SEC:Merging}, we have 
$$D_n^{(\alpha)}(\boldsymbol{\theta})= \frac{1}{n} \sum_{i=1}^n 
KLD\left(g_i, \widetilde{q}^{(i,\alpha)}(\cdot|\boldsymbol{\theta})\right),$$
and hence the definition of information denseness can be simplified for the INH models as follows.

\begin{definition}
The prior $\pi$ under the INH model is said to be information dense at  $\boldsymbol{G}_n=(G_1, \ldots, G_n)$ 
with respect to $\mathcal{F}_{n,\alpha}=\otimes_{i=1}^n \mathcal{F}^i_\alpha$, 
if  there exists a finite measure $\widetilde{\pi}$ satisfying (\ref{EQ:infodense_prior}) such that\\
$\widetilde{\pi}\left(\left\{\boldsymbol{\theta} : \limsup_{n\rightarrow\infty}\frac{1}{n}\sum_{i=1}^n 
KLD\left(g_i,\widetilde{q}^{(i,\alpha)}(\cdot|\boldsymbol{\theta})\right) < \epsilon \right\}\right)> 0$,
for all $\epsilon>0$.
\label{DEF:info_dense_inh}
\end{definition} 

When $f_{i,\boldsymbol{\theta}}=f_{\boldsymbol{\theta}}$ is independent of $i$, 
then the INH set-up coincides with the IID set-up and the information denseness in Definition \ref{DEF:info_dense_inh}
coincides with that in Definition \ref{DEF:info_dense_iid}.
Further, Definition \ref{DEF:info_dense_inh} is also equivalent to the general Definition \ref{DEF:info_dense}
and hence implies that $G^n$ and $M_n^{(\alpha)}$ merge in probability.
Then, we have the following simplified  results for the INH set-up.

\begin{proposition}
Consider the set-up of INH models with $\Theta_n=\Theta$ and assume that the prior $\pi$ is independent of $n$ 
and information dense at $\boldsymbol{G}_n$ with respect to $\mathcal{F}_{n,\alpha}$ as per Definition \ref{DEF:info_dense_inh}.
Then, for any sequence of measurable parameter sets $A_n \subset \Theta$,
$\pi_n^{(\alpha)}\left(A_n^c|\underline{\boldsymbol{X}}_n\right)$ is exponentially small with $P$-probability one,
if and only if there exists sequences of measurable parameter sets $B_n, C_n \subset \Theta$
such that $A_n \cup B_n \cup C_n =\Theta$, 
$\frac{M_n^{(\alpha)}(\boldsymbol{\chi}_n, B_n)}{M_n^{(\alpha)}(\boldsymbol{\chi}_n, \Theta_n)}\leq e^{-nr}$
for  $r>0$ and a UEC test for $G^n$ against 
$\left\{ {Q_n^{(\alpha)}(\cdot|\boldsymbol{\theta})}/{Q_n^{(\alpha)}(\chi_n|\boldsymbol{\theta})} : \boldsymbol{\theta} \in C_n \right\}$
exists.
\label{PROP:PostConv_inh1}
\end{proposition}

However, the existence of the required UEC in Proposition \ref{PROP:PostConv_inh1} is equivalent to 
the existence of a UEC test for $G_i$ against 
$\left\{\frac{Q^{(i,\alpha)}(\cdot|\boldsymbol{\theta})}{Q^{(i,\alpha)}(\chi^i|\boldsymbol{\theta})} 
: \boldsymbol{\theta} \in C_n \right\}$
uniformly over $i=1,\ldots,n$. Following 
the discussions of Section \ref{SEC:AsympR_IID}, 
this holds if Assumption (B) is satisfied for $\widetilde{B}_{\epsilon}^{T_n} 
= \left\{\boldsymbol{\theta} : \frac{1}{n}\sum_{i=1}^{n}
d_1(f_{i,\boldsymbol{\theta}}, \widetilde{q}^{(i,\alpha)}(\cdot|\boldsymbol{\theta})^{T_n}) > \epsilon \right\}$
in place of ${B}_{\epsilon}^{T_n}$.
This leads to following simplification.

\begin{theorem}
Consider the INH models with $\Theta_n=\Theta$ and assume that the prior $\pi$ is independent of $n$ 
and information dense at $G_n$ with respect to $\mathcal{F}_{n,\alpha}$ as per Definition \ref{DEF:info_dense_inh}.
If Assumption (B) holds for $\widetilde{B}_{\epsilon}^{T_n}$ in place of ${B}_{\epsilon}^{T_n}$ for every $\epsilon>0$, 
the $R^{(\alpha)}$-posterior probability
$\pi_n^{(\alpha)}\left( \left\{\boldsymbol{\theta} : \frac{1}{n}\sum_{i=1}^{n} d_1(g_i, f_{i,\boldsymbol{\theta}}) \geq \epsilon \right\}	
|\underline{\boldsymbol{x}}_n\right)$ 
is exponentially small with $P$-probability one for $\epsilon>0$.
\label{PROP:PostConv_inh2}
\end{theorem}

%

We note that the Bernstein-von Mises type asymptotic results for the $R^{(\alpha)}$-posterior 
distribution under the INH set-up would be extremely important to provide contraction rates for 
our new robust pseudo-posterior; similar results for IID models were discussed in \cite{Ghosh/Basu:2016}.
However, considering the length of the present paper and to keep its focus clear on the exponential convergence results,
we propose to present the results on contraction rates for INH models
in a sequel paper; for the time being, they are made available in the ArXiv version  \citep{Majumder/etc:2019}.

\subsection{Robust Bayes Estimation under Fixed Design Regression Models}
\label{SEC:reg}

As noted in Example \ref{SEC:Asymp_Post_setup}.2,
the most common example of the general INH set-up is the fixed design regression models.
We consider the important example of model (\ref{EQ:Reg_Gen}) 
with $n$ fixed $k$-variate design points $\boldsymbol{t}_1, \ldots, \boldsymbol{t}_n$ and 
$f_{i, \boldsymbol{\theta}}(x) = \frac{1}{\sigma}f\left(\frac{x - \psi(\boldsymbol{t}_i, \boldsymbol{\beta})}{\sigma}\right)$
for some univariate density $f$.
The corresponding $\alpha$-likelihood is given by 
$q_n^{(\alpha)}(\underline{\boldsymbol{x}}_n|(\boldsymbol{\beta},\sigma)) 
= \sum_{i=1}^n q_{i,(\boldsymbol{\beta}, \sigma)}^{(\alpha)}(x_i)$ with 
$$q_{i,(\boldsymbol{\beta}, \sigma)}^{(\alpha)}(x_i)
= \frac{1}{\alpha\sigma^\alpha} f\left(\frac{x_i - \psi(\boldsymbol{t}_i, \boldsymbol{\beta})}{\sigma}\right)^\alpha
- \frac{M_{f,\alpha}}{(1+\alpha)\sigma^{\alpha}} - \frac{1}{\alpha},$$
where $M_{f,\alpha} = \int f^{1+\alpha}$.
Consider a prior density $\pi(\boldsymbol{\beta}, \sigma)$ 
for the parameters $(\boldsymbol{\beta}, \sigma)$ over the space $\Theta=\mathbb{R}^k\times(0,\infty)$ [$p=k+1$].
This prior can be chosen to be the conjugate prior or any subjective or objective prior; a common objective prior 
is the Jeffrey's prior given by $\pi(\boldsymbol{\beta}, \sigma) ={\sigma}^{-1}$.
Then, the $R^{(\alpha)}$-posterior density of $(\boldsymbol{\beta}, \sigma)$ is given by (\ref{EQ:R_post_densityGen})
which now simplifies as 
\begin{equation}
\pi_n^{(\alpha)}\left((\boldsymbol{\beta}, \sigma)|\underline{\boldsymbol{x}}_n\right) 
= \frac{ \prod_{i=1}^n \exp\left[\frac{1}{\alpha\sigma^\alpha} 
	f\left(\frac{x_i - \psi(\boldsymbol{t}_i, \boldsymbol{\beta})}{\sigma}\right)^\alpha
	- \frac{M_{f,\alpha}}{(1+\alpha)\sigma^{\alpha}}\right]\pi(\boldsymbol{\beta}, \sigma)}{
	\int\int \prod_{i=1}^n \exp\left[\frac{1}{\alpha\sigma^\alpha} 
	f\left(\frac{x_i - \psi(\boldsymbol{t}_i, \boldsymbol{\beta})}{\sigma}\right)^\alpha
	- \frac{M_{f,\alpha}}{(1+\alpha)\sigma^{\alpha}}\right]\pi(\boldsymbol{\beta}, \sigma) d\boldsymbol{\beta}d\sigma}.
\label{EQ:R_post_densityReg}
\end{equation}  
If $\sigma$ is known as in the Poisson or logistic regression models
(or can be assumed to be known with properly scaled variables), 
we consider a prior only on $\boldsymbol{\beta}$ given by, say, $\pi(\boldsymbol{\beta})$
which is either the objective uniform prior or the conjugate prior or some other proper prior.
In such cases, we can get the simplified form for the $R^{(\alpha)}$-posterior density of $\boldsymbol{\beta}$ as given by
\begin{equation}
\pi_n^{(\alpha)}\left(\boldsymbol{\beta}|\underline{\boldsymbol{x}}_n\right) 
= \frac{ \prod_{i=1}^n \exp\left[\frac{1}{\alpha\sigma^\alpha} 
	f\left(\frac{x_i - \psi(\boldsymbol{t}_i, \boldsymbol{\beta})}{\sigma}\right)^\alpha\right]\pi(\boldsymbol{\beta})}{
	\int \prod_{i=1}^n \exp\left[\frac{1}{\alpha\sigma^\alpha} 
	f\left(\frac{x_i - \psi(\boldsymbol{t}_i, \boldsymbol{\beta})}{\sigma}\right)^\alpha\right]\pi(\boldsymbol{\beta}) d\boldsymbol{\beta}}.
\label{EQ:R_post_densityReg2}
\end{equation}  
One can obtain the $R^{(\alpha)}$-Bayes estimators of $\boldsymbol{\beta}$, $\sigma$
under any suitable loss. 
We now study the exponential convergence for some regression examples
providing  simplifications for the required assumptions.

\subsection{Example: Normal Linear Regression Model with known variance}
\label{SEC:LRM1}

We consider the normal regression model, 
a particular member of the class of regression models considered in Section \ref{SEC:reg}, 
where $\psi(\boldsymbol{t}_i, \boldsymbol{\beta})=\boldsymbol{t}_i^T\boldsymbol{\beta}$
with $f$ being a standard normal density.
For simplicity, here we assume that the error variance $\sigma$ is known;
the unknown $\sigma$ case is considered later.
In this case, we can simplify the $R^{(\alpha)}$-posterior from (\ref{EQ:R_post_densityReg2})
and compute the expected $R^{(\alpha)}$-posterior estimator (ERPE) of $\boldsymbol{\beta}$;
however the resulting  $R^{(\alpha)}$-posterior has no explicit form and 
hence the corresponding ERPE needs to be computed numerically (see Sections \ref{SEC:simulation}, \ref{SEC:Computation}).

However, being a particular case of the INH set-up,  the exponential consistency of 
the $R^{(\alpha)}$-posterior of $\boldsymbol{\beta}$ directly holds under the assumptions of  
Proposition \ref{PROP:PostConv_inh1}. 
We now verify the required conditions for this present case normal linear regression models
with known $\sigma$. For this purpose, let us denote 
$\boldsymbol{D}=[\boldsymbol{t}_1, \ldots, \boldsymbol{t}_n]^T$, the fixed-design matrix, and 
$\boldsymbol{x}=(x_1, \ldots, x_n)^T$. 
Recall that, provided $\boldsymbol{D}$ has full column rank, 
the ordinary least square estimate of $\boldsymbol{\beta}$ is
$\widehat{\boldsymbol{\beta}}=(\boldsymbol{D}^T\boldsymbol{D})^{-1}\boldsymbol{D}^T\boldsymbol{x}$,
which is also the ordinary Bayes estimator under the uniform prior and 
has the variance ${n}^{-1}(\boldsymbol{D}^T\boldsymbol{D})^{-1}$.
We assume the following intuitive assumptions on the fixed design matrix $\boldsymbol{D}$ of the linear regression models.
\begin{itemize}
	\item[(R1)] The design points $\boldsymbol{t}_i=(t_{i1}, \ldots, t_{ik})^T$, $i=1, \ldots, n$, 
	are such that, for all $ j, l, s =1, \ldots, k$,
	\begin{equation}
	\displaystyle{\sup _{n >1 }} \displaystyle{\max _{1 \leq i \leq n }} |t_{ij}|=O(1), 
	~~~~ \displaystyle{\max _{1 \leq i \leq n }} |t_{ij}||t_{il}|=O(1),
~~~~
	\dfrac{1}{n}\sum_{i=1}^{n}|t_{ij}t_{il}t_{is}|=O(1).
	\end{equation}
	\item[(R2)] The matrix $\boldsymbol{D}$ satisfies
	$
	\displaystyle{\inf_{n }}\hspace*{0.1 in}[\text{min eigenvalue of } n^{-1}(\boldsymbol{D}^T\boldsymbol{D})] > 0,
	$
	which also implies the matrix $\boldsymbol{D}$ has full column rank, and
	$
	\displaystyle{\max _{1 \leq i \leq n }} [\boldsymbol{t}_i^T(\boldsymbol{D}^T\boldsymbol{D})^{-1}\boldsymbol{t}_i] 
	=O(n^{-1}).
	$
\end{itemize}

Note that these Assumptions (R1)--(R2) imply the (weak) consistency of 
the corresponding (frequentist) MDPDE of $\boldsymbol{\beta}$
obtained by minimizing the negative of the associated $\alpha$-likelihood function \citep{Ghosh/Basu:2013}.
They are easy to verify for any given design matrix; 
in particular they hold if $\boldsymbol{t}_i$'s are generated from some non-singular $k$-variate distributions.
It is shown in \cite{Majumder/etc:2019} that these two conditions indeed ensure 
a Bernstein-von Mises type result for the associated $R^{(\alpha)}$-posterior. 

It is really fascinating to see that, despite the complex natures of our earlier assumptions for general INH models,
these two simple Assumptions (R1)--(R2) imply the exponential consistency of the $R^{(\alpha)}$-posterior probability
at any $\alpha\geq 0$ for the example of linear regression (along with some mild conditions on the prior). 
The result is presented in the following theorem.

\begin{theorem}\label{THM:exp_linreg}
Consider the normal linear regression set up with known error variance. 
Assume that the true  parameter value is $\boldsymbol{\beta}_0$, i.e.,  $g_i=f_{i,\boldsymbol{\beta}_0}$ for all $i$,
and the prior on $\boldsymbol{\beta}$ is continuous and positive at $\boldsymbol{\beta}_0$. Take any $\alpha\geq 0$.
Then, under Assumptions  (R1)--(R2), given any $\epsilon>0$, there exists $r>0$ such that
$$\lim\limits_{n\rightarrow \infty}P\left[\pi_n^{(\alpha)}\left(\left\{\boldsymbol{\beta}: 
\frac{1}{n}\sum_{i=1}^{n}d_1(g_i,f_{i,\boldsymbol{\beta}}) \geq  \epsilon\right\}
\bigg|\underline{\boldsymbol{x}}_n\right)< e^{-nr}\right] =1,
$$
or equivalently, 
$$
~~~~~\lim\limits_{n\rightarrow \infty}P\left[\pi_n^{(\alpha)}\left(\left\{\boldsymbol{\beta}: 
\frac{1}{n}\sum_{i=1}^{n}\boldsymbol{t}_i^T|\boldsymbol{\beta}-\boldsymbol{\beta}_0| \geq  \epsilon\right\}
\bigg|\underline{\boldsymbol{x}}_n\right)< e^{-nr}\right] =1,
$$
i.e., the $R^{(\alpha)}$-posterior probabilities asymptotically concentrates 
on the neighborhoods of the true regression line at a exponential rate of convergence.
\end{theorem}


\subsection{Example: Normal Linear Regression Model with unknown variance}
\label{SEC:LRM2}

We now consider an extended version of the previous example of normal linear regression
with unknown error variance. Consider the set-up and notation of the previous subsection
with $\psi(\boldsymbol{t}_i, \boldsymbol{\beta})=\boldsymbol{t}_i^T\boldsymbol{\beta}$
and $f$ being a normal density with mean 0 and variance $\sigma$; 
but now we consider $\sigma^2$ to be also an unknown parameter along with the regression coefficient $\boldsymbol{\beta}$.
Given a prior $\pi(\boldsymbol{\beta}, \sigma)$ in this case, 
the $R^{(\alpha)}$-posterior distribution is given by (\ref{EQ:R_post_densityReg})
with $M_f=(2\pi)^{-\alpha/2}(1+\alpha)^{-1/2}$.

In this case as well, we have simplified the required conditions for 
the exponential convergence of the $R^{(\alpha)}$-posterior probabilities,
which is presented in the following theorem; 
interestingly, the same sets of conditions as in the known $\sigma$ case suffice.

\begin{theorem}\label{THM:exp_linreg2}
Consider the normal linear regression set up with unknown error variance. 
Assume that the true  parameter value is $\boldsymbol{\theta}_0 = (\boldsymbol{\beta}_0, \sigma_0^2)$, 
i.e.,  $g_i=f_{i,\boldsymbol{\theta}_0}$ for all $i$,
and the prior on $\boldsymbol{\theta}$ is continuous and positive at $\boldsymbol{\theta}_0$. 
Take any $\alpha\geq 0$.
Then, under Assumptions  (R1)--(R2), given any $\epsilon>0$, there exists $r>0$ such that
$$
\lim\limits_{n\rightarrow \infty}P\left[\pi_n^{(\alpha)}\left(\left\{\boldsymbol{\theta}: 
\frac{1}{n}\sum_{i=1}^{n}d_1(g_i,f_{i,\boldsymbol{\theta}}) \geq  \epsilon\right\}
	\bigg|\underline{\boldsymbol{x}}_n\right)< e^{-nr}\right] =1.
	$$
\end{theorem}

\subsection{Example: Logistic Regression Model}
\label{SEC:Logistic}

We now consider the important logistic regression model, 
which does not belong to the class of location-scale type regressions in Section \ref{SEC:reg}.
In the notation of Example \ref{SEC:Asymp_Post_setup}.2,
given fixed-design points $\boldsymbol{t}_1, \ldots, \boldsymbol{t}_n$, 
the logistic regression model considers binary response variables $x_i$, respectively,  
having  Bernoulli distribution with expectation 
$\psi(\boldsymbol{t}_i, \boldsymbol{\beta}) =\dfrac{e^{\boldsymbol{t}_i^{T}\beta}}{1+e^{\boldsymbol{t}_i^{T}\beta}}$,
for $i=1, \ldots, n$. 
As in Example \ref{SEC:Asymp_Post_setup}.2, this 
model clearly belongs 
to the INH set-up with the only parameter being the regression coefficient $\boldsymbol{\theta}=\boldsymbol{\beta}$; 
there is no scale parameter here. Thus, the $\alpha$-likelihood
$q_n^{(\alpha)}(\underline{\boldsymbol{x}}_n|\boldsymbol{\beta})$ of $\boldsymbol{\beta}$ 
is given by (\ref{EQ:alpha-likelihood_INH})
with $f_{i, \boldsymbol{\theta}}$ being the probability mass function of 
Bernoulli($\psi(\boldsymbol{t}_i, \boldsymbol{\beta})$) distribution and the integral being 
the sum over its support $\chi^i = \{0,1\}$; the underlying measure is the counting measure.
The $R^{(\alpha)}$ is obtained by using (\ref{EQ:R_post_densityGen}) given any prior $\pi(\boldsymbol{\beta})$,
which does not have a closed form and needs to be computed numerically;
see Section \ref{SEC:simulation} for illustrations. 

Let us now simplify the conditions required for the exponential consistency of the $R^{(\alpha)}$-posterior
for the logistic regression model. For this purpose, we recall the Assumption (R1) on the fixed design points
and consider the new condition (R3) in terms of the matrix 
$
\boldsymbol{\Psi}_n(\boldsymbol{\beta}) = n^{-1} E_{g_i}\left[
\frac{\partial^2}{\partial\boldsymbol{\beta}\partial\boldsymbol{\beta}^T}
q_n^{(\alpha)}(\underline{\boldsymbol{x}}_n|\boldsymbol{\beta}) \right].
$

(R3) $\displaystyle{\inf_{n }}~ [\text{min eigenvalue of }~ \boldsymbol{\Psi}_n(\boldsymbol{\beta})] > 0$, 
for all $\boldsymbol{\beta}$.

The matrix $\boldsymbol{\Psi}_n(\boldsymbol{\beta})$ appears in the asymptotic variance of the (frequentist) MDPDE
of $\boldsymbol{\beta}$ under the fixed-design logistic regression model \cite{Ghosh/Basu:2016b},
as well as in the Bernstein-von Mises type results for 
the corresponding $R^{(\alpha)}$-posterior distribution \cite{Majumder/etc:2019}.
Thus, in view of those results,  Assumption (R3) is extremely intuitive 
and easy to verify for any given design-matrix.
We have shown that, Assumptions (R1) and (R3) also imply the exponential convergence 
of our generalized $R^{(\alpha)}$-posterior probability in this logistic regression set-up,
as presented in the following theorem.

\begin{theorem}\label{exp_logistic}
Consider the fixed-design regression set up as above. 
Assume that the true  parameter value is $\boldsymbol{\beta}_0$, i.e.,  $g_i=f_{i,\boldsymbol{\beta}_0}$ for all $i$,
and the prior  $\pi(\boldsymbol{\beta})$ is continuous and positive at $\boldsymbol{\beta}_0$. 
Take any $\alpha\geq 0$.
Then, under Assumptions  (R1) and (R3), given any $\epsilon>0$, there exists $r>0$ such that
$$
\lim\limits_{n\rightarrow \infty}P\left[\pi_n^{(\alpha)}\left(\left\{\boldsymbol{\theta}: 
\frac{1}{n}\sum_{i=1}^{n}d_1(g_i,f_{i,\boldsymbol{\beta}}) \geq  \epsilon\right\}
\bigg|\underline{\boldsymbol{x}}_n\right)< e^{-nr}\right] =1.
$$
\end{theorem}

\section{Numerical Illustrations: Simulations}
\label{SEC:simulation}

\subsection{Performance of ERPE in Normal Linear Regression Model}
\label{SEC:Sim_LRM}

Let us now reconsider the regression model described in Sections \ref{SEC:LRM1}--\ref{SEC:LRM2},
and examine the finite sample performance of the expected $R^{(\alpha)}$-posterior estimator (ERPE) 
of the parameters.

We first assume that the error variance $\sigma$ is known and equals one.
The corresponding  $R^{(\alpha)}$-posterior is given by (\ref{EQ:R_post_densityReg2}), 
as discussed in Section \ref{SEC:LRM1}, and has no closed form solution. 
So, we have  computed the ERPE through an importance sampling Monte-Carlo.
We first simulate $n$ observations $t_{11}, \ldots, t_{1n}$ independently from $N(5,1)$
to fix the predictor values $\boldsymbol{t}_i = (1, t_{1i})^T$.
Then, $n$ independent error values $\epsilon_1, \ldots, \epsilon_n$ are generated from $N(0,1)$ (note $\sigma=1$)
and the responses are obtained through the linear regression structure 
$x_i = \boldsymbol{t}_i^T\boldsymbol{\beta} + \epsilon_i$ for $i=1, \ldots, n$,
with the true value of $\boldsymbol{\beta}$ being $\boldsymbol{\beta}_0 = (5, 2)^T$.
We have considered different sample sizes $n=20, 50, 100$, and 
different contamination proportions $\epsilon_C=$0\% (pure data), 5\%, 10\%, 20\%
to examine the finite sample robustness properties of our proposal. 
For contaminated samples, $[n\epsilon_C]$ error values are contaminated by generating them from $N(5,1)$ instead of $N(0, 1)$.
In each case, given a prior, the ERPE at different $\alpha\geq 0$ are
computed using 20000 steps in the importance sampling Monte-Carlo with the proposal density 
$N_k\left(\widehat{\boldsymbol{\beta}}, {n}^{-1}(\boldsymbol{D}^T\boldsymbol{D})^{-1}\right)$.
We replicate the above procedure 1000 times to compute the empirical bias and MSE of the ERPE 
for two different priors, namely the non-informative uniform prior 
and the conjugate normal prior, which are presented in the Online Supplement (Figures 1 and 2) due to page restriction.
The figures show that, under pure data, the bias and the MSE are the least 
for the usual Bayes estimator of $\boldsymbol{\beta}$ at $\alpha=0$, 
but their inflations are not very significant for the ERPEs with moderate $\alpha>0$.
Under contamination, the usual Bayes estimator (at $\alpha=0$)
has severely inflated bias and MSE and becomes highly unstable.
Our ERPEs with $\alpha>0$ are much more stable under contamination in terms of both bias and MSE;
the maximum  stability is observed for tuning parameters $\alpha \in [0.4, 0.6]$
yielding significantly improved robust Bayes estimators.

Next we consider the case of unknown error variance $\sigma$ in the above linear regression model,
as discussed in Section \ref{SEC:LRM2}. We repeat the above simulation exercise for the unknown $\sigma$ case as well,
by taking the true value of $\sigma_0=1$ and the conjugate prior on $(\boldsymbol{\beta}, \sigma)$ given by
$\pi(\boldsymbol{\beta},\sigma)=\pi(\boldsymbol{\beta}|\sigma)\pi(\sigma)$,
where $\pi(\beta|\sigma)$ is taken to be  $N_2(\boldsymbol{\beta}_0,\sigma^2 I_2)$ density
and $\pi(\sigma)$ is the density of the square root of Inverse chi-square distribution with 5 degrees of freedom
(i.e., prior for $\sigma^2$ is Inverse-$\chi^2_{5}$). 
However, in this case the computation of the ERPE could not be done efficiently using the simple importance sampling method
as in the case of known $\sigma$; 
alternatively we have used the  Metropolis-Hastings algorithm.

\noindent
\underline{\textbf{Algorithm 1: Computation of ERPE in LRM with unknown variance:}}\\
We generate 20000 sample observation from $R^{(\alpha)}$ posterior  distribution of 
$\boldsymbol{\theta}=(\boldsymbol{\beta},\sigma)$ as follows.
\begin{enumerate}
	\item[Step 1.] Start with $\boldsymbol{\theta}^{(0)}=(0,0,2)^T$. Set $k=1$.
	\item[Step 2.] After generating  $\boldsymbol{\theta}^{(k-1)}=(\boldsymbol{\beta}^{(k-1)}, \sigma^{(k-1)})$ 
	in the $(k-1)$-th step, at the $k^{th}$ step, generate 
	$\boldsymbol{\beta}^*$ and $\sigma^\ast$ from the proposal densities 
	$g_1 \equiv \mathcal{N}_2(\boldsymbol{\beta}^{(k-1)},I_2)$ and  $g_2\equiv$Exponential$(\sigma^{(k-1)})$, respectively.
	\item[Step 3.] Generate $U \sim U(0,1)$ and compute 
$
	\gamma=\frac{\exp[q_n^{(\alpha)}(\underline{\boldsymbol{x}}_n|\boldsymbol{\beta}^{\ast},\sigma^{\ast})
		g_1(\boldsymbol{\beta}^{\ast})g_2(\sigma^{\ast})}
	{\exp[q_n^{(\alpha)}(\underline{\boldsymbol{x}}_n|\boldsymbol{\beta}^{(k-1)},\sigma^{(k-1)})
		g_1(\boldsymbol{\beta}^{(k-1)})g_2(\sigma^{(k-1)})}.
$ 
	\item[Step 4.] If $U< \gamma$, set $\boldsymbol{\beta}^{(k)}=\boldsymbol{\beta}^*$ and $\sigma^{(k)}=\sigma^\ast$. 
	Otherwise, set $\boldsymbol{\beta}^{(k)}=\boldsymbol{\beta}^{(k-1)}$ and $\sigma^{(k)}=\sigma^{(k-1)}$.
	\item[Step 5.] Set $k=k+1$, and go to Step 2.
\end{enumerate}
In each cases, the first 5000 values generated are rejected as burn-in and 
the remaining 15000 parameter values are averaged to get a good approximation of the ERPE of $(\boldsymbol{\beta},\sigma)$. 
\hfill{$\square$}

\begin{figure}[!h]
	\centering
	\subfloat[$n=20$]{
		\begin{tabular}{c}
			\includegraphics[width=0.3\textwidth]{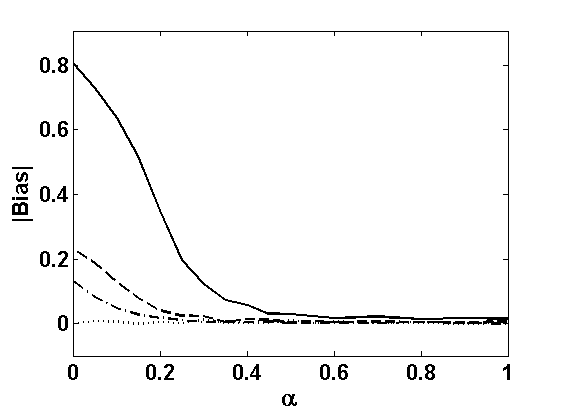}
			\\
			\includegraphics[width=0.3\textwidth]{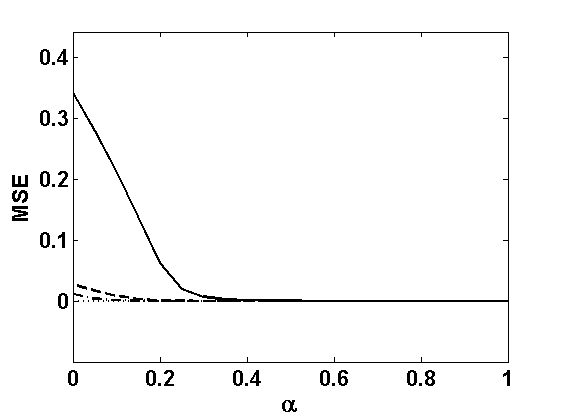}
		\end{tabular}				
		\label{FIG:LRM_n20CP}}
	\subfloat[$n=50$]{
		\begin{tabular}{c}
			\includegraphics[width=0.3\textwidth]{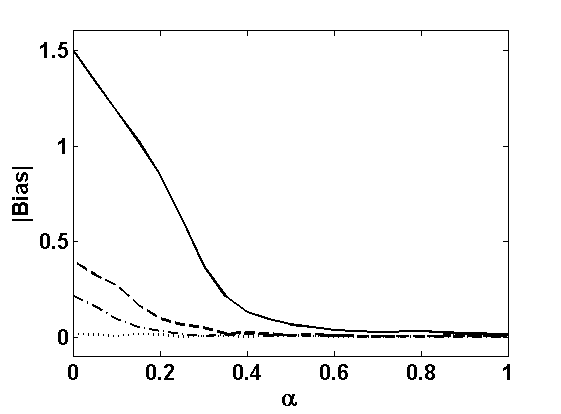}\\
			\includegraphics[width=0.3\textwidth]{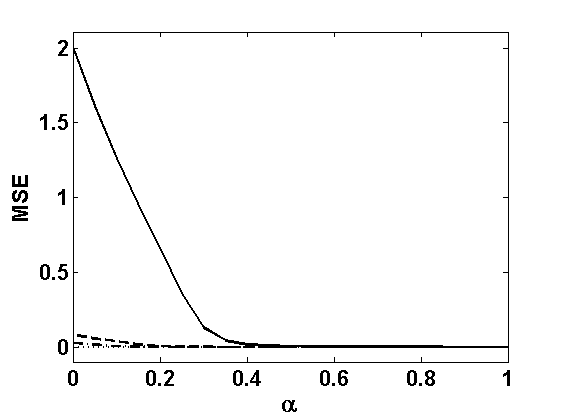}
		\end{tabular}			
		\label{FIG:LRM_n50cp}}
	\subfloat[$n=100$]{
		\begin{tabular}{c}		
			\includegraphics[width=0.3\textwidth]{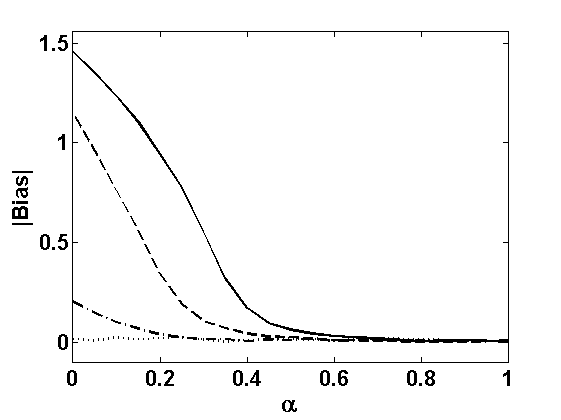}\\
			\includegraphics[width=0.3\textwidth]{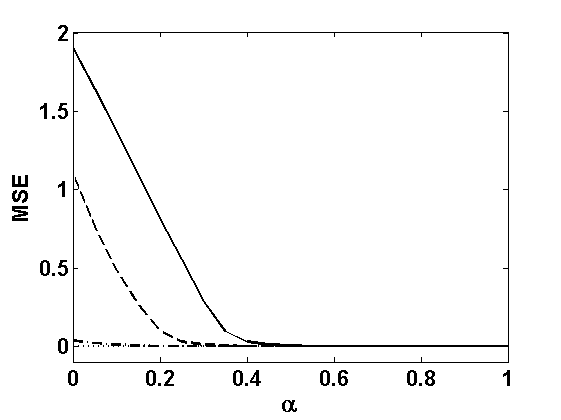}
		\end{tabular}	
		\label{FIG:LRM_n100cp}}
	\caption{Empirical total absolute bias and total MSE of the ERPE of $\boldsymbol{\beta}$ in the linear regression model
		with unknown $\sigma$ and the conjugate priors. 
		[Dotted line: $\epsilon_C=0\%$, Dash-Dotted line: $\epsilon_C=5\%$, Dashed line: $\epsilon_C=10\%$, 
		Solid line: $\epsilon_C=20\%$] (See additional discussions in the Online Supplement)}
	\label{FIG:MSE_normalRegCP1}
\end{figure}

\begin{figure}[!h]
	\centering
	\subfloat[$n=20$]{
\begin{tabular}{c}
\includegraphics[width=0.3\textwidth]{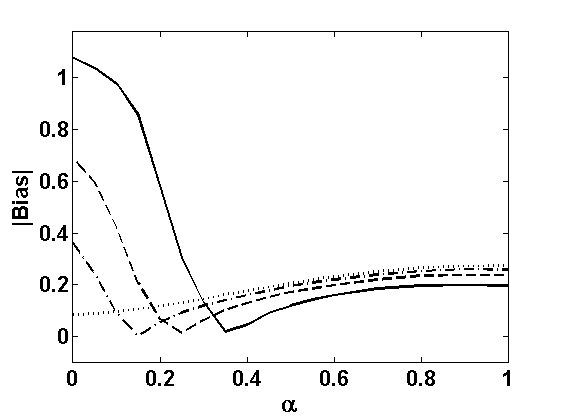}
		\\
		\includegraphics[width=0.3\textwidth]{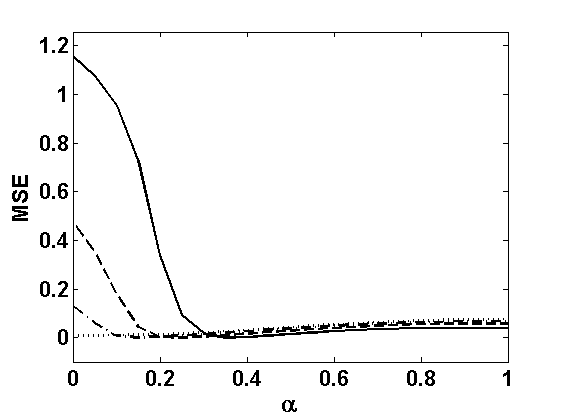}
\end{tabular}				
\label{FIG:LRM_n20CP}}
	\subfloat[$n=50$]{
\begin{tabular}{c}
		\includegraphics[width=0.3\textwidth]{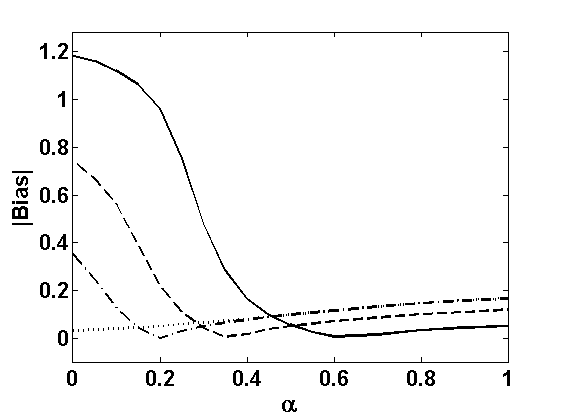}\\
		\includegraphics[width=0.3\textwidth]{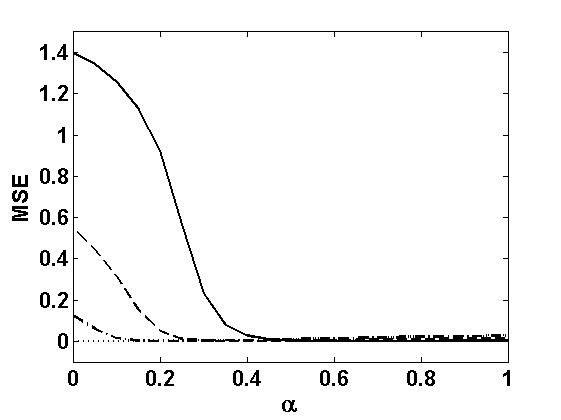}
\end{tabular}			
\label{FIG:LRM_n50cp}}
	\subfloat[$n=100$]{
\begin{tabular}{c}		\includegraphics[width=0.3\textwidth]{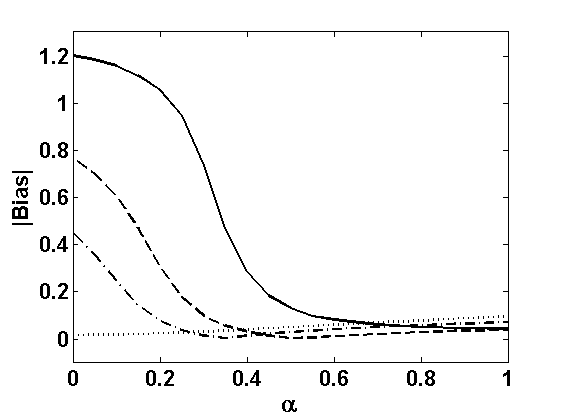}\\
		\includegraphics[width=0.3\textwidth]{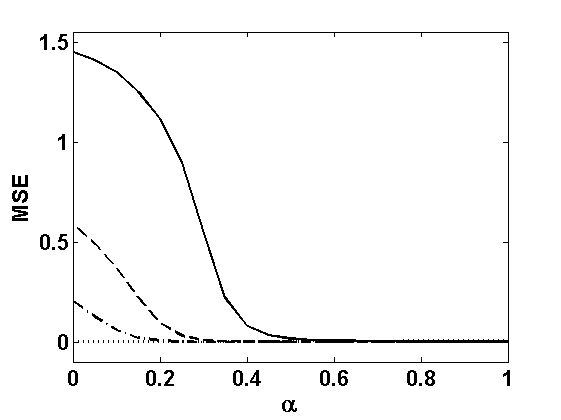}
\end{tabular}	
\label{FIG:LRM_n100cp}}
	\caption{Empirical absolute bias and MSE of the ERPE of $\sigma$ in the linear regression model
		with unknown $\sigma$ and the conjugate priors. 
		[Dotted line: $\epsilon_C=0\%$, Dash-Dotted line: $\epsilon_C=5\%$, Dashed line: $\epsilon_C=10\%$, 
		Solid line: $\epsilon_C=20\%$] (See additional discussions in the Online Supplement)}
	\label{FIG:MSE_normalRegCP1S}
\end{figure}

The process is replicated 1000 times to compute the empirical biases and MSEs of the ERPEs
of $\boldsymbol{\beta}$ and $\sigma$ at different $\alpha$ for the previous simulation set-up.
The resulting values of total absolute bias and the total MSE over the two components of $\boldsymbol{\beta}$ 
as well as the absolute bias and MSE of the ERPE of $\sigma$ are presented 
in Figures \ref{FIG:MSE_normalRegCP1} and \ref{FIG:MSE_normalRegCP1S}, respectively. 

The performances of the ERPE of regression coefficient and error variance 
are again the same as before in that the proposed ERPE with larger $\alpha$ provide extremely stable estimates
even under contamination up to 20\%. Under pure data the usual Bayes estimators give minimum absolute bias and MSEs,
but the ERPEs with $\alpha>0$ are also not very far away. 
However, under data contamination, the usual Bayes estimates (at $\alpha=0$) become extremely non-robust
yielding significantly higher bias and MSEs even though we are using strong conjugate prior. 
As the contamination proportion increases, we need larger values of $\alpha$ in the proposed ERPE 
to produce smaller biases and MSEs close to the pure data scenarios; 
in particular, $\alpha\geq 0.5$ always has excellent robust performance.


\subsection{Performance of ERPE in Logistic Regression Model}
\label{SEC:Sim_logistic}

We now consider the fixed-design logistic regression model as in Section \ref{SEC:Logistic}
and study the finite sample properties of the ERPE, 
the expectation of the regression coefficient $\boldsymbol{\beta}$ 
under the proposed $R^{(\alpha)}$-posterior distribution. 
Since the corresponding $R^{(\alpha)}$-posterior has no closed form solution, 
we have computed the ERPE numerically  in our simulation exercise. 

We first simulate $n$ values $t_{11}, \ldots, t_{1n}$ independently from $U(-5,5)$
and fix the design points as $\boldsymbol{t}_i = (1, t_{1i})^T$.
Then, the $n$ response values $x_1, \ldots, x_n$ are obtained through the logistic regression structure 
with $x_i$ generated from Bernoulli distribution with mean parameter
$\psi(\boldsymbol{t}_i, \boldsymbol{\beta}) =\dfrac{e^{\boldsymbol{t}_i^{T}\beta}}{1+e^{\boldsymbol{t}_i^{T}\beta}}$
for each $i=1, \ldots, n$; the true parameter value is taken as $\boldsymbol{\beta}_0 = (0, 5)^T$.
Again we have considered different sample sizes $n=20, 50, 100$, and 
different contamination proportions $\epsilon_C=$0\% (pure data), 5\%, 10\%, 20\%.
The contaminated observations, $[n\epsilon_C]$ many in a sample of size $n$, 
are  forced through misspecification of the response values, i.e., by changing $x_i$ to $(1-x_i)$,
and the prior is taken as the (bivariate) normal distribution as $\pi(\boldsymbol{\beta})\equiv N_2(\boldsymbol{\beta}_0, I_2)$.
However, in this case also, the importance sampling is seen to fail to provide a good approximation to the ERPE
and we have alternatively used the Metropolis-Hastings method. 
Note that, the target density, i.e $R^{(\alpha)}$ posterior density here is proportional to 
$g(\boldsymbol{\beta})=\exp[q_n^{(\alpha)}(\underline{\boldsymbol{x}}_n|\boldsymbol{\beta})]
\pi(\boldsymbol{\beta})d\boldsymbol{\beta}$.  

\noindent
\underline{\textbf{Algorithm 2: Computation of ERPE in logistic Regression:}}\\
We generate 20000 sample observation from $R^{(\alpha)}$ posterior  distribution of $\boldsymbol{\beta}$ as follows.
\begin{enumerate}
	\item[Step 1.] Start with $\boldsymbol{\beta}^{(0)}=(0,0)^T$.
	\item[Step 2.] After generating  $\boldsymbol{\beta}^{(k-1)}$ in the $(k-1)$-th step,
	at the $k^{th}$ step, generate $\boldsymbol{\beta}^*$ from 
	$\mathcal{N}_2(\boldsymbol{\beta}^{(k-1)},I_2)$. 
	\item[Step 3.] Generate $U \sim U(0,1)$ and compute $\gamma=g(\boldsymbol{\beta}^*)/g(\boldsymbol{\beta}^{(k-1)})$. 
	\item[Step 4.] If $U< \gamma$, set $\boldsymbol{\beta}^{(k)}=\boldsymbol{\beta}^*$. 
	Otherwise, set $\boldsymbol{\beta}^{(k)}=\boldsymbol{\beta}^{(k-1)}$.
	\item[Step 5.] Set $k=k+1$, and go to Step 2.
\end{enumerate}
In each case, the first 5000 values generated are rejected as burn-in and 
the remaining 15000 parameter values are averaged to get a good approximation of the ERPE. 
\hfill{$\square$}

\begin{figure}[!h]
	\centering
	\subfloat[$n=20$]{
\begin{tabular}{c}		
	\includegraphics[width=0.3\textwidth]{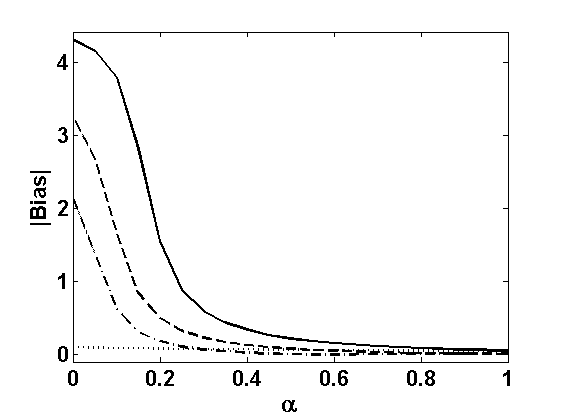}\\
		\includegraphics[width=0.3\textwidth]{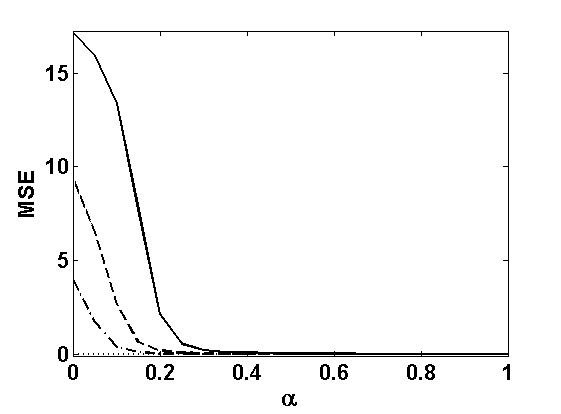}
\end{tabular}		
\label{FIG:Log_n20CP}}
	\subfloat[$n=50$]{
\begin{tabular}{c}		
	\includegraphics[width=0.3\textwidth]{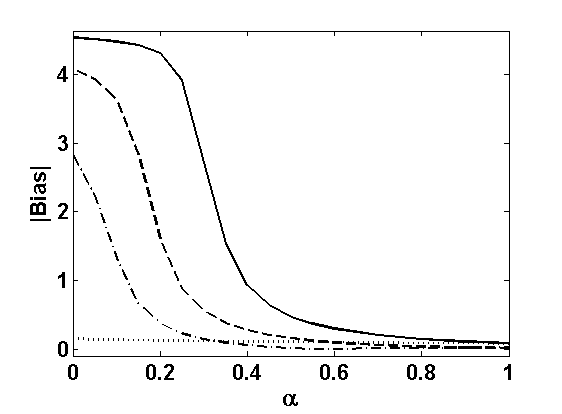}\\
		\includegraphics[width=0.3\textwidth]{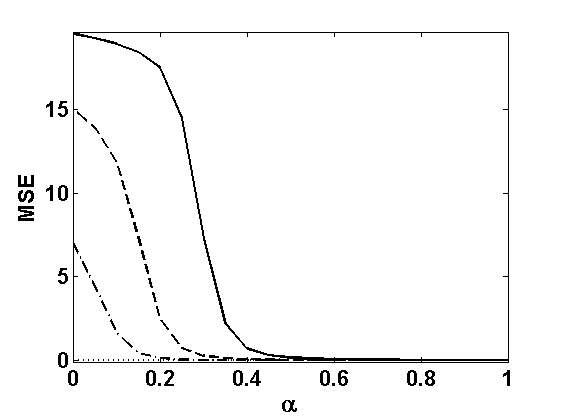}
\end{tabular}
		\label{FIG:Log_n50cp}}
	\subfloat[$n=100$]{
\begin{tabular}{c}	
	\includegraphics[width=0.3\textwidth]{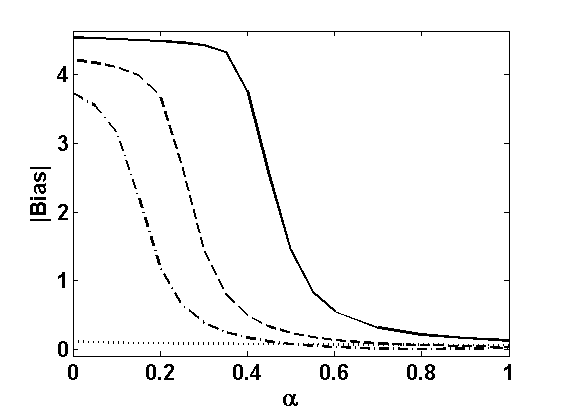}\\
		\includegraphics[width=0.3\textwidth]{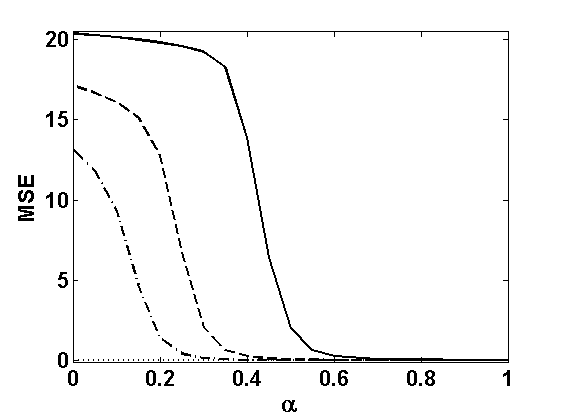}
\end{tabular}
		\label{FIG:Log_n100cp}}
	\caption{Empirical total absolute bias and total MSE of the ERPE of $\boldsymbol{\beta}$ in the logistic regression model
		with normal prior. 
		[Dotted line: $\epsilon_C=0\%$, Dash-Dotted line: $\epsilon_C=5\%$, Dashed line: $\epsilon_C=10\%$, 
		Solid line: $\epsilon_C=20\%$] (See additional discussions in the Online Supplement)}
	\label{FIG:MSE_logistic_normalP}
\end{figure}

The simulation exercise is replicated 1000 times to compute 1000 ERPEs of $\boldsymbol{\beta}$.
Their empirical biases and MSEs are presented in Figure \ref{FIG:MSE_logistic_normalP}.
Here also, 
it is clearly observed that the moderately larger values of $\alpha$
produce highly robust estimates under contaminations with only a slight loss in efficiency under pure data.
Under contamination, the MSE of the ERPEs remain stable for $\alpha\geq 0.5$; 
however, we need slightly larger $\alpha\geq 0.7$ to get smaller biases under heavy contamination of 20\%.

\section{Practical Aspects}
\label{SEC:Computation}

\subsection{On the Computation of the $R^{(\alpha)}$-Bayes Estimators}

A complex and challenging aspect of the proposed $R^{(\alpha)}$-Bayes estimators is their computation. 
This is, in fact, a common problem with all pseudo-posteriors 
that replace the likelihood  with some robust loss function.
In a frequentist sense, using a suitable optimization algorithm to derive a point estimator from some robust loss function results  
in scalable computation for many applications. 
In contrast, the computation of the whole pseudo-posterior is challenging for complicated models
and needs careful attention (even for the usual Bayes methods). 

For our $R^{(\alpha)}$-posterior also, no closed form expressions exist in most applications
and hence we need to compute the corresponding $R^{(\alpha)}$-Bayes estimators numerically.
One such possible approach could be the use of the importance sampling technique, 
which is seen to work well in our illustrations for 
normal means (\cite{Ghosh/Basu:2016}) or linear models with known $\sigma$ (Section \ref{SEC:Sim_LRM}).
But, this simple  approach can be useful only when it is possible to utilize some conjugacy structure;
in our cases, the standard posterior distribution is used as the proposal distribution due to their conjugacy. 
However, when the model is more complicated and we do not have a good proposal distribution, 
importance sampling fails to provide good approximations to the proposed $R^{(\alpha)}$-Bayes estimators;
this is because the $\alpha$-likelihood parts do not enjoy some conjugacy when the model is little bit more complicated, 
for example, the linear regression with unknown variance or the logistic regression models. 
In such cases, we propose to use a suitable Metropolis-Hastings algorithm
that is seen to work very well for the computations of the proposed ERPE 
under the above-mentioned two cases;
the corresponding algorithms are given in Sections \ref{SEC:Sim_LRM} and \ref{SEC:Sim_logistic}, respectively. 
We have also supplied the relevant R codes for the computations of the ERPEs for our examples in the Online Supplement.

We hope that, with the advance in modern computers, it would be possible to develop similar algorithms for 
the computation of the $R^{(\alpha)}$-posterior and the $R^{(\alpha)}$-Bayes estimators for other useful models.
However, if the model becomes too complex, the usual Bayes computation also becomes challenging
and we have to develop appropriate computation algorithms more carefully. 
An alternative approach can be to approximate the $R^{(\alpha)}$-Bayes estimators for larger sample sizes 
using asymptotic expansions like Laplace's one; such an approximation for our $R^{(\alpha)}$-posterior 
and its expectations are provided in \cite{Majumder/etc:2019} for general non-homogeneous (but independent) observations.  
These computational aspects of our robust pseudo-posterior would surely form a sequence of interesting future works.

\subsection{On the Choice of the tuning parameter $\alpha$}

We have proposed a class of robust pseudo-posteriors, indexed by the tuning parameter $\alpha>0$,
which coincides with the non-robust but (asymptotically) most efficient ordinary Bayes posterior as $\alpha\rightarrow0$.
In all our illustrations in Section \ref{SEC:simulation} it is observed that, with  increasing values of $\alpha>0$,
the asymptotic performance of the proposed $R^{(\alpha)}$-Bayes estimators deteriorate slightly  under pure data, 
but their robustness under data contamination improves significantly compared to the usual Bayes estimates (at $\alpha=0$).
Thus, a natural and practical question arises 
-- which $\alpha$ should one use for a given data set? 
As we have observed numerically that, with conjugate prior,
any $\alpha\geq 0.5$ provides extremely  robust inference under contamination,
whereas the empirically suggested range for the cases with uniform prior is $\alpha\in(0.4, 0.07)$;
thus, from our simulations presented here (along with numerous others not presented for brevity) $\alpha\approx0.5$
seems to be a good choice in most cases.

However, a more systematic procedure for selection of this tuning parameter depending on the given data at hand
would surely be useful for reliable applications of our proposal. 
In this regard, we note that the asymptotic distribution of the proposed ERPE at any $\alpha\geq 0$ 
is the same as that of the corresponding frequentist MDPDE for both IID and INH cases \cite{Ghosh/Basu:2016, Majumder/etc:2019}.
Therefore, finding the optimal tuning parameter for the ERPE becomes an asymptotically equivalent problem of choosing 
an $\alpha$ for the optimal control between robustness and efficiency of the MDPDE. 
The second one has received some attention  in the literature; 
one such approach chooses $\alpha$ by minimizing an asymptotic MSE of the MDPDE, 
with respect to $\alpha\in[0,1]$, given by
\begin{eqnarray}
\widehat{\mbox{AMSE}}(\alpha) = (\widehat{\boldsymbol{\theta}}_\alpha-\boldsymbol{\theta}^P)^T(\widehat{\boldsymbol{\theta}}_\alpha-\boldsymbol{\theta}^P)
+\frac{1}{n} \mbox{Trace}(\Sigma_\alpha(\widehat{\boldsymbol{\theta}}_\alpha)),
\label{EQ:AMSE}
\end{eqnarray}

where $\widehat{\boldsymbol{\theta}}_\alpha$ is the MDPDE at $\alpha$, 
$\Sigma_\alpha$ is the asymptotic variance of $\sqrt{n}\widehat{\boldsymbol{\theta}}_\alpha$
and $\boldsymbol{\theta}^P$ is some suitable pilot estimator.
The details can be found in \cite{Warwick/Jones:2005} and \cite{Ghosh/Basu:2015} for IID and INH set-ups, respectively,
where some suggestions regarding the choice of pilot $\boldsymbol{\theta}^P$ are also provided.

Since the asymptotic MSE of the MDPDE is indeed the same as the frequentist MSE of our ERPE,
the same process can be used to chose optimum $\alpha$ for the ERPE here when using improper non-informative priors
with $\widehat{\boldsymbol{\theta}}_\alpha$ being replaced by the corresponding ERPE, 
say $\widehat{\boldsymbol{\theta}}_\alpha^\ast$, at any given $\alpha$.
However, if we have a proper subjective prior, say $\pi(\boldsymbol{\theta})$, 
then we can improve this approach appropriately  by taking the pilot $\boldsymbol{\theta}^P$ 
as a random variable following $\pi(\boldsymbol{\theta})$ and then taking expected bias in (\ref{EQ:AMSE});
the modified criterion is  then given by
\begin{eqnarray}
\widehat{\mbox{AMSE}}^\ast(\alpha) = \int (\widehat{\boldsymbol{\theta}}_\alpha^\ast-\boldsymbol{\theta})^T
(\widehat{\boldsymbol{\theta}}_\alpha^\ast-\boldsymbol{\theta})\pi(\boldsymbol{\theta})d\boldsymbol{\theta}
+\frac{1}{n} \mbox{Trace}(\Sigma_\alpha(\widehat{\boldsymbol{\theta}}_\alpha^\ast)),
\label{EQ:AMSE_B}
\end{eqnarray}

which we can minimize with respect to $\alpha$, possibly through a grid search over $[0,1]$, 
to chose an appropriate tuning parameter value. 
However, this proposal clearly needs further detailed investigation which, considering the length of the current paper, 
we hope to do in a future work.

\section{Real Data Applications}
\label{SEC:data}

\subsection{Hertzsprung-Russell star cluster data}

As our first application, let us consider the famous star cluster (CYG OB1) data from the Hertzsprung-Russell diagram
about the logarithms of the light intensity ($L/L_0$) and the effective temperature ($T_e$) at the surface of 47 stars in the direction of Cygnus 
(Table 3, Chapter 2, \cite{Rousseeuw/Leroy:1987}). 
These data are studied by several authors (e.g., \cite{Rousseeuw/Leroy:1987,Ghosh/Basu:2013}) for demonstration of robust methods
through a simple linear regression with ($L/L_0$) being the response and $T_e$ as the covariate;
it has been observed there that four stars in the data (with indices 11, 20, 30 and 34) are indeed significantly different from the remaining stars
and produce the non-robust outlier effects while using classical estimation methods.

Here we have performed the Bayesian analyses of the simple linear regression model with different conjugate and improper priors.
As in the common practice, we assume the error variance $\sigma^2$ to be unknown.
For brevity, we present only the results for the extreme case of uniform  priors $\pi(\boldsymbol{\beta}, \sigma)=\sigma^{-1}$;
the resulting values of the ERPE (with and without the outliers) are presented in Table \ref{TAB:HR}.
It can be clearly observed that the usual Bayes estimates (at $\alpha=0$)
are extremely non-robust producing regression coefficients of opposite sign due to the presence of outliers.
however, our proposed $R^{(\alpha)}$-Bayes approach and the corresponding ERPEs remain extremely stable 
for moderately large values of $\alpha$ and successfully counter the effect of outliers. 

\begin{table}[h]
	\centering
	\caption{The ERPEs of the coefficients and error variance $\sigma^2$ in the simple linear regression models 
		for the Hertzsprung-Russell data with uniform prior.}
	\begin{tabular}{r|rrr|rrr} \hline
		& \multicolumn{3}{c|}{Original Data} & \multicolumn{3}{c}{Without Four Outliers}\\\noalign{\smallskip}\hline\noalign{\smallskip}
		$\alpha$ & Intercept &	Slope	&	$\sigma$	&	Intercept	&	Slope	&	$\sigma$ \\\noalign{\smallskip}\hline\noalign{\smallskip}
		0	&	7.33	&	$-$0.54	&	0.55	&	$-$3.38	&	1.89	&	0.41	\\
		0.1	&	6.83	&	$-$0.42	&	0.58	&	$-$4.90	&	2.24	&	0.42	\\
		0.25	&	$-$8.91	&	3.14	&	0.41	&	$-$5.78	&	2.43	&	0.41	\\
		0.4	&	$-$6.13	&	2.51	&	0.42	&	$-$8.73	&	3.10	&	0.39	\\
		0.5	&	$-$6.60	&	2.62	&	0.43	&	$-$7.75	&	2.88	&	0.38	\\
		0.6	&	$-$7.19	&	2.75	&	0.41	&	$-$9.68	&	3.31	&	0.39	\\
		0.8	&	$-$7.22	&	2.76	&	0.42	&	$-$7.76	&	2.88	&	0.42	\\
		\hline
	\end{tabular}
	\label{TAB:HR}
\end{table}

\subsection{Skin Data}

Let us now consider another popular example of logistic regression models having outlier issue,
namely a controlled study on the occurrence of  ``vaso constrictions" in the skin of digits due to air inspiration after a single deep breath  \cite{f47}. 
This Skin  dataset was analyzed by several authors including the recent  work by \cite{Ghosh/Basu:2016b}
where the logistic regression parameters are robustly estimated by the MDPDEs. 
Here the important covariates to model the vaso constriction occurrences are 
the logarithms of the volume of inspired air (``log.Vol") and the rate of inspiration (``log.Rate").
One can observe by plotting these data (see, for example, \cite{Ghosh/Basu:2016b}) that the $4$-th
and $18$-th observations are indeed the outliers making it difficult to separate the responses;
the MLE of the corresponding regression coefficients in the logistic regression model also changes significantly to have the values 
($-$2.88, 4.56, 5.18) in the presence of outliers and  ($-$24.58, 31.94, 39.55) after removal of the outliers. 

Here we have considered the Bayesian modeling of the same regression model with different types of priors. 
Again for brevity, we present only the case of uniform prior over the cube $[-50, 50]^3$ having the most extreme effect of outliers. 
The resulting ERPE for different values of $\alpha$ under the full data (including outliers)
as well as under the outlier deleted data are given in Table \ref{TAB:10est_skin_Wout};
note that the values corresponding to $\alpha=0$ gives the usual Bayes estimator (posterior mean). 
Clearly, the usual Bayes estimates get highly affected by the presence of only two outliers
whereas our $R^{(\alpha)}$-Bayes estimators, the ERPEs, with $\alpha$ around 0.5 provides extremely stable 
results even in the presence of outliers.

\begin{table}[h]
	\centering
	\caption{The ERPEs of the coefficients in a logistic regression for the Skin data with uniform prior.}
		\begin{tabular}{r|rrr|rrr} \hline
			& \multicolumn{3}{c|}{Original Data} & \multicolumn{3}{c}{Without Outliers ($4^{\rm th}$ and $18^{\rm th}$ obs.)}\\\noalign{\smallskip}\hline\noalign{\smallskip}
			$\alpha$ & Intercept &	log(Rate)	&	log(Vol)	&	Intercept	&	log(Rate)	&	log(Vol) \\\noalign{\smallskip}\hline\noalign{\smallskip}
0	&	$-$4.68	&	7.26	&	7.23	&	$-$22.35	&	35.17	&	29.58	\\
0.1	&	$-$5.73	&	9.02	&	8.46	&	$-$22.32	&	34.96	&	29.62	\\
0.25	&	$-$19.45	&	30.21	&	26.03	&	$-$22.53	&	34.91	&	30.02	\\
0.4	&	$-$22.38	&	34.15	&	29.94	&	$-$22.91	&	34.92	&	30.61	\\
0.5	&	$-$22.94	&	34.54	&	30.72	&	$-$23.18	&	34.88	&	31.02	\\
0.6	&	$-$23.29	&	34.61	&	31.20	&	$-$23.41	&	34.79	&	31.37	\\
0.8	&	$-$23.63	&	34.45	&	31.72	&	$-$23.71	&	34.54	&	31.80	\\
\hline
	\end{tabular}
	\label{TAB:10est_skin_Wout}
\end{table}

\section{Concluding Remark}
\label{SEC:Conclusion}

This paper presents a general Bayes pseudo-posterior under general parametric set-up
that produces pseudo-Bayes estimators which incorporate prior belief in the general spirit of Bayesian philosophy 
but are also robust against data contamination.
The exponential consistency of the proposed pseudo-posterior probabilities 
and the corresponding estimators are proved and illustrated for the cases of independent 
stationary and non-homogeneous models; separate attention is given to the case of discrete priors
with stationary models. Further applications of the proposed pseudo-Bayes estimators
are described in the context of 
linear and logistic regression models. 
All results of \cite{Barron:1988} turn out to be special cases of our results 
when the tuning parameter $\alpha$ is set to 0.

On the whole, we trust that this paper opens up a new and interesting area of research on robust hybrid inference
that has the flexibility to incorporate prior belief and inherits optimal properties from the Bayesian paradigm
along with the frequentists' robustness against data contamination
and hence could be very helpful in different complex practical problems.
In this sense, all Bayesian inference methodologies can be extended with this new pseudo-posterior.
In particular, a detailed study of the examples discussed in Section \ref{SEC:Asymp_Post_setup} 
should be an interesting future work for different applications. 
Extended versions of the Bayes testing and model selection criteria based on this new pseudo-posterior
can also be developed to achieve greater robustness for inference under data contamination. 


\newpage
\begin{center}
{\Large\bf Supplementary Material}
\end{center}

\appendix
\section{Proofs of the Results of Section 3 in Main paper}

\subsection{Proof of Theorem 3.1}

First note that, by help of Equation (6) of the main paper,
we can rewrite the $R^{(\alpha)}$-Bayes joint distribution as
\begin{eqnarray}
L_n^{(\alpha) Bayes}\left(d\boldsymbol{\theta}, d\underline{\boldsymbol{x}}_n\right) 
= \frac{M_n^{(\alpha)}(d\underline{\boldsymbol{x}}_n, d\boldsymbol{\theta})}{
	M_n^{(\alpha)}(\boldsymbol{\chi}_n, \Theta_n)}
=\widetilde{q}_n^{(\alpha)}(\underline{\boldsymbol{x}}_n|\boldsymbol{\theta})
\widetilde{\pi}_n^{(\alpha)}(d\boldsymbol{\theta})\lambda^n(d\underline{\boldsymbol{x}}_n),
\end{eqnarray}
which has a density  $\widetilde{q}_n^{(\alpha)}(\underline{\boldsymbol{x}}_n|\boldsymbol{\theta})$
with respect to $\widetilde{\pi}_n^{(\alpha)} \times \lambda^n$. On the other hand, the frequnetist approximation 
$L_n^{*(\alpha)}$ has the density function 
$e^{-nD_n^{(\alpha)}(\boldsymbol{\theta})}g_n\left(\underline{\boldsymbol{x}}_n\right)/c_n$
and hence we get
\begin{eqnarray}
KLD\left(L_n^{*(\alpha)}, L_n^{(\alpha)Bayes}\right)
&=& E_{L_n^{*(\alpha)}}\left[\log \frac{e^{-nD_n^{(\alpha)}(\boldsymbol{\theta})}g_n\left(\underline{\boldsymbol{X}}_n\right)/c_n}{
	\widetilde{q}_n^{(\alpha)}(\underline{\boldsymbol{X}}_n|\boldsymbol{\theta})}\right]\nonumber\\
&=& E_{\pi_n^{*(\alpha)}}E_{G_n}
\left[ -nD_n^{(\alpha)}(\boldsymbol{\theta}) + \log \frac{g_n\left(\underline{\boldsymbol{X}}_n\right)}{	
	\widetilde{q}_n^{(\alpha)}(\underline{\boldsymbol{X}}_n|\boldsymbol{\theta}) }- \log c_n\right]\nonumber\\
&=& E_{\pi_n^{*(\alpha)}}\left[ -nD_n^{(\alpha)}(\boldsymbol{\theta}) + E_{G_n} \log \frac{g_n\left(\underline{\boldsymbol{X}}_n\right)}{	
	\widetilde{q}_n^{(\alpha)}(\underline{\boldsymbol{X}}_n|\boldsymbol{\theta}) }\right]- \log c_n\nonumber\\
&=& E_{\pi_n^{*(\alpha)}}\left[ -nD_n^{(\alpha)}(\boldsymbol{\theta}) + nD_n^{(\alpha)}(\boldsymbol{\theta})\right]- \log c_n\nonumber\\
&=& - \log c_n 
= -\log \left[\int e^{-nD_n^{(\alpha)}(\boldsymbol{\theta})}\widetilde{\pi}_n^{(\alpha)}(d\boldsymbol{\theta})\right]
\nonumber
\end{eqnarray}

Therefore, for any $\epsilon>0$, we get 
\begin{eqnarray}
\frac{1}{n}KLD\left(L_n^{*(\alpha)}, L_n^{(\alpha)Bayes}\right)
&=& - \frac{1}{n}\log
\left[\int e^{-nD_n^{(\alpha)}(\boldsymbol{\theta})}\widetilde{\pi}_n^{(\alpha)}(d\boldsymbol{\theta})\right]
\nonumber\\
&\leq & \frac{\epsilon}{2} - \frac{1}{n} \log \widetilde{\pi}_n^{(\alpha)}\left(
\left\{\boldsymbol{\theta} : D_n^{(\alpha)}(\boldsymbol{\theta}) < \epsilon \right\}\right)
\nonumber\\
&\leq & \frac{\epsilon}{2} + \frac{\epsilon}{2} = \epsilon,  ~~~~\mbox{ for all but finitely many}~n,
\end{eqnarray}
by applying Assumption (M1) with $r=\frac{\epsilon}{2}$. 
Since $\epsilon>0$ is arbitrary,  this completes the proof of the first equation in Part (a) of the theorem.

The second equation of Part (a) and the first equation in Part (b) of the theorem follows by the relation 
$$
\frac{1}{n} KLD\left(L_n^{*(\alpha)}, L_n^{(\alpha)Bayes}\right) 
= \frac{1}{n} E_{G^n}\left[KLD\left(\pi_n^{*(\alpha)}(\cdot), \pi_n^{(\alpha)}(\cdot|\underline{\boldsymbol{X}}_n)\right)\right] 
+ \frac{1}{n} KLD\left(g^n, m_n^{(\alpha)}\right).
$$

\noindent
Finally, to proof the the last part of (b) in the theorem, we note that the Kullback-Leibler divergence satisfies the relation
$$
E\left|\log \frac{g_n(\underline{\boldsymbol{X}}_n)}{m_n^{(\alpha)}(\underline{\boldsymbol{X}}_n)}\right|
\leq  KLD\left(g^n, m_n^{(\alpha)}\right) + \frac{2}{e}.
$$
Therefore, by the first part of (b), we get 
$\lim\limits_{n\rightarrow\infty}  
E\left|\log \frac{g_n(\underline{\boldsymbol{X}}_n)}{m_n^{(\alpha)}(\underline{\boldsymbol{X}}_n)}\right|
=0$ and hence $G^n$ and $M_n^{(\alpha)}$ merge in probability by using Markov inequality.
\hfill{$\square$}

\subsection{Proof of Theorem 3.2}

To show Assumption (M1), let us fix $\epsilon, r >0$ and define
$\rho_n(\boldsymbol{\theta}) = e^{nr} \frac{d\pi_n}{d\widetilde{\pi}}(\boldsymbol{\theta})$.
Then, using Fatou's Lemma, we get
\begin{eqnarray}
\liminf_{n\rightarrow\infty} e^{nr} 
\pi_n\left(\left\{\boldsymbol{\theta} : D_n^{(\alpha)}(\boldsymbol{\theta}) < \epsilon \right\}\right)
&=& \liminf_{n\rightarrow\infty} e^{nr} 
\int I\left(\left\{\boldsymbol{\theta} : D_n^{(\alpha)}(\boldsymbol{\theta})<\epsilon \right\}\right){\pi}_n(d\boldsymbol{\theta})
\nonumber\\
&=&\liminf_{n\rightarrow\infty}\int I\left(\left\{\boldsymbol{\theta} : D_n^{(\alpha)}(\boldsymbol{\theta})<\epsilon \right\}\right)
\rho_n(\boldsymbol{\theta})\widetilde{\pi}(d\boldsymbol{\theta})
\nonumber\\
&=&\int \liminf_{n\rightarrow\infty}I\left(\left\{\boldsymbol{\theta} : D_n^{(\alpha)}(\boldsymbol{\theta})<\epsilon \right\}\right)
\rho_n(\boldsymbol{\theta})\widetilde{\pi}(d\boldsymbol{\theta})
\nonumber\\
&=&\int I\left(\left\{\boldsymbol{\theta} : \bar{D}^{(\alpha)}(\boldsymbol{\theta})<\epsilon \right\}\right)\widetilde{\pi}(d\boldsymbol{\theta})
\nonumber\\
&=& \widetilde{\pi}\left(\left\{\boldsymbol{\theta} : D_n^{(\alpha)}(\boldsymbol{\theta})<\epsilon \right\}\right),
\end{eqnarray}
which is strictly positive by the information denseness with respect to $\mathcal{F}_{n,\alpha}$ 
(Definition 3.2 of the main paper).
This implies Assumption (M1) and we are done in the view of Theorem 3.1.
\hfill{$\square$}

\subsection{Proof of Theorem 3.3}

We use an argument similar to that used by \cite{Barron:1988}.
Let us consider the following two assumptions in addition to Assumptions (A1)--(A3) and (A3)$^\ast$.
\begin{itemize}
	\item[(A4)] The true distribution $G^n$ and the $R^{(\alpha)}$-marginal distribution $M_n^{(\alpha)}$ satisfy 
	$$
	\lim\limits_{n\rightarrow\infty} P\left(\frac{m_n^{(\alpha)}(\underline{\boldsymbol{X}}_n)}{g^n(\underline{\boldsymbol{X}}_n)}
	\geq a_n \right) = 1.
	$$	
	
	\item[(A4)$^\ast$] The true distribution $G^n$ and the $R^{(\alpha)}$-marginal distribution $M_n^{(\alpha)}$ satisfy
	$$
	P\left(\frac{m_n^{(\alpha)}(\underline{\boldsymbol{X}}_n)}{g^n(\underline{\boldsymbol{X}}_n)}
	< a_n ~\mbox{ i.o.}\right) = 0.
	$$	
\end{itemize}

Note that, if Conditions (A4) and (A4)$^\ast$ hold with $a_n = e^{-n\epsilon}$ for every $\epsilon>0$, 
they indicate that the true distribution $G^n$ and the $R^{(\alpha)}$-marginal distribution $M_n^{(\alpha)}$ merge 
in probability or with probability one respectively.

\noindent
Now, we start with two primary results on the convergence of the $R^{(\alpha)}$-posterior probabilities.

\begin{lemma}
	\label{LEM:PostConv_Suff}
	Suppose Assumptions (A1)--(A3) and (A4) hold with $\lim b_n = \lim c_n = 0$ 
	such that $r_n :=(b_n + c_n)/a_n$ is finitely defined. 
	Then, for all $\delta>0$, we have 
	\begin{eqnarray}
	\limsup_{n\rightarrow\infty} P\left(\pi_n^{(\alpha)}\left(A_n^c|\underline{\boldsymbol{X}}_n\right)> \frac{r_n}{\delta} \right)
	\leq \delta.
	\label{EQ:Postconv_Suff1}
	\end{eqnarray}
	Further, if additionally Assumptions (A3)$^\ast$ and (A4)$^\ast$ are satisfied, 
	then for any summable sequence $\delta_n>0$ we have 
	\begin{eqnarray}
	P\left(\pi_n^{(\alpha)}\left(A_n^c|\underline{\boldsymbol{X}}_n\right)> \frac{r_n}{\delta_n} 
	\mbox{ i.o.}\right)=0.
	\label{EQ:Postconv_Suff2}
	\end{eqnarray}
\end{lemma}
\noindent\textbf{Proof:}
Note that, with $G^\infty$ probability one, the $R^{(\alpha)}$-posterior  probability can be re-expressed as
\begin{equation}
\pi_n^{(\alpha)}\left(A_n^c|\underline{\boldsymbol{X}}_n\right) 
= \frac{m_n^{(\alpha)}(\underline{\boldsymbol{X}}_n, A_n^c)/M_n^{(\alpha)}(\boldsymbol{\chi}_n, \Theta_n)g^n(\underline{\boldsymbol{X}}_n)}{
	m_n^{(\alpha)}(\underline{\boldsymbol{X}}_n)/g^n(\underline{\boldsymbol{X}}_n)},
\label{EQ:Postconv_Suff_P1}
\end{equation}  
since $g^n(\underline{\boldsymbol{X}}_n)$ is non-zero for each $n$ with $G^\infty$ probability one. 
Let us first consider the numerator in (\ref{EQ:Postconv_Suff_P1}) 
and define $E_n$ to be the event that the numerator is greater than $(b_n+c_n)/\delta$.
Note that, $G^n(E_n) \leq G^n(E_n\cap S_n^c) + G^n(S_n)$ for any sequence of measurable sets $S_n \in \mathcal{B}_n$.
So, taking $S_n$ to be the critical sets of Assumption (A3), we get 
\begin{eqnarray}
&& G^n(E_n\cap S_n^c) = \int_{E\cap S_n^c} G^n(d\underline{\boldsymbol{x}}_n)\nonumber\\
&\leq& \frac{\delta}{(b_n + c_n)} \int_{S_n^c}\frac{m_n^{(\alpha)}(\underline{\boldsymbol{x}}_n, A_n^c)}{
	M_n^{(\alpha)}(\boldsymbol{\chi}_n, \Theta_n)g^n(\underline{\boldsymbol{x}}_n)}G^n(d\underline{\boldsymbol{x}}_n)
~~~~\mbox{[by Markov's inequality and definition of $E_n$]}\nonumber\\
&=& \frac{\delta}{(b_n + c_n)M_n^{(\alpha)}(\boldsymbol{\chi}_n, \Theta_n)} 
\int_{S_n^c} m_n^{(\alpha)}(\underline{\boldsymbol{x}}_n, A_n^c)d\underline{\boldsymbol{x}}_n \nonumber\\
&=& \frac{\delta}{(b_n + c_n)M_n^{(\alpha)}(\boldsymbol{\chi}_n, \Theta_n)} 
\int_{S_n^c} \int_{A_n^c}\exp(q_n^{(\alpha)}(\underline{\boldsymbol{x}}_n|\boldsymbol{\theta}))\pi_n(\boldsymbol{\theta})d\boldsymbol{\theta} 
d\underline{\boldsymbol{x}}_n\nonumber\\
&=&\frac{\delta}{(b_n + c_n)M_n^{(\alpha)}(\boldsymbol{\chi}_n, \Theta_n)} 
\int_{A_n^c} Q_n^{(\alpha)}(S_n^c|\boldsymbol{\theta})\pi_n(\boldsymbol{\theta})d\boldsymbol{\theta} 
~~~~~~~~\mbox{[by Fubini Theorem]}\nonumber\\
&\leq&\frac{\delta}{(b_n + c_n)M_n^{(\alpha)}(\boldsymbol{\chi}_n, \Theta_n)} \left[
\int_{B_n} Q_n^{(\alpha)}(S_n^c|\boldsymbol{\theta})\pi_n(\boldsymbol{\theta})d\boldsymbol{\theta} 
+ \int_{C_n} Q_n^{(\alpha)}(S_n^c|\boldsymbol{\theta})\pi_n(\boldsymbol{\theta})d\boldsymbol{\theta} 
\right]\nonumber\\
&& ~~~~~~~\mbox{[for the sets $B_n$ and $C_n$ from Assumptions (A1)--(A3)]}\nonumber\\
&\leq& \frac{\delta}{(b_n + c_n)M_n^{(\alpha)}(\boldsymbol{\chi}_n, \Theta_n)} \left[
\int_{B_n} Q_n^{(\alpha)}(\boldsymbol{\chi}_n|\boldsymbol{\theta})\pi_n(\boldsymbol{\theta})d\boldsymbol{\theta} 
+ \int_{C_n} \frac{Q_n^{(\alpha)}(S_n^c|\boldsymbol{\theta})}{Q_n^{(\alpha)}(\boldsymbol{\chi}_n|\boldsymbol{\theta})}
Q_n^{(\alpha)}(\boldsymbol{\chi}_n|\boldsymbol{\theta})\pi_n(\boldsymbol{\theta})d\boldsymbol{\theta} 
\right]\nonumber\\
&\leq& \frac{\delta}{(b_n + c_n)M_n^{(\alpha)}(\boldsymbol{\chi}_n, \Theta_n)} \left[
M_n^{(\alpha)}(\boldsymbol{\chi}_n, B_n) + \sup_{\boldsymbol{\theta}\in C_n} \frac{Q_n^{(\alpha)}(S_n^c|\boldsymbol{\theta})}{Q_n^{(\alpha)}(\boldsymbol{\chi}_n|\boldsymbol{\theta})}
M_n^{(\alpha)}(\boldsymbol{\chi}_n, C_n)\right]\nonumber\\
&\leq& \frac{\delta}{(b_n + c_n)} \left[b_n + c_n 
\frac{M_n^{(\alpha)}(\boldsymbol{\chi}_n, C_n)}{M_n^{(\alpha)}(\boldsymbol{\chi}_n, \Theta_n)}\right]
~~~~~~~~\mbox{[by Assumptions (A2) and (A3)]}\nonumber\\
&\leq& \delta.\nonumber
\end{eqnarray} 
Hence, $G^n(E_n) \leq \delta + G^n(S_n)$ and using Assumption (A3) we get 
$\displaystyle\limsup_{n\rightarrow\infty} G^n(E_n) \leq \delta$.
Further, by Assumption (A4) the denominator in (\ref{EQ:Postconv_Suff_P1}) is less than $a_n$  
has probability tending to zero. Combining the numerator and denominator probabilities 
(using the bound by the union of events related to numerator and denominator), 
we get the desired result (\ref{EQ:Postconv_Suff1}).

To prove the second part (\ref{EQ:Postconv_Suff2}), we proceed as before by noting that 
$P(\underline{\boldsymbol{X}}_n\in E_n~\mbox{ i.o.}) \leq P(\underline{\boldsymbol{X}}_n\in E_n\cap S_n^c ~\mbox{ i.o.})
+ P(\underline{\boldsymbol{X}}_n\in S_n~\mbox{ i.o.})$. Then, defining $E_n$ with any summable sequence $\delta_n$
and proceeding as before, we get $P(\underline{\boldsymbol{X}}_n\in E_n\cap S_n^c ~\mbox{ i.o.})=0$ by Borel-Cantelli Lemma. 
Next, by Assumption (A3)$^\ast$, we have $P(\underline{\boldsymbol{X}}_n\in S_n~\mbox{ i.o.})=0$ and hence
$P(\underline{\boldsymbol{X}}_n\in E_n~\mbox{ i.o.})=0$. Then, the desired result (\ref{EQ:Postconv_Suff2}) follows 
by noting that the  denominator in (\ref{EQ:Postconv_Suff_P1}) is less than $a_n$ infinitely often 
with probability zero by Assumption (A4)$^\ast$.
\hfill{$\square$}

\begin{lemma}
	Suppose, for some sequence of constants $r_n$, we have	
	\begin{eqnarray}
	\displaystyle\lim_{n\rightarrow\infty} P\left(\pi_n^{(\alpha)}\left(A_n^c|\underline{\boldsymbol{X}}_n\right)\leq r_n \right)=1.
	\label{EQ:Postconv_Nec1}
	\end{eqnarray}
	Then, for any sequences $b_n$ and $c_n$ satisfying $b_nc_n\geq r_n$, 
	there exists parameter sets $B_n, C_n \subset \Theta_n$ such that Conditions (A1)--(A3) hold. 
	
	Moreover, if additionally we have 
	\begin{eqnarray}
	P\left(\pi_n^{(\alpha)}\left(A_n^c|\underline{\boldsymbol{X}}_n\right)> r_n \mbox{ i.o.}\right)=0,
	\label{EQ:Postconv_Nec2}
	\end{eqnarray}
	then Conditions (A1), (A2) and (A3)$^\ast$ hold.
	\label{LEM:PostConv_Nec}
\end{lemma}
\noindent\textbf{Proof:}
Let us define $S_n =\left\{\underline{\boldsymbol{x}}_n : \pi_n^{(\alpha)}\left(A_n^c|\underline{\boldsymbol{x}}_n\right) > r_n \right\}$
so that $\displaystyle\lim_{n\rightarrow\infty} G^n(S_n) =0$ by Assumption (\ref{EQ:Postconv_Nec1}). 
Next, for any sequence $c_n$, we construct the parameter sets 
$$
C_n =\left\{\boldsymbol{\theta} : \frac{Q_n^{(\alpha)}(S_n^c|\boldsymbol{\theta})}{Q_n^{(\alpha)}(\boldsymbol{\chi}_n|\boldsymbol{\theta})}
\leq c_n \right\}, ~~ 
B_n =\left\{\boldsymbol{\theta}\in A_n^c : 
\frac{Q_n^{(\alpha)}(S_n^c|\boldsymbol{\theta})}{Q_n^{(\alpha)}(\boldsymbol{\chi}_n|\boldsymbol{\theta})}> c_n \right\}.
$$
Then, Conditions (A1) and (A3) hold by constructions of $C_n$ and $B_n$. 
Finally, to show Condition (A2), note that 
$m_n^{(\alpha)}(\underline{\boldsymbol{x}}_n, A_n^c) 
\leq r_n M_n^{(\alpha)}(\boldsymbol{\chi}_n, \Theta_n)m_n^{(\alpha)}(\underline{\boldsymbol{x}}_n)$
for all $\underline{\boldsymbol{x}}_n\in S_n^c$ by its definition. Then, 
\begin{eqnarray}
\frac{M_n^{(\alpha)}(\boldsymbol{\chi}_n, B_n)}{M_n^{(\alpha)}(\boldsymbol{\chi}_n, \Theta_n)}
&=& \frac{1}{M_n^{(\alpha)}(\boldsymbol{\chi}_n, \Theta_n)}
\int_{B_n} Q_n^{(\alpha)}(\boldsymbol{\chi}_n|\boldsymbol{\theta})\pi_n(\boldsymbol{\theta}) d\boldsymbol{\theta}
\nonumber\\
&\leq& \frac{1}{M_n^{(\alpha)}(\boldsymbol{\chi}_n, \Theta_n)c_n}
\int_{A_n^c} Q_n^{(\alpha)}(S_n^c|\boldsymbol{\theta})\pi_n(\boldsymbol{\theta}) d\boldsymbol{\theta}\nonumber\\
&&~~~~~~~~~~~~~~~~\mbox{[by Definition of $B_n$ and Markov's inequality]}\nonumber\\
&\leq& \frac{1}{M_n^{(\alpha)}(\boldsymbol{\chi}_n, \Theta_n)c_n}
\int_{A_n^c} \int_{S_n} \exp(q_n^{(\alpha)}(\underline{\boldsymbol{x}}_n|\boldsymbol{\theta})) d\underline{\boldsymbol{x}}_n
\pi_n(\boldsymbol{\theta}) d\boldsymbol{\theta}
\nonumber\\
&\leq& \frac{1}{M_n^{(\alpha)}(\boldsymbol{\chi}_n, \Theta_n)c_n}
\int_{S_n} m_n^{(\alpha)}(\underline{\boldsymbol{x}}_n, A_n^c)d\underline{\boldsymbol{x}}_n
~~~~~~~~\mbox{[by Fubini Theorem]}\nonumber\\
&\leq& \frac{1}{M_n^{(\alpha)}(\boldsymbol{\chi}_n, \Theta_n)c_n}
\int_{S_n} r_n M_n^{(\alpha)}(\boldsymbol{\chi}_n, \Theta_n)m_n^{(\alpha)}(\underline{\boldsymbol{x}}_n)d\underline{\boldsymbol{x}}_n
~~\mbox{[by the construction of $S_n$]}\nonumber\\
&\leq& \frac{r_n}{c_n} \int_{S_n} m_n^{(\alpha)}(\underline{\boldsymbol{x}}_n)d\underline{\boldsymbol{x}}_n
\nonumber\\
&\leq& \frac{r_n}{c_n} 
~~~~~~~~~~~~~~~~~~~~~~~\mbox{[$\int_{S_n} m_n^{(\alpha)}(\underline{\boldsymbol{x}}_n)d\underline{\boldsymbol{x}}_n
	\leq \int_{\boldsymbol{\chi}_n} m_n^{(\alpha)}(\underline{\boldsymbol{x}}_n)d\underline{\boldsymbol{x}}_n = 1$]} \nonumber\\
&\leq& b_n ~~~~~~~~~~~~~~~~~~~~~~~\mbox{[for any sequence $b_n$ satisfying $b_n c_n \geq r_n$]} \nonumber
\end{eqnarray}

For the second part of the Lemma , we use the same definitions of sets as above. 
Then, by Assumption (\ref{EQ:Postconv_Nec2}), we have $P(\underline{\boldsymbol{X}}_n \in S_n~{ i.o.}) = 0$
and hence Condition (A3)$^\ast$ holds by the construction of $C_n$. Other two conditions then hold similarly as before. 
\hfill{$\square$}

\bigskip
\noindent\textbf{Proof of Theorem 	3.3:}\\
Theorem 3.3 now follows directly from the above two lemmas.

The sufficiency part of the theorem follows from Lemma \ref{LEM:PostConv_Suff}
by taking $b_n=e^{-nr_1}$, $c_n = e^{-nr_2}$, $a_n = e^{-n\epsilon}$ and $\delta_n=e^{-n\Delta}$ (for Part 2) 
with $\epsilon, \Delta>0$ and $\epsilon+\Delta < \min\{r_1, r_2\}$ 
Then, $r_n$ and $r_n' = r_n/\delta_n$ tend to zero exponentially fast.

The Necessity part of the theorem follows from Lemma \ref{LEM:PostConv_Nec} with $r_n = e^{-nr}$
and then letting $b_n = e^{-nr_1}$, $c_n = e^{-nr_2}$ for any $r_1, r_2 >0$ with $r_1 + r_2 \leq r$.
\hfill{$\square$}

\section{Proofs of the Results of Section 4 in Main paper}

\subsection{Proof of Theorem 4.4}
Note that, by the definition of $\widehat{\boldsymbol{\theta}}_\alpha$, it is sufficient to show that 
\begin{eqnarray}
\sup_{\boldsymbol{\theta}\in A_n^c}\widetilde{\pi}_n^{(\alpha)}({\boldsymbol{\theta}})
\widetilde{q}_n^{(\alpha)}(\underline{\boldsymbol{X}}_n|{\boldsymbol{\theta}})
< \sup_{\boldsymbol{\theta}}\widetilde{\pi}_n^{(\alpha)}(\boldsymbol{\theta}) \widetilde{q}_n^{(\alpha)}(\underline{\boldsymbol{X}}_n|\boldsymbol{\theta})
e^{-n\delta_n}~ ~a.s.[G],~~~\mbox{ for all large }n.
\label{EQ:P1}
\end{eqnarray} 
Now, by the information denseness assumption, Theorem 3.2 of the main paper
implies that $G^n$ and $M_n^{(\alpha)}$ merge in probability.
Therefore, the exponential convergence of $\pi_n^{(\alpha)}\left(A_n^c|\underline{\boldsymbol{X}}_n\right)$
is equivalent to 
$$
\sum_{\boldsymbol{\theta}\in A_n^c} \widetilde{\pi}_n^{(\alpha)}(\boldsymbol{\theta}) 
\widetilde{q}_n^{(\alpha)}(\underline{\boldsymbol{X}}_n|\boldsymbol{\theta}) 
\leq m_n^{(\alpha)}(\underline{\boldsymbol{X}}_n)e^{-nr_1}
< g_n(\underline{\boldsymbol{X}}_n)e^{-nr}~ ~a.s.[G],~~~\mbox{ for all large }n, ~\mbox{for some }r_1, r>0.
$$
Let us now choose a $\boldsymbol{\theta}^\ast \in \Theta$ such that 
$KLD(g, \widetilde{q}^{(\alpha)}(\cdot|\boldsymbol{\theta}^\ast)) < r/4$.
Then, using SLLN along with Assumption (17) of the main paper, we get 
$$
g_n(\underline{\boldsymbol{X}}_n) < \widetilde{\pi}_n^{(\alpha)}(\boldsymbol{\theta}^\ast) \widetilde{q}^{(\alpha)}(\underline{\boldsymbol{X}}_n|\boldsymbol{\theta}^\ast)e^{nr/2}
~ ~a.s.[G],~~~\mbox{ for all large }n.
$$
Therefore, for all large $n$, we have with $a.s.[G]$,
\begin{eqnarray}
\sup_{\boldsymbol{\theta}\in A_n^c}\widetilde{\pi}_n^{(\alpha)}({\boldsymbol{\theta}})
\widetilde{q}_n^{(\alpha)}(\underline{\boldsymbol{X}}_n|{\boldsymbol{\theta}})
&\leq& \sum_{\boldsymbol{\theta}\in A_n^c}\widetilde{\pi}_n^{(\alpha)}({\boldsymbol{\theta}})
\widetilde{q}_n^{(\alpha)}(\underline{\boldsymbol{X}}_n|{\boldsymbol{\theta}})
\nonumber\\
&<&  g_n(\underline{\boldsymbol{X}}_n)e^{-nr}
\nonumber\\
&<& \widetilde{\pi}_n^{(\alpha)}(\boldsymbol{\theta}^\ast) \widetilde{q}^{(\alpha)}(\underline{\boldsymbol{X}}_n|\boldsymbol{\theta}^\ast)e^{-nr/2}
\nonumber\\
&<&\sup_{\boldsymbol{\theta}}\widetilde{\pi}_n^{(\alpha)}(\boldsymbol{\theta}) \widetilde{q}_n^{(\alpha)}(\underline{\boldsymbol{X}}_n|\boldsymbol{\theta})
e^{-n\delta_n}.\nonumber
\end{eqnarray}

This completes the proof that $ \widehat{\boldsymbol{\theta}}_\alpha \in A_n$ $a.s.[G]$, for all sufficiently large $n$.
\hfill{$\square$}

\subsection{Proof of Theorem 4.5}
Using the equivalence of $d_1$ and $d_H$ (the Hellinger metric), it is enough to show that
$\pi_n^{(\alpha)}\left( A^c|\underline{\boldsymbol{X}}_n\right)$ is exponentially small with probability one,
with $A=\left\{\boldsymbol{\theta} : d_H(g, f_{\boldsymbol{\theta}}) \geq \epsilon \right\}$
for each fixed $\epsilon>0$. Note that,  $G^n$ and $M_n^{(\alpha)}$ merge in probability
by applying Theorem 3.2 of the main paper.
So, we will use Theorem 3.3 by constructing suitable 
parameter sets $B_n$ and $C_n$ with $A\cup B_n \cup C_n = \Theta$. 

Put $B_n = \left\{ \boldsymbol{\theta} : \pi_n (\boldsymbol{\theta}) <e^{-n\epsilon/4} \right\}$
and  $C_n = \left\{ \boldsymbol{\theta}\in A^c : \pi_n (\boldsymbol{\theta}) \geq e^{-n\epsilon/4} \right\}$.
Then, clearly $A\cup B_n \cup C_n = \Theta$. Further, for some $\tau\in (0,1)$,
\begin{eqnarray}
\frac{M_n^{(\alpha)}(\boldsymbol{\chi}_n, B_n)}{M_n^{(\alpha)}(\boldsymbol{\chi}_n, \Theta_n)}
= \sum_{\boldsymbol{\theta}\in B_n} 
\frac{Q_n^{(\alpha)}(\chi_n|\boldsymbol{\theta})}{M_n^{(\alpha)}(\boldsymbol{\chi}_n, \Theta_n)} 
\pi_n(\boldsymbol{\theta})
\leq \frac{e^{-\frac{n(1-\tau)\epsilon}{4}}}{M_n^{(\alpha)}(\boldsymbol{\chi}_n, \Theta_n)} 
\sum_{\boldsymbol{\theta}\in B_n} Q_n^{(\alpha)}(\chi_n|\boldsymbol{\theta})\pi_n(\boldsymbol{\theta})^\tau
\nonumber
\end{eqnarray}
But, since the prior sequence $\pi_n$ satisfies Assumption (17) of the main paper, we get,
for all sufficiently large $n$, (assuming all the relevant quantities exists finitely)
\begin{eqnarray}
M_n^{(\alpha)}(\boldsymbol{\chi}_n, \Theta_n)
&=& \sum_{\boldsymbol{\theta}\in \Theta} Q_n^{(\alpha)}(\chi_n|\boldsymbol{\theta})\pi_n(\boldsymbol{\theta})
\nonumber\\
&\geq & e^{-n(1-\tau)\epsilon/8}\sum_{\boldsymbol{\theta}\in \Theta} 
Q_n^{(\alpha)}(\chi_n|\boldsymbol{\theta})\pi_n(\boldsymbol{\theta})^\tau
\nonumber\\
&\geq & e^{-n(1-\tau)\epsilon/8}\sum_{\boldsymbol{\theta}\in B_n} 
Q_n^{(\alpha)}(\chi_n|\boldsymbol{\theta})\pi_n(\boldsymbol{\theta})^\tau,
\nonumber
\end{eqnarray}
and hence 
$$
\frac{M_n^{(\alpha)}(\boldsymbol{\chi}_n, B_n)}{M_n^{(\alpha)}(\boldsymbol{\chi}_n, \Theta_n)}
\leq e^{-n(1-\tau)\epsilon/8}.
$$
Thus, the first two conditions of Theorem 3.3 hold.
For the third condition related to $C_n$, note that $\sum_{\boldsymbol{\theta}\in C_n}\pi_n(\boldsymbol{\theta}) \leq 1$
and so the number of points in $C_n$ is less than $e^{n\epsilon/4}$. 
Then, consider the likelihood ratio test for $g_n$ against $\left\{\frac{\exp(q_n^{(\alpha)}(\cdot|\boldsymbol{\theta}))}{Q_n^{(\alpha)}(\boldsymbol{\chi}_n| \boldsymbol{\theta})} 
: \boldsymbol{\theta} \in C_n \right\}$ having the critical sets
$$
S_n = \left\{\underline{\boldsymbol{x}}_n : \max_{\boldsymbol{\theta}\in C_n} 
\frac{\exp(q_n^{(\alpha)}(\underline{\boldsymbol{x}}_n|\boldsymbol{\theta}))}{
	Q_n^{(\alpha)}(\boldsymbol{\chi}_n| \boldsymbol{\theta})} > g_n(\underline{\boldsymbol{x}}_n) \right\}.
$$ 
We will show that this $S_n$ serves as the desired set in the required condition (A3) on $C_n$.
For note that, $S_n=\cup_{\boldsymbol{\theta}\in C_n}S_{n,\boldsymbol{\theta}}$, where
$
S_{n,\boldsymbol{\theta}} 
= \left\{\underline{\boldsymbol{x}}_n : \left[
\frac{\exp(q_n^{(\alpha)}(\underline{\boldsymbol{x}}_n|\boldsymbol{\theta}))}{
	Q_n^{(\alpha)}(\boldsymbol{\chi}_n| \boldsymbol{\theta})}\right]^{1/2} > g_n(\underline{\boldsymbol{x}}_n)^{1/2} \right\}.
$
But for each of these sets, we get from Markov inequality that, 
$$
G_n(S_{n,\boldsymbol{\theta}} ) \leq \left[1 - \frac{1}{2}d_H(g, f_{\boldsymbol{\theta}}) \right]^n < e^{-n\epsilon/2},
$$
and hence $G_n(S_n)< e^{-n\epsilon/4}$. Similarly, we can also show that 
$$
\frac{Q_n^{(\alpha)}(S_n^c| \boldsymbol{\theta})}{Q_n^{(\alpha)}(\boldsymbol{\chi}_n| \boldsymbol{\theta})}
\leq \frac{Q_n^{(\alpha)}(S_{n,\boldsymbol{\theta}}^c|\boldsymbol{\theta})}{
	Q_n^{(\alpha)}(\boldsymbol{\chi}_n| \boldsymbol{\theta})} <e^{-n\epsilon/2},
$$
uniformly over $\boldsymbol{\theta}\in C_n$.
Hence, all the required conditions of Theorem 3.3 hold
and we get the first part of the present theorem.

The second part then follows directly from Theorem 4.4  of the main paper.
\hfill{$\square$}

\subsection{Proof of Proposition 4.6}
Fix $\epsilon>0$ and $\underline{\boldsymbol{x}}_n$. 
Let $a = \pi_n^{(\alpha)}\left(A_{n,\epsilon}|\underline{\boldsymbol{x}}_n\right)$. 
If $a=0$, the result is trivial. So, assume $a>0$ and consider the distribution 
$\pi_n^{(\alpha)}\left(d\boldsymbol{\theta}|\underline{\boldsymbol{x}}_n, A_{n,\epsilon}\right)$ 
obtained from $\pi_n^{(\alpha)}\left(A_{n,\epsilon}|\underline{\boldsymbol{x}}_n\right)$ by conditioning on 
$\boldsymbol{\theta}\in A_{n,\epsilon}$. Then, we have
\begin{eqnarray}
~~~~~~~~~~~~~
L_n\left(P, \widehat{P}_\alpha^E\right) &\leq& L_n\left(P, \int_{A_{n,\epsilon}}F_{\boldsymbol{\theta}} \pi_n^{(\alpha)}\left(d\boldsymbol{\theta}|\underline{\boldsymbol{x}}_n\right)\right) 
~~~~~~~~\mbox{[By Monotonicity]}\nonumber\\
&=& L_n\left(P, \int_{A_{n,\epsilon}}(a F_{\boldsymbol{\theta}})
\pi_n^{(\alpha)}\left(d\boldsymbol{\theta}|\underline{\boldsymbol{x}}_n, A_{n,\epsilon}\right)\right) \nonumber\\
&\leq& \int_{A_{n,\epsilon}} L_n(P, a F_{\boldsymbol{\theta}}) \pi_n^{(\alpha)}\left(d\boldsymbol{\theta}|\underline{\boldsymbol{x}}_n\right)
~~~~~~~~\mbox{[By Convexity]} \nonumber\\
&\leq& \int_{A_{n,\epsilon}} \left[L_n(P, F_{\boldsymbol{\theta}}) + \rho(a)\right] 
\pi_n^{(\alpha)}\left(d\boldsymbol{\theta}|\underline{\boldsymbol{x}}_n\right)
~~~~~~~~\mbox{[By Scaling]} \nonumber\\
&\leq& \epsilon + \rho(a).\nonumber
~~~~~~~~~~~~~~~~~~~~~~~~~~~~~~~~~~~~~~~~~~~~~~~~~~~~~~~~~~~~~~~~~~\hfill{\square}
\end{eqnarray}

\section{Proofs of the Results of Section 5 in Main paper}

\subsection{Proof of Theorem 5.3}

By straightforward calculation, it turns out that 
$$
d_1(g_i,f_{i,\boldsymbol{\beta}})=4\Phi\Big(\frac{|\boldsymbol{t}_i^T(\boldsymbol{\beta}-\boldsymbol{\beta}_0)|}{2}\Big)-2,
$$ 
where $\Phi$ is the cumulative distribution function of the standard Normal distribution. 
Therefore, 
\begin{eqnarray}
&&\pi_n^{(\alpha)}\left(\left\{\boldsymbol{\beta}: \frac{1}{n}\sum_{i=1}^{n}d_1(g_i,f_{i,\boldsymbol{\beta}}) \geq \epsilon\right\}
\Bigg|\underline{\boldsymbol{x}}_n\right)\nonumber\\
&&~~~~~~~=\pi_n^{(\alpha)}\left(\left\{\boldsymbol{\beta}: \frac{1}{n}\sum_{i=1}^{n}\Phi\Big(\frac{|\boldsymbol{t}_i^T(\boldsymbol{\beta}-\boldsymbol{\beta}_0)|}{2}\Big) \geq \frac{\epsilon}{4}+\frac{1}{2}\right\}
\Bigg|\underline{\boldsymbol{x}}_n\right)
\nonumber\\
&&~~~~~~~\leq \sum_{i=1}^{n} \pi_n^{(\alpha)}\Big(\{\boldsymbol{\beta}: \Phi\Big(\frac{|\boldsymbol{t}_i^T(\boldsymbol{\beta}-\boldsymbol{\beta}_0)|}{2} \Big)\geq \frac{\epsilon}{4}+\frac{1}{2}\}|\underline{\boldsymbol{x}}_n\Big)
\nonumber\\
&&~~~~~~~=\sum_{i=1}^{n} \pi_n^{(\alpha)}\Big(\{\boldsymbol{\beta}: \frac{|\boldsymbol{t}_i^T(\boldsymbol{\beta}-\boldsymbol{\beta}_0)|}{2} \geq \Phi^{-1}(\frac{1}{2}+\frac{\epsilon}{4}) \}|\underline{\boldsymbol{x}}_n\Big).\nonumber
\end{eqnarray}
For notational simplicity, let us denote the set $A_{i,n}^{c} =
\left\{\boldsymbol{\beta}: \frac{|\boldsymbol{t}_i^T(\boldsymbol{\beta}-\boldsymbol{\beta}_0)|}{2} \geq \Phi^{-1}(\frac{1}{2}+\frac{\epsilon}{4})\right\}$, 
for each $i=1,\ldots,n$ and let $\widehat{\boldsymbol{\beta}}_n$ be the minimum DPD estimator (MDPDE) of $\boldsymbol{\beta}$ under the same model. 
After some basic algebra, using the consistency of the MDPDE under (R1)--(R2) \citep{Ghosh/Basu:2013}, it can be shown that the event $A_{i,n}$ implies 
$$
-2\Phi^{-1}\left(\frac{1}{2}+\frac{\epsilon}{4}\right)-\boldsymbol{t}_i^T(\widehat{\boldsymbol{\beta}}_n-\boldsymbol{\beta}_0) 
\leq \boldsymbol{t}_i^T(\boldsymbol{\beta}-\widehat{\boldsymbol{\beta}}_n) 
\leq 2\Phi^{-1}\left(\frac{1}{2}+\frac{\epsilon}{4}\right)-\boldsymbol{t}_i^T(\widehat{\boldsymbol{\beta}}_n-\boldsymbol{\beta}_0).
$$ 
Now, we use Theorem 2.1 of \cite{Majumder/etc:2019} to approximate the above probability as follows.
\begin{eqnarray}
\pi_n^{(\alpha)}(A_{i,n}^{c}|\underline{\boldsymbol{x}}_n)&=&1-\pi_n^{(\alpha)}(A_{i,n}|\underline{\boldsymbol{x}}_n)
\nonumber\\
&=& 1-\Bigg[\Phi\left(\frac{2\sqrt{n}\Phi^{-1}(\frac{1}{2}+\frac{\epsilon}{4})
	-\sqrt{n}\boldsymbol{t}_i^T(\widehat{\boldsymbol{\beta}}_n-\boldsymbol{\beta}_0)}{\boldsymbol{t}_i^T\boldsymbol{\Psi}_n^{-1}\boldsymbol{t}_i}\right)
\nonumber\\
&& ~~~~~~~~~-\Phi\left(\frac{-2\sqrt{n}\Phi^{-1}(\frac{1}{2}+\frac{\epsilon}{4})
	-\sqrt{n}\boldsymbol{t}_i^T(\widehat{\boldsymbol{\beta}}_n-\boldsymbol{\beta}_0)}{\boldsymbol{t}_i^T\boldsymbol{\Psi}_n^{-1}\boldsymbol{t}_i}\right)
\Bigg]+o_p(1)
\nonumber\\
&=&\Phi(u_n)+\Phi(v_n)+o_p(1), \nonumber
\end{eqnarray} 
where $\boldsymbol{\Psi}_n^{-1} = n[\boldsymbol{D}^T\boldsymbol{D}]$, 
$u_n=\frac{-2\sqrt{n}\Phi^{-1}(\frac{1}{2}+\frac{\epsilon}{4})+\sqrt{n}\boldsymbol{t}_i^T(\widehat{\boldsymbol{\beta}}_n-\boldsymbol{\beta}_0)}{
	\boldsymbol{t}_i^T\boldsymbol{\Psi}_n^{-1}\boldsymbol{t}_i}$ 
and $v_n=\frac{-2\sqrt{n}\Phi^{-1}(\frac{1}{2}+\frac{\epsilon}{4})-\sqrt{n}\boldsymbol{t}_i^T(\widehat{\boldsymbol{\beta}}_n-\boldsymbol{\beta}_0)}{
	\boldsymbol{t}_i^T\boldsymbol{\Psi}_n^{-1}\boldsymbol{t}_i}.$

Now, using the consistency of the MDPDE under (R1)--(R2) \citep{Ghosh/Basu:2013}, along with (R1),
we get 
$\boldsymbol{t}_i^T(\widehat{\boldsymbol{\beta}}_n-\boldsymbol{\beta}_0) \xrightarrow{p}0 $, for all $i=1,...,n$. 
So, 
$$
\lim\limits_{n \rightarrow \infty}
P\left(\left|\boldsymbol{t}_i^T(\widehat{\boldsymbol{\beta}}_n-\boldsymbol{\beta}_0)\right|<\Phi^{-1}\left(\frac{1}{2}+\frac{\epsilon}{4}\right)\right)=1.
$$
Further, by (R2), $\displaystyle{\max_{1 \leq i \leq n}} \boldsymbol{t}_i^T\boldsymbol{\Psi}_n^{-1}\boldsymbol{t}_i
=\displaystyle{\max_{1 \leq i \leq n}}n\boldsymbol{t}_i^T[\boldsymbol{D}^T\boldsymbol{D}]^{-1}\boldsymbol{t}_i= O(1)$
and hence we get 
$$P\left(u_n<\frac{-\sqrt{n}\Phi^{-1}(\frac{1}{2}+\frac{\epsilon}{4})}{\boldsymbol{t}_i^T\boldsymbol{\Psi}_n^{-1}\boldsymbol{t}_i}\right) \rightarrow 1$$ 
and 
$$P\left(v_n<\frac{-\sqrt{n}\Phi^{-1}(\frac{1}{2}+\frac{\epsilon}{4})}{\boldsymbol{t}_i^T\boldsymbol{\Psi}_n^{-1}\boldsymbol{t}_i}\right) \rightarrow 1,$$ 
as $n\rightarrow \infty$. 
Since $\Phi$ is continuous, in turn, we have 
$$ 
P\left(\Phi(u_n)+\Phi(v_n)
<2\Phi\left(\frac{-\sqrt{n}\Phi^{-1}(\frac{1}{2}+\frac{\epsilon}{4})}{\boldsymbol{t}_i^T\boldsymbol{\Psi}_n^{-1}\boldsymbol{t}_i}\right)\right) 
\rightarrow 1,
~~~~~\mbox{ as }~ n \rightarrow \infty.
$$ 
Let us now denote $l_n=\frac{\sqrt{n}\Phi^{-1}(\frac{1}{2}+\frac{\epsilon}{4})}{\boldsymbol{t}_i^T\boldsymbol{\Psi}_n^{-1}\boldsymbol{t}_i}$. 
Note that $l_n>0$ for all $n$, and $l_n \rightarrow \infty$ as $n \rightarrow \infty$. 
So, $\Phi(-l_n) \sim \dfrac{\phi(l_n)}{l_n}$, where $\phi$ is the density of the standard normal distribution. 
Hence, for large $n$, $\Phi(-l_n)\leq 2\dfrac{\phi(l_n)}{l_n}$ so that we get, with probability tending to one, 
$$
\Phi(u_n)+\Phi(v_n) < 4\dfrac{\phi(l_n)}{l_n} =\dfrac{4\boldsymbol{t}_i^T\boldsymbol{\Psi}_n^{-1}\boldsymbol{t}_i}{\sqrt{2\pi}\sqrt{n}\Phi^{-1}(\frac{1}{2}+\frac{\epsilon}{4})}
e^{-\dfrac{n(\Phi^{-1}(\frac{1}{2}+\frac{\epsilon}{4}))^2}{2(\boldsymbol{t}_i^T\boldsymbol{\Psi}_n^{-1}\boldsymbol{t}_i)^2}}.
$$
Next, by Assumption (R2), there exists a constant $C_0>0$ satisfying $\boldsymbol{t}_i^T\boldsymbol{\Psi}_n^{-1}\boldsymbol{t}_i \leq C_0$.
Thus, with probability tending to one, we have
$$
\Phi(u_n)+\Phi(v_n) < \dfrac{C_1}{\eta \sqrt{n}}e^{-nC_2}
$$ 
for some constants $C_1,C_2 >0$ and $\eta=\Phi^{-1}(\frac{1}{2}+\frac{\epsilon}{4})$. 
Also, since $ C_2>0$, we get $ \dfrac{C_1\sqrt{n}}{\eta } < e^{\frac{nC_2}{2}} $ for all sufficiently large $n$.
Therefore, with probability tending to one, we have 
$$
\pi_n^{(\alpha)}\Big(\{\boldsymbol{\beta}: \frac{1}{n}\sum_{i=1}^{n}d_1(g_i,f_{i,\boldsymbol{\beta}}) \geq \epsilon\}|\underline{\boldsymbol{x}}_n\Big) 
< \dfrac{C_1\sqrt{n}}{\eta }e^{-nC_2} < e^{\frac{-nC_2}{2}},
$$
%
and hence the theorem holds with $r=\frac{C_2}{2}$. 
\hfill{$\square$}

\subsection{Proof of Theorem 5.4}

In this case, the parameter is $\boldsymbol{\theta}=(\boldsymbol{\beta},\sigma^2)$ and, after some basic algebra, we find that 
$$
d_1(g_i,f_{i,\boldsymbol{\theta}}) \leq 2\dfrac{|\sigma-\sigma_0|}{\sigma} + 4\Phi\left(\frac{|\boldsymbol{t}_i^T(\boldsymbol{\beta}-\boldsymbol{\beta}_0)|}{2\sigma_0}\right)-2. 
$$
Hence, for any $\epsilon\in(0,1)$, it follows that
\begin{eqnarray}
&&
\pi_n^{(\alpha)}\left(\left\{\boldsymbol{\theta}: \frac{1}{n}\sum_{i=1}^{n}d_1(g_i,f_{i,\boldsymbol{\beta}}) \geq \epsilon\right\} 
\Bigg| \underline{\boldsymbol{x}}_n\right) 
\nonumber\\
&\leq& \pi_n^{(\alpha)}\left(\left\{\boldsymbol{\beta}: \frac{4}{n}\sum_{i=1}^{n}\Phi\Big(\frac{|\boldsymbol{t}_i^T(\boldsymbol{\beta}-\boldsymbol{\beta}_0)|}{2\sigma_0}\Big)-2 \geq \frac{\epsilon}{2}\right\}
\Bigg| \underline{\boldsymbol{x}}_n \right)
+ \pi_n^{(\alpha)}\left( \left\{\sigma: 2\dfrac{|\sigma-\sigma_0|}{\sigma} \geq \frac{\epsilon}{2}\right\}| \underline{\boldsymbol{x}}_n\right)
\nonumber\\
& \leq& \sum_{i=1}^{n} \pi_n^{(\alpha)}\left(\left\{\boldsymbol{\beta}: \frac{|\boldsymbol{t}_i^T(\boldsymbol{\beta}-\boldsymbol{\beta}_0)|}{2\sigma_0} 
\geq \Phi^{-1}\left(\frac{1}{2}+\frac{\epsilon}{8}\right) \right\} \Bigg| \underline{\boldsymbol{x}}_n \right)
+ \pi_n^{(\alpha)}\left( \left\{\sigma: |\sigma-\sigma_0|\geq \frac{\epsilon/4}{1+\epsilon/4}\right\} \Big| \underline{\boldsymbol{x}}_n \right).
\nonumber\\
\label{Eq:5.4.1}
\end{eqnarray}

The first term in above (\ref{Eq:5.4.1}) is exponentially small in probability under (R1)--(R2), as proved in Theorem 5.3. 
However, if $\widehat{\sigma}_n^2$ denotes the MDPDE of $\sigma$, using its consistency under (R1)--(R2) \citep{Ghosh/Basu:2013}, we have
$$
\pi_n^{(\alpha)}\Big( \{\sigma: |\sigma-\sigma_0|\geq \eta\}|\underline{\boldsymbol{x}}_n\Big)
\leq \pi_n^{(\alpha)}\Big( \{\sigma: |\sigma-\widehat{\sigma}_n|\geq \eta/2\}|\underline{\boldsymbol{x}}_n\Big)+ o_p(1).
$$

But 
$\pi_n^{(\alpha)}\Big( \{\sigma: |\sigma-\widehat{\sigma}_n|\geq \eta/2\}|\underline{\boldsymbol{x}}_n\Big)
=\pi_n^{(\alpha)}\Big( \{\sigma: |\sigma^2-\widehat{\sigma}_n^2|\geq \eta'\}|\underline{\boldsymbol{x}}_n \Big)$ for some $\eta'>0$. 
And, it follows from Theorem 2.1 of \cite{Majumder/etc:2019} that 
the posterior distribution of $\sqrt{n}(\sigma^2-\widehat{\sigma}_n^2)$ is $\mathcal{N}(0,\zeta_\alpha)$, 
where $\zeta_\alpha$ is some function of $\sigma_0$ and $\alpha$. 
So, combining them, we get the following approximation by a similar technique used to prove Theorem 5.1 above:
$$
\pi_n^{(\alpha)}\Big( \{\sigma: |\sigma-\widehat{\sigma}_n|\geq \eta/2\}|\underline{\boldsymbol{x}}_n,\boldsymbol{D}\Big)= 2\Phi\Big(\frac{-\sqrt{n}\eta'}{\zeta_\alpha}\Big)+o_p(1).
$$

Note that the right-hand side of the above equation is again exponentially small implying the same for its left-hand side. 
Thus the last term in (\ref{Eq:5.4.1}) is also exponentially small, completing the proof of the theorem.
\hfill{$\square$}

\subsection{Proof of Theorem 5.5}

In the present case of logistic set-up, the $L_1$ distance between the true density $g_i$ and model density $f_{i,\boldsymbol{\beta}}$ turns out to be 
$$
d_1(g_i,f_{i,\boldsymbol{\beta}})=2\lvert p_i(\boldsymbol{\beta}_0)-p_i(\boldsymbol{\beta}) \rvert,
$$
where $p_i(\boldsymbol{\beta})=\dfrac{e^{\boldsymbol{t}_i'\boldsymbol{\beta}}}{1+e^{\boldsymbol{t}_i'\boldsymbol{\beta}}}$. Recall that $\boldsymbol{\beta}_0$ is the true parameter. Now, applying the mean value theorem on the function $g(t)=\frac{e^t}{1+e^t}$,  we get
$$
d_1(g_i,f_{i,\boldsymbol{\beta}})=2\lvert \boldsymbol{t}_i'\boldsymbol{\beta}_0-\boldsymbol{t}_i'\boldsymbol{\beta}  \rvert \dfrac{e^{t_i}}{(1+e^{t_i})^2}\leq 2\lvert \boldsymbol{t}_i'(\boldsymbol{\beta}-\boldsymbol{\beta}_0)  \rvert.
$$
Hence, we get
\begin{eqnarray}
\pi_n^{(\alpha)}\left(\left\{\boldsymbol{\beta}: \frac{1}{n}\sum_{i=1}^{n}d_1(g_i,f_{i,\boldsymbol{\beta}}) \geq \epsilon\right\}
\Bigg| \underline{\boldsymbol{x}}_n \right)
&\leq& \pi_n^{(\alpha)}\left(\left\{\boldsymbol{\beta}: \frac{2}{n}\sum_{i=1}^{n}\lvert \boldsymbol{t}_i'(\boldsymbol{\beta}-\boldsymbol{\beta}_0) \rvert 
\geq \epsilon\right\} \Bigg|\underline{\boldsymbol{x}}_n \right)
\nonumber\\
&\leq& \sum_{i=1}^{n} \pi_n^{(\alpha)}\left(\left\{\boldsymbol{\beta}: \lvert \boldsymbol{t}_i'(\boldsymbol{\beta}-\boldsymbol{\beta}_0) \rvert 
\geq \frac{\epsilon}{2}\right\} \Big|\underline{\boldsymbol{x}}_n \right).
\nonumber
\end{eqnarray}
Denote the set $\{\boldsymbol{\beta}: \lvert \boldsymbol{t}_i'(\boldsymbol{\beta}-\boldsymbol{\beta}_0) \rvert \geq \frac{\epsilon}{2}\}=A_{i,n}^c$
and let $\widehat{\boldsymbol{\beta}}_n$ is the MDPDE of $\boldsymbol{\beta}$. 
After some basic algebra, it can be shown that the event $A_{i,n}$ implies
$$
-\dfrac{\epsilon}{2}-\boldsymbol{t}_i'(\widehat{\boldsymbol{\beta}}_n-\boldsymbol{\beta}_0)\leq \boldsymbol{t}_i'(\boldsymbol{\beta}-\widehat{\boldsymbol{\beta}}_n) \leq \dfrac{\epsilon}{2}-\boldsymbol{t}_i'(\widehat{\boldsymbol{\beta}}_n-\boldsymbol{\beta}_0).
$$
Now, by the consistency of the MDPDE under (R1) and (R3) \citep{Majumder/etc:2019}, 
we have $\widehat{\boldsymbol{\beta}}_n-\boldsymbol{\beta}_0 \xrightarrow[]{p}0$, 
and the supremum of the elements of the vectors $\boldsymbol{t}_i$ are bounded by (R3).
Hence,  $\displaystyle{\max_{i}}~ \boldsymbol{t}_i'(\widehat{\boldsymbol{\beta}}_n-\boldsymbol{\beta}_0)\xrightarrow[]{p}0$, implying 
$ -\epsilon \leq \boldsymbol{t}_i'(\boldsymbol{\beta}-\widehat{\boldsymbol{\beta}}_n) \leq \epsilon$ with probability tending to one.
Then, following similar line of calculations as in the proof of Theorem 5.1, we find that 
$$
\pi_n^{(\alpha)} \Big( -\epsilon \leq \boldsymbol{t}_i'(\boldsymbol{\beta}-\widehat{\boldsymbol{\beta}}_n) \leq \epsilon| \underline{\boldsymbol{x}}_n\Big)
=2\Phi\left(\dfrac{\sqrt{n}\epsilon}{\boldsymbol{t}_i'\boldsymbol{\Psi}_n(\boldsymbol{\beta})^{-1}\boldsymbol{t}_i}\right)-1+o_p(1).
$$
So, with probability tending to one, we get 
$$
\pi_n^{(\alpha)} \Big( A_{i,n}^c| \underline{\boldsymbol{x}}_n, \boldsymbol{D} \Big)\leq 2\left[1-\Phi\left(\dfrac{\sqrt{n}\epsilon}{\boldsymbol{t}_i'\boldsymbol{\Psi}_n(\boldsymbol{\beta})^{-1}\boldsymbol{t}_i}\right)\right]
=2\Phi\left(-\dfrac{\sqrt{n}\epsilon}{\boldsymbol{t}_i'\boldsymbol{\Psi}_n(\boldsymbol{\beta})^{-1}\boldsymbol{t}_i}\right).
$$ 
Using boundedness of $\boldsymbol{t}_i'\boldsymbol{\Psi}_n(\boldsymbol{\beta})^{-1}\boldsymbol{t}_i$ and the fact that $\Phi(-l_n) \sim \dfrac{\phi(l_n)}{l_n}$ as $n\rightarrow \infty$, 
we can derive finally that
$$
\pi_n^{(\alpha)} \Big( A_{i,n}^c| \underline{\boldsymbol{x}}_n, \boldsymbol{D} \Big)\leq \dfrac{C_0}{\sqrt{n}}e^{-nC_1}
$$
for some positive constants $C_0$ and $C_1$. Thus 
$$
\pi_n^{(\alpha)}\Big(\{\boldsymbol{\beta}: \frac{1}{n}\sum_{i=1}^{n}d_1(g_i,f_{i,\boldsymbol{\beta}}) \geq \epsilon\}|\underline{\boldsymbol{x}}_n,\boldsymbol{D}\Big)
\leq \sqrt{n}C_0e^{-nC_1}< e^{-nr}
$$
for some $r>0$, with probability tending to one, proving the theorem.

\section{Additional Description of Figures 1, 2 and 3 of the Main Paper}

In Figures 1, 2 and 3 of the main paper it may appear that the MSE is practically zero at $\alpha = 1$ 
(in case of pure data, as well as for each contaminated scenario). 
If that is so, why should we not set $\alpha = 1$ all the time, rather than going through the exercise of choosing an optimal $\alpha$? 
Actually this is a false impression created by the scale of these figures. 
Under contamination, the inflations in the MSE of the estimators corresponding to low values of $\alpha$ are of such a high magnitude, 
that in trying to accommodate them within the same frame of the figure, 
the MSEs of several of the stable estimators corresponding to large values of $\alpha$ 
(and not just for $\alpha = 1$) appear to be zero, or to be very close to it. 
The phenomenon can be better explained by looking at blow-ups of these MSE curves 
by restricting the Y-axis (the MSE axis) to a small range around zero. 
We provide such a representative figure (Figure \ref{FIG:MSE_scaled} below), 
which corresponds to such a blown up MSE curve for the estimation of $\sigma$ with a sample size of $n=100$ 
in the linear regression model with unknown $\sigma$ and the conjugate priors. 
Notice that under pure data, the MSE curve in Figure \ref{FIG:MSE_scaled} is steadily increasing, 
indicating that under the model the performance progressively deteriorates with increasing $\alpha$. 
With increasing contamination the optimal value of $\alpha$ (the value which minimizes the MSE) keeps getting shifted upward; 
the optimal values of $\alpha$ for $\epsilon_C = 0.05, 0.1$ and 0.2 are, approximately, 0.35, 0.5 and 1. 
In this example, therefore, it is clear that the choice of the optimal $\alpha$ is very much a function of the amount of the anomaly in the data, 
and the blanket selection of $\alpha = 1$ does not necessarily provide the best solution in all cases. 
A tuning parameter selection does remain an important component of our methodology as described in Section 7.2 of the main paper. 

Similar phenomena occur for the graphs of MSEs under all the cases reported in Figures 1, 2 and 3 of the main paper;
so we do not repeat them for brevity.

\begin{figure}[!h]
	\centering
	\includegraphics[width=0.7\textwidth]{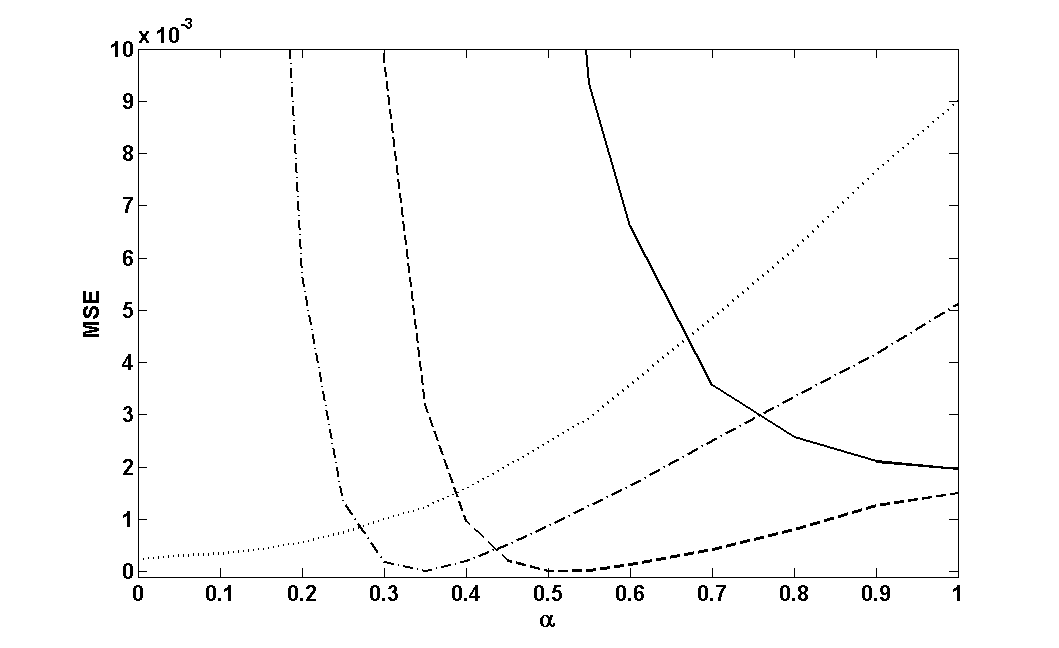}
	\caption{Empirical MSE of the ERPE of $\sigma$ for $n=100$ in the linear regression model with unknown $\sigma$ and the conjugate priors. 
		[Dotted line: $\epsilon_C=0\%$, Dash-Dotted line: $\epsilon_C=5\%$, Dashed line: $\epsilon_C=10\%$, 
		Solid line: $\epsilon_C=20\%$] (rescaled version of Fig. 2(c), lower panel, of the main paper)}
	\label{FIG:MSE_scaled}
\end{figure}

\section{Additional Simulation results for Normal Linear Regression Model with fixed $\sigma$}

\begin{figure}
	\centering
	\subfloat[$n=20$]{
		\includegraphics[width=0.4\textwidth]{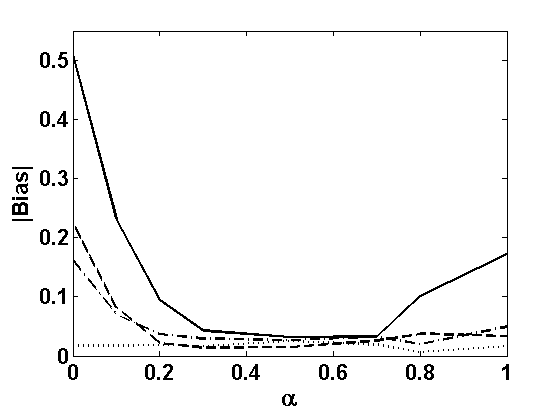}
		\includegraphics[width=0.4\textwidth]{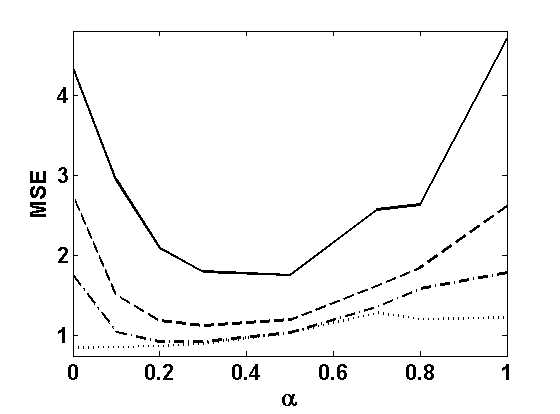}
		\label{FIG:LRM_n20}}
	\\
	\subfloat[$n=50$]{
		\includegraphics[width=0.4\textwidth]{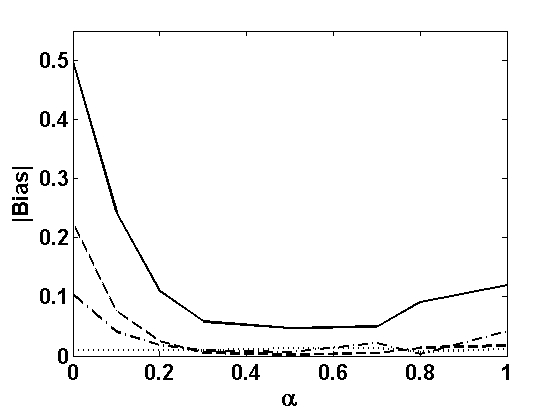}
		\includegraphics[width=0.4\textwidth]{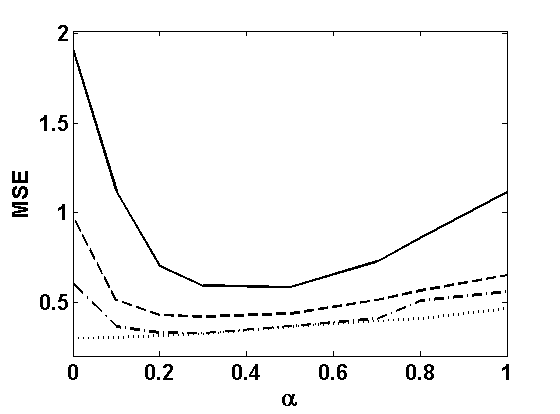}
		\label{FIG:LRM_n50}}
	\\
	\subfloat[$n=100$]{
		\includegraphics[width=0.4\textwidth]{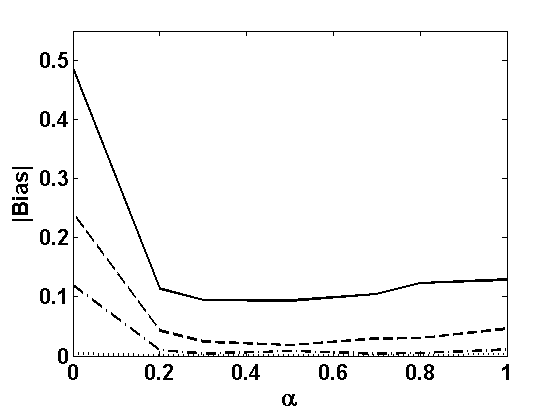}
		\includegraphics[width=0.4\textwidth]{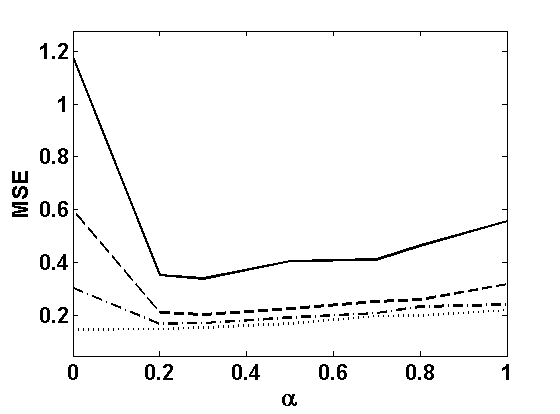}
		\label{FIG:LRM_n100}}
	\caption{Empirical  total absolute bias and total MSE of the ERPE of $\boldsymbol{\beta}$ in the linear regression model 
		with known $\sigma=1$ and the uniform prior. 
		{ [Dotted line: $\epsilon_C=0\%$, Dash-Dotted line: $\epsilon_C=5\%$, Dashed line: $\epsilon_C=10\%$, 
			Solid line: $\epsilon_C=20\%$] }
	}
	\label{FIG:MSE_normalReg}
\end{figure}

\begin{figure}[!h]
	\centering
	\subfloat[$n=20$]{
		\includegraphics[width=0.4\textwidth]{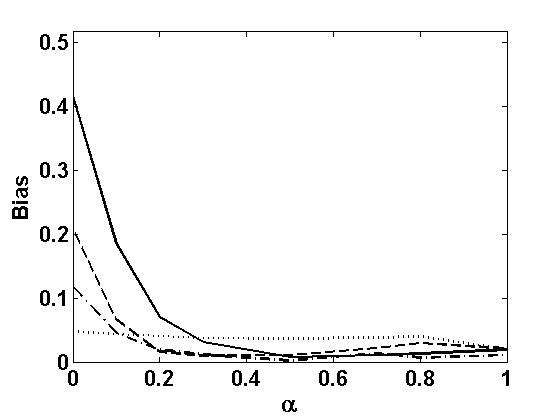}
		\includegraphics[width=0.4\textwidth]{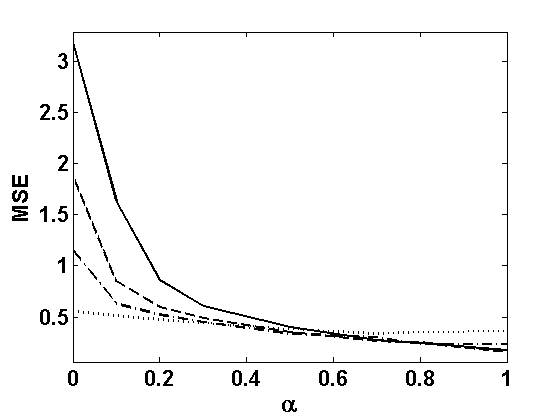}
		\label{FIG:LRM_n20CP}}
	\\
	\subfloat[$n=50$]{
		\includegraphics[width=0.4\textwidth]{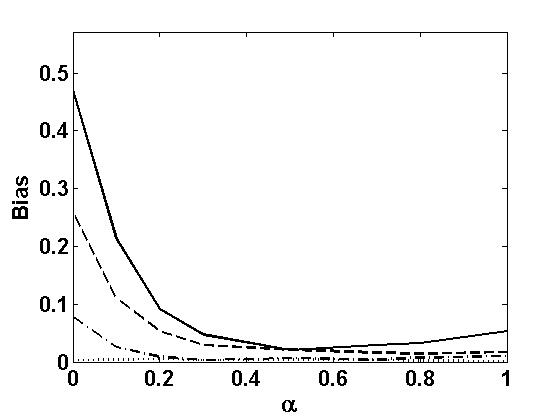}
		\includegraphics[width=0.4\textwidth]{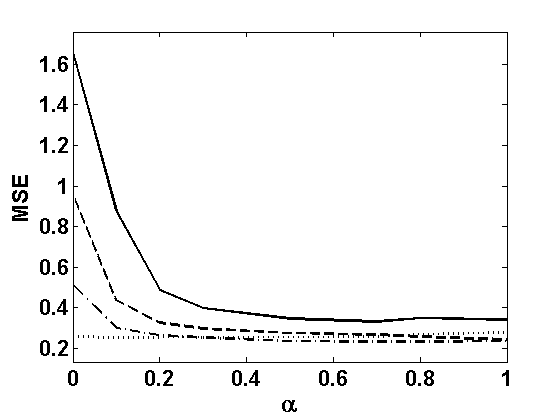}
		\label{FIG:LRM_n50cp}}
	\\
	\subfloat[$n=100$]{
		\includegraphics[width=0.4\textwidth]{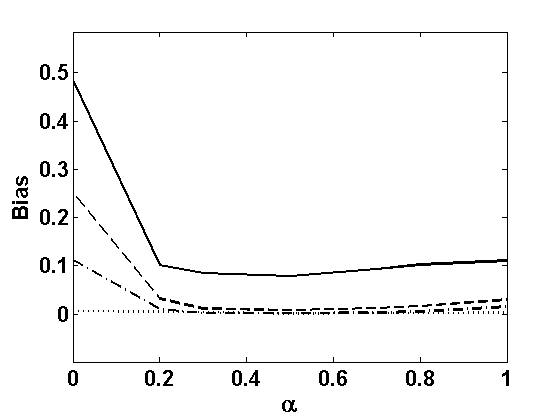}
		\includegraphics[width=0.4\textwidth]{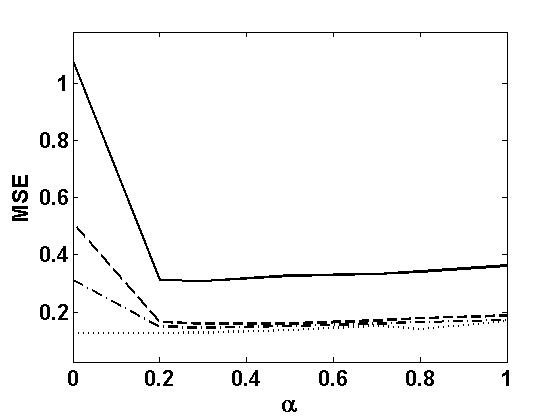}
		\label{FIG:LRM_n100cp}}
	\caption{Empirical total absolute bias and total MSE of the ERPE of $\boldsymbol{\beta}$ in the linear regression model
		with known $\sigma=1$ and the Conjugate normal prior. 
		[Dotted line: $\epsilon_C=0\%$, Dash-Dotted line: $\epsilon_C=5\%$, Dashed line: $\epsilon_C=10\%$, 
		Solid line: $\epsilon_C=20\%$] }
	\label{FIG:MSE_normalRegCP}
\end{figure}

Recall the simulation set-up and notation as in Section 6.1 of the main paper
to illustrate the performance of the ERPE for fixed-design linear regression model with known error variance $\sigma$
and two different choices for the prior.
As the first choice of the prior $\pi(\boldsymbol{\beta})$,
we consider the non-informative uniform prior $\pi(\boldsymbol{\beta}) \equiv 1$.
Secondly, we consider the conjugate normal prior 
$\pi(\boldsymbol{\beta}) \equiv N_k (\boldsymbol{\beta}_0, \tau^2I_k)$
which signifies that the prior belief about our true parameter value 
is quantified by a symmetric structure with uncertainty quantified by $\tau$.
The resulting values of the total absolute bias and the total MSE (over the two components of $\boldsymbol{\beta}$)
are shown in Figure \ref{FIG:MSE_normalReg} and \ref{FIG:MSE_normalRegCP}, respectively.
As noted in the main paper, these results clearly demonstrate the significant improvement 
for the ERPE over the usual Bayes estimators under data contamination
with only a slight loss in case of pure data. 

\clearpage

\clearpage

\newpage
\section{R Codes for Computation of the ERPEs}
%
%
\subsection{Fixed-Design Linear Regression with Unknown Error Variance and Jefferey's Prior}

\begin{lstlisting}[language=R]
library("MASS")
library("parallel")

## Computation of ERPE for Fixed-design Linear Regression Model 
## with unknonw error variance and
## Non informative Jefferey's Prior pi(\beta, \sigma) = 1/sigma^2 


####### Auxiliary Functions  --------------
## Defining the function Q ##
Q=function(beta,sigma,X,y,alpha)   
{
  n=nrow(X)
  p=ncol(X)
  M=X%*%beta
  L=matrix(0,n,1,byrow=TRUE)
  for(i in 1:n)
  {
    L[i]=(1/alpha)*(dnorm(y[i],M[i],sigma))^(alpha)
             -1/(((1+alpha)^(3/2))*((sqrt(2*pi)*sigma)^(alpha)))- (1/alpha)
  }
  L_0=matrix(0,n,1,byrow=TRUE)
  for(i in 1:n)
  {
    L_0[i]=log(dnorm(y[i],M[i],sigma))
  }
  if(alpha>0)
  {
    return(sum(L))
  }
  else
  {
    return(sum(L_0))
  }
}

## Kernel of the prior density ##
prior=function(beta,sigma,true_beta)     
{ 
  return (1/sigma^2)
}

## Kernel of the log of robust posterior density ##
robust_posterior=function(beta,sigma,X,y,alpha,true_beta)   
{
  z1=Q(beta,sigma,X,y,alpha)
  z2=log(prior(beta,sigma,true_beta))
  z=z1+z2
  return (z)
}

## Proposal Density ##
proposalfunction=function(beta,sigma)    
{
  p=length(beta)
  return (c(rnorm(p,mean=beta,sd=rep(1,p)),rexp(1,1/sigma)))
}

## Metropolis Hastings ##
run_metropolis_MCMC <- function(startvalue_beta,startvalue_sigma, 
                                  iterations,X,y,alpha,true_beta)    
{
  p=length(startvalue_beta)
  chain_beta = array(dim = c(iterations+1,p))
  chain_sigma = array(dim = c(iterations+1,1))
  chain_beta[1,] = startvalue_beta
  chain_sigma[1]=startvalue_sigma
  for (i in 1:iterations){
    proposal = proposalfunction(chain_beta[i,],chain_sigma[i])
    proposal_beta=proposal[1:p]
    proposal_sigma=proposal[p+1]
    probab1 
     = exp(robust_posterior(proposal_beta,proposal_sigma,X,y,alpha,true_beta)
        -robust_posterior(chain_beta[i,],chain_sigma[i],X,y,alpha,true_beta))
    logprobab2=log(dexp(chain_sigma[i],(1/proposal_sigma)))
                -log(dexp(proposal_sigma,(1/chain_sigma[i])))
    probab2=exp(logprobab2)
    probab=probab1*probab2
    if (runif(1) < probab){
      chain_beta[i+1,] = proposal_beta
      chain_sigma[i+1] = proposal_sigma
    }else{
      chain_beta[i+1,] = chain_beta[i,]
      chain_sigma[i+1,] = chain_sigma[i,]
    }
  }
  chain=cbind(chain_beta,chain_sigma)
}


####### Computation of the ERPE --------------
data=read.table(file="data.txt")   ## Call the data-file  ##
y=data[,2]      ## Response ##
x=data[,-2]      ## Covariate ##
X=cbind(1,x)    ## Design matrix ##
n=nrow(X)       ## Number of observations ##
p=ncol(X)       ## Dimension of regression parameters excluding sigma ##

## Fit linear regression ## 
res=lm(y~x)     
summary(res)
beta_0=res$coefficients    ## MLE ##

burn=25000      ## burn in ##
max_iter=50000  ## Tortal number of iterations including burn in ##
alpha = 0.5     ## tuning parameter alpha in alpha-likelihood
#set.seed(12345)
output=run_metropolis_MCMC(startvalue_beta=c(-7,3),
                           startvalue_sigma=1,iterations = max_iter,
                           X=X,y=y,alpha=alpha,true_beta=beta_0)
output2=output^2

#Estimated ERPE, the means of the the R^\alpha-posterior distribution
est_mean=colMeans(output[(burn+1):max_iter+1,])
#Estimated variance of the R^\alpha-posterior distribution
est_var=colMeans(output2[(burn+1):max_iter+1,])-(est_mean)^2



\end{lstlisting}

\bigskip\bigskip
\subsection{Fixed-Design Linear Regression with Unknown Error Variance and Conjugate Prior}

\begin{lstlisting}[language=R]
library("MASS")
library("parallel")

## Computation of ERPE for Fixed-design Linear Regression Model 
## with unknonw error variance and
## Conjugate normal-Inverse Gamma prior with prior mean for beta is given ##


####### Auxiliary Functions  --------------
## Defining the function Q ##
Q=function(beta,sigma,X,y,alpha)
{
  n=nrow(X)
  p=ncol(X)
  M=X%*%beta
  L=matrix(0,n,1,byrow=TRUE)
  for(i in 1:n)
  {
    L[i]=(1/alpha)*(dnorm(y[i],M[i],sigma))^(alpha)
             -1/(((1+alpha)^(3/2))*((sqrt(2*pi)*sigma)^(alpha)))- (1/alpha)
  }
  L_0=matrix(0,n,1,byrow=TRUE)
  for(i in 1:n)
  {
    L_0[i]=log(dnorm(y[i],M[i],sigma))
  }
  if(alpha>0)
  {
    return(sum(L))
  }
  else
  {
    return(sum(L_0))
  }
}

## Kernel of the prior density ##
prior=function(beta,sigma,true_beta)      
{
  z1=prod(dnorm(beta,true_beta,sigma))
  z2=sigma^(-5)*exp(-0.5/(sigma^2))
  return (z1*z2)
}

## Kernel of the log of robust posterior density ##
robust_posterior=function(beta,sigma,X,y,alpha,true_beta)  
{
  z1=Q(beta,sigma,X,y,alpha)
  z2=log(prior(beta,sigma,true_beta))
  z=z1+z2
  return (z)
}

## Proposal Density ##
proposalfunction=function(beta,sigma)
{
  p=length(beta)
  return (c(rnorm(p,mean=beta,sd=rep(1,p)),rexp(1,1/sigma)))
}

## Metropolis Hastings ##
run_metropolis_MCMC <- function(startvalue_beta,startvalue_sigma, 
                                iterations,X,y,alpha,true_beta)
{
  p=length(startvalue_beta)
  chain_beta = array(dim = c(iterations+1,p))
  chain_sigma = array(dim = c(iterations+1,1))
  chain_beta[1,] = startvalue_beta
  chain_sigma[1]=startvalue_sigma
  for (i in 1:iterations){
    proposal = proposalfunction(chain_beta[i,],chain_sigma[i])
    proposal_beta=proposal[1:p]
    proposal_sigma=proposal[p+1]
    probab1 
     = exp(robust_posterior(proposal_beta,proposal_sigma,X,y,alpha,true_beta)
       -robust_posterior(chain_beta[i,],chain_sigma[i],X,y,alpha,true_beta))
    logprobab2=log(dexp(chain_sigma[i],(1/proposal_sigma)))
                 -log(dexp(proposal_sigma,(1/chain_sigma[i])))
    probab2=exp(logprobab2)
    probab=probab1*probab2
    if (runif(1) < probab){
      chain_beta[i+1,] = proposal_beta
      chain_sigma[i+1] = proposal_sigma
    }else{
      chain_beta[i+1,] = chain_beta[i,]
      chain_sigma[i+1,] = chain_sigma[i,]
    }
  }
  chain=cbind(chain_beta,chain_sigma)
}


####### Computation of the ERPE --------------
data=read.table(file="data.txt")   ## Call the data-file  ##
y=data[,2]      ## Response ##
x=data[,-2]      ## Covariate ##
X=cbind(1,x)    ## Design matrix ##
n=nrow(X)       ## Number of observations ##
p=ncol(X)       ## Dimension of regression parameters excluding sigma ##

beta_0=c(-8.03,2.95)  ## Given value for the mean of the normal prior ## 
burn=25000
max_iter=50000
alpha = 0.5     ## tuning parameter alpha in alpha-likelihood
#set.seed(12345)
output=run_metropolis_MCMC(startvalue_beta=beta_0,
                           startvalue_sigma=1,iterations = max_iter,
                           X=X,y=y,alpha=alpha,true_beta=beta_0)
output2=output^2
#Estimated ERPE, the means of the the R^\alpha-posterior distribution
est_mean=colMeans(output[(burn+1):max_iter+1,])
#Estimated variance of the R^\alpha-posterior distribution
est_var=colMeans(output2[(burn+1):max_iter+1,])-(est_mean)^2




\end{lstlisting}

\bigskip\bigskip\bigskip\newpage
\subsection{Fixed-Design Logistic Regression with Normal Prior for Regression Coefficient}

\begin{lstlisting}[language=R]
library("MASS")
library("mvtnorm")
library("PACBO")

## Computation of ERPE for Fixed-design Logistic Regression Model 
## with normal prior for regression coefficient beta 
## with prior mean  given ##


####### Auxiliary Functions  --------------
## Logistic mass function ##
dlogistic=function(t,y)      
{
  z1=exp(t*y)/(1+exp(t))
  return (z1)
}

## Defining the function Q ##
Q=function(beta,X,y,alpha)  
{
  n=nrow(X)
  p=ncol(X)
  M=X%*%beta
  L=matrix(0,n,1,byrow=TRUE)
  for(i in 1:n)
  {
    L[i]=(1/alpha)*(dlogistic(M[i],y[i]))^(alpha)
              -(1/(1+alpha))*((dlogistic(M[i],0))^(1+alpha)
                              +(dlogistic(M[i],1))^(1+alpha))-(1/alpha)
  }
  L_0=matrix(0,n,1,byrow=TRUE)
  for(i in 1:n)
  {
    L_0[i]=log(dlogistic(M[i],y[i]))
  }
  if(alpha>0)
  {
    return(sum(L))
  }
  else
  {
    return(sum(L_0))
  }
}

## Normal Prior with mean true_beta ##
prior=function(beta,true_beta)   
{
  z=prod(dnorm(beta,true_beta,1))
  return (z)
}

## Kernel of the log of robust posterior density ##
robust_posterior=function(beta,X,y,alpha,true_beta)  
{
  z1=Q(beta,X,y,alpha)
  z2=log(prior(beta,true_beta))
  z=z1+z2
  return (z)
}

## Proposal Density ##
proposalfunction=function(beta)    
{
  p=length(beta)
  return (rnorm(p,mean=beta,sd=rep(1,p)))
}

## Metropolis Hastings ##
run_metropolis_MCMC <- function(startvalue,iterations,
                                X,y,alpha,true_beta)  
{
  p=length(startvalue)
  chain = array(dim = c(iterations+1,p))
  chain[1,] = startvalue
  for (i in 1:iterations){
    proposal = proposalfunction(chain[i,])
    
    probab = exp(robust_posterior(proposal,X,y,alpha,true_beta)
                 -robust_posterior(chain[i,],X,y,alpha,true_beta))
    if (runif(1) < probab){
      chain[i+1,] = proposal
    }else{
      chain[i+1,] = chain[i,]
    }
    
  }
  return(chain)
}


####### Computation of the ERPE --------------
data=read.table(file="data.txt")   ## Call the data-file  ##
y=data[,2]      ## Response ##
x=data[,-2]      ## Covariate ##
X=cbind(1,x)    ## Design matrix ##
n=nrow(X)       ## Number of observations ##
p=ncol(X)       ## Dimension of regression parameters excluding sigma ##

beta_0=c(-23,39.4,31.8)   ## Given value for the mean of the normal prior ## 

burn=25000
max_iter=50000
alpha = 0.5     ## tuning parameter alpha in alpha-likelihood
#set.seed(12345)
output=run_metropolis_MCMC(startvalue=beta_0,iterations = max_iter,
                           X=X,y=y,alpha=alpha,true_beta=beta_0)
output2=output^2
#Estimated ERPE, the means of the the R^\alpha-posterior distribution
est_mean=colMeans(output[(burn+1):max_iter+1,])
#Estimated variance of the R^\alpha-posterior distribution
est_var=colMeans(output2[(burn+1):max_iter+1,])-(est_mean)^2


\end{lstlisting}


\end{document}